\documentclass[13pt]{amsart}
\usepackage{amssymb,amscd, epic, eepic,  latexsym, newcommand}

\newtheorem{mainthm}{Main Theorem}

\newtheorem{theorem}{Theorem}[section]
\newtheorem{proposition}{Proposition}[section]
\newtheorem{definition}{Definition}[section]
\newtheorem{lemma}{Lemma}[section]
\newtheorem{corollary}{Corollary}[section]

\numberwithin{equation}{section}
\numberwithin{figure}{section}

\begin{document}
\title[Fixed point formulas and  loop group actions]{ Fixed point formulas and  loop group actions }
\author[S. Chang]{Sheldon X. Chang}
\thanks{Supported in part by a Sloan Fellowship}
\maketitle
 \centerline{\it Dedicated to the memory of Professor  Fredrick Almgren, Jr.} 

\section{Introduction}

In this paper we present  a new fixed point formula associated
with  loop group actions on infinite dimensional manifolds.
This     formula     provides information for  certain infinite dimensional
 situations similarly as    the well known Atiyah-Bott-Segal-Singer's formula does  in finite dimension. A generalization of the  latter to orbifolds     will be used as an intermediate step.

There exist extensive literature  on loop groups, loop
algebras and their representations. 
What set the present work and [C1]  apart  is the focus on 
representations of the central extensions of the loop groups, induced   from Hamiltonian actions;   getting  explicit formulas
which determine  the multiplicities of the irreducible highest weight components,
hence the structure of the induced represenations.

The results in  [C1] show
that the highest weight vector occurring in  an induced representation  are carried by
 a compact variety, provided the setup meets certain general  conditions.
In particular 
 geometric quantization is generalized  to the setting of loop group actions.
 The  results here 
 link the representations with  local data on  fixed
points. In that  direction, we also obtain a new multiplicity formula.

The current  paper can be read either as a sequel to [C1] or on its own. 

The original motivation was to understand in geometric setting a conjecture by Verlinde [V].  From  it the better known Verlinde formula was derived [V, MS].
As the project progressed, it became clear that  Verlinde's conjecture is
the tip of an iceberg.
Representations of the central extensions of  loop groups,  induced from Hamiltonian  actions 
on infinite dimensional  manifolds,  can be related 
to local geometry at the fixed point sets.

Application
in a forthcoming paper will include a direct proof
of the aforementioned conjecture.

More interestingly the we construct a class of $G$-orbifolds  called 
fusion product  which are geometric dual to the product in Verlinde
fusion algebra.

Much has been done  about Verlinde formula, we refer the readers to [F,Be] for
references.

\subsection{Some notations}

 Let $G$ be a connected  and simply connected compact simple  Lie
group, and $T$ be a maximal  torus.  Let $W$ be its Weyl group, $\ft_+$ be  the positive Weyl
chamber and $P, P_+$ the sets of weights and dominant weights respectively,   after fixing a set of simple roots. Let $D=\prod_{\alpha>0}(1-e^{-\alpha})$
denote the Weyl denominator,  
 $\theta$ be the highest root of $\fg$.  
On $\fg$, fix an invariant  bilinear form $(\cdot|\cdot)$, so that 
$(\theta^\vv|\theta^\vv)=2$ where $\theta^\vv$ is the  coroot corresponding to
$\theta$.
The bilinear  form induces a map $\nu :\ft\rightarrow \ft^*$ and 
a bilinear form on $\ft^*$.
The lattice in $\ft$  generated  by $\{W(\theta^\vv)\}$
is denoted by  $M$, and  $$M^*=\{t\in \ft\big |\, (t| n)\in \Z, \forall
n\in M\}.$$ 
Let $h^\vv$ be the dual Coxeter number, $\rho$ the half sum of positive roots.
Define the affine alcove and the set of dominant weights in there respectively
as  $$C=\ft_+\cap \{a|\, (a|\theta)\leq 1\},\quad P_+^k =\{\lambda\in P_+|\, (\lambda|\theta)\leq k\}= P_+\cap kC.$$
The face $C^\aff$ in $C$ defined by $(\cdot|\theta)=1$ is
special, in particular $$\partial (\cup_{w\in W} wC)=\cup_{w\in W} wC^\aff.$$ 
The lattice $\latk$ induces a finite subgroup of $T$, $\exp(2\pi \nu^{-1}\latk)$. The
subset
 \begin{equation}\label{life}
\explife
\end{equation}
 plays an important role in this paper.
An useful observation is that each element $\tau$ in that subset is 
a regular element of $G$.

\subsection{Assumptions on $X$}
Let $X$ be a Banach manifold on which $LG$ acts, $\omega$ be   an invariant
symplectic  form and $\mu$ be the associated  moment map at level $k\in \Z_+$, 
i.e., 
$$\mu:X\rightarrow \lfg^0\times\{k\}$$
where elements of $\lfg^0$ have one degree less of differentiability  than
those  in
$\lfg$.  Being a moment map,  $\mu$ is equivariant with respect to the $LG$-action on $X$
and the co-adjoint action on $\lfg^0\times\{k\}$. The co-adjoint action 
loses one degree of differentiability, hence   the  target of $\mu$ is  $\lfg^0\times\{k\}$. 
We assume in this paper the following unless stated otherwise:   

{\bf H1}: $\mu$ proper;

 {\bf H2}: $\mu(X)$ is transversal to $\ft\times\{k\}$.
More details are given in Section 2.

\subsection{The highest weight modules of the induced representations
}

In [C1], we studied holomorphic actions by a loop group $LG$,
  when $X$ is complex.
Let $N^+\subset LG^\C$  denote the  subgroup whose Lie algebra  is
$sum_{\alpha>0} (\lfg)^\C_\alpha$. It consists of  non-constant  boundary values of holomorphic maps from the unit disk to $G^\C$, together with constant maps
with values in the positve nilpotent group of $G^\C$.

 Suppose $L$ is an holomorphic $LG$-line bundle over $X$.
 We proved in [C1]  that there is  a  compact complex $T$-orbifold
$X_N$ and an orbifold $T$-line bundle $L_N$, so that $H^0(X_N,L_N)$
carries all the highest weight vectors of $H^0(X,L)$.
In other words,  as $T$-modules:
 $$H^0(X_N,L_N)\simeq H^0(X,L)^{N^+}.$$
The orbifold $X_N$ naturally can be viewed as the compactification of
the quotient of $X/N^+$. Therefore this  compact  model of $X$ carries  the same amount of
 information as $X$, in terms of understanding the induced  representation.

 To $X_N$ obviously one can apply the fixed point formula, a
certain generalization of
Atiyah-Bott-Segal-Singer results to orbifold, and get informations about the
 $T$-equivariant Riemann-Roch
$\RR(X_N,L_N)$, thus 
the structure of $H^0(X_N,L_N)$.
 But compactification involves adding
   certain locus,  and additional     fixed 
points sets on the locus.
To understand those new fixed points is not an easy
issue, even when dealing with  finite dimensional groups, e.g. 
 the  compactification of symmetric spaces.

However here  we  will find a solution to resolve this problem  here, utilizing the affine Weyl group $W^\aff$. As it  will be shown by
 examples, this solution is the best one can hope for.

\subsection{Description of the main result}

The main result in this paper does not require $X$ is complex, although
that situation motivates the construction of $X_N$ and later consideration.

 Let $Y=G\times_T X_N$,  the line bundle $L_N$ induces  $L_Y$   on $Y$.
If the original line bundle is of level $k \in \Z$ (or the moment map $\mu$ is
of level $k$), which means the central part of $\cLG$, $S^1$,
acts on $L$ with character $k$,
  the following function on $T$ will be  uniquely determined by its
restriction to $\explife$:
  $$ \sum_{w\in W}w\frac{1}{D}\RR(X_N, L_N)=\RR(Y,L_G).$$ 
The function $\RR(Y ,L_G)$ can be defined directly from $(X, L)$.
It is given by $$\sum_{a\in P_+^k}\RR(\calM_a,L_a)\chi_a$$
where $\calM_a$ is the reduced space      of $X$ at $a$, $L_a$ is 
the line bundle induced from $L$, $\chi_a$ is the $G$-character function
 of the highest weight representation defined by $a$. 

What  geometric data are needed  to  determine this function? 
Before answering that question, let's motivate the discussion by first
detailing the holomorphic case.
As mentioned earlier,  the quotient  $X/N^+$, after throwing away some bad orbits,
can be compactified by $X_N$. The compactification locus has its image  given by
 the boundary of the affine chamber $C$. The chamber is 
a simplex if $G$ is simple, and is a product of simplices if $G$ is semi-simple.  
For generic $X$, $X_N$ is an orbifold and strata in  the compactification  locus are  in 1-1
correspondence with the sub-faces of $C$.  Particularly interesting here are
the strata $\{X_Q\}$ whose images are on the affine wall $C^\aff$.
Each $X_Q$  has a corresponding  subvariety in $Y=G\times_T X_N$, $Y_Q$. 
The collection $\{Y_Q\}$ is    part of the compactifying strata 
in the $G$-space $Y$. 

Each  $\tau\in\explife$ is a regular element in $T$. Thus its fixed point
sets $\{V\}$ in $X$ has images under $\mu$  in $\ft$. Each component  $V$ will induce a subvariety $V_\Delta$ in $Y_\Delta=wX_N\subset Y$.
 And it also induces a subvariety $V_Q=V_\Delta\cap Y_Q$. 
The collections $\{V_\Delta\},\{V_Q\}$  will be used to determine
$\RR(Y)(\tau)$.  

We emphasize that in general  $\tau$ has lots more fixed point
sets than $\{V_\Delta\},   \{V_Q\}$ on 
 the compactification strata in $Y$.  Not all of the fixed points
in the compactification are in  the closure of the interior ones. 
So the important feature  of the 
main result is that   only the  closure of the interior $\tau$-fixed points
  $\{V_\Delta\}$,  and their intersection with strata in the
compactification $ \{V_Q\}$,  matter
in  determining $\RR(Y)(\tau)$. This  
feature manifests  the underlying  affine Weyl group symmetries, 
and  is not known to hold in finite dimension.
\begin{mainthm} \label{main}
At $\tau\in\explife$,
the following holds:
\begin{equation*}
\begin{split}
&\RR(Y)(\tau)\\
&=\sum_{a\in P_+^k}\RR(\calM_a,L_a)\chi_a(\tau)\\
&=\sum_{V_\Delta}\Big(\int_{V_\Delta}
\frac{\tTd( V_\Delta)\tCh(L_{V_\Delta})}
{ \det_{\no(V_\Delta, Y)}(1-t^{-1}e^{- \Omega})}\\
 &\quad\quad + \sum_{V_Q\subset V_\Delta}\frac{1}{|W_Q^\aff||I_{V_Q}|} \sum_{t\in \tau
I_{V_Q} } \int_{V_Q}\frac{\tTd(
V_Q)\tCh\left(L_{V_Q}
\oplus
\Lambda^{\max}\no(V_Q,V_\Delta )|_{V_Q} \right)}{ \det_{\no(V_Q, Y_Q)}(1-t^{-1}e^{-\Omega})}\Big)(\tau)
\end{split}
\end{equation*}
where  $W_Q^\aff\subset W^\aff$ is the  subgroup  preserving  $Q$, $I_{V_Q}$
is the isotropy group associated with $V_Q$ and $\tau I_{V_Q}$ is the set
of liftings of $\tau$.
Furthermore, each integral  above can be localized to the $T$-fixed point sets  $F$ in $V_\Delta, V_Q$ respectively to yield:
\begin{equation}
\RR(Y)(\tau) =\sum_{\{F|\mu(F)\in W(C^\inte)\}} \FC(F)(\tau) +\RT(\tau)
\end{equation}
where 
\begin{equation*}\begin{split}
&\FC_F(\tau)=\int_F\frac{\tTd( F)\tCh(L_F)}
{ \det_{\no(F,Y)}(1-t^{-1}e^{-\Omega})}(\tau);\\
&\RT(\tau)=
\sum_{\{F|\phi(F)\in W(C^\aff)\}}\frac{1}{|W_\phi^\aff||I_{F}|}  \sum_{t\in \tau I_F}\int_{F}\frac{\tTd(
F)\tCh\left(L_F
\oplus
\Lambda^{\max}\no(Y_Q,Y )_F \right)}{ \det_{\no(F,Y_Q)}(1-t^{-1}e^{-\Omega})}(\tau).
\end{split} \end{equation*}
where $C^\inte$ is the interior of $C$, $W_\phi^\aff$ is the subgroup 
of $W^\aff$ preserving $\phi(F)$, and $I_{F}$ the isotropy group
of  $F$.
\end{mainthm}

{\it Remark:} 1). Similar to an earlier comment, the interesting feature in the second expression above is that  only those $F$ on the intersection $V_\Delta\cap Y_Q$ matters.  Other
$T$-fixed points on the compactification locus $Y_Q$ do exist and there are lots of
them, but they do not contribute to $\RR(Y)(\tau)$ as it will be shown. 

2).  
 The presence of $I_F, \tau I_F, I_{V_Q}, \tau I_{V_Q}$  in fixed
point formula is a common feature in  orbifold setting, this  
has been known for a while.

Consequences on $\cLG$-modules are given in Section 10.

\subsection{Riemann-Roch of the reduced spaces}
  The previous  result provides a way of  computing 
the Riemann-Roch numbers
 of the  reduced space $\calM_a=\mu^{-1}(a)/(LG)_a$ via
certain fixed point sets. 
\begin{corollary}\label{multiplicity}
\begin{equation*}\begin{split}
&\RR(\calM_a,L_a)=\\
&\frac{(-1)^l}{\big|\frac{M^*}{(k+h^\vv)M}\big|}\sum_{\tau\in \explife}
 \chi_{\bar{a}}(\tau)D^2(\tau)\Big(\sum_{\{F|\mu(F)\in kW(C^\inte)\}}  \FC(F)(\tau) +\RT(\tau)\Big)
\end{split}\end{equation*}
where $\bar{a}=w_L(-a)$ with  $w_L$ being the longest element in $W$ (or the
highest weight  in the contragredient representation  to the one defined by $a$).
\end{corollary}

{\it Remark:}
The above extends Verlinde's formula for moduli space of flat connections
over Riemann-surfaces.

If $X$ is holomorphic and  
 $H^i(\calM_a,L_a)=0$ for $ i>0$, then the above formula gives the
multiplicities of the irreducible component with highest weight $(a,k)$ for 
the $\cLG$ representation on $H^0(X,L)$. 

In finite dimension, the Riemann-Roch number of the reduced space can also
be expressed in terms of the fixed points. The expression is obtained only
recently as  an application of
Atiyah-Bott-Singer-Segal's formula and the recent in [M].

\subsection{What's  in the proof?}
The theorem is proved based on two main ingredients, both utilize the
Weyl and affine Weyl group symmetries.

After applying directly fixed point formula for orbifolds to $Y$, 
one ends up with many terms of contributions from the fixed point
sets on the compactifying locus. There are two kinds of them, 
those of the first kind have their  images on $W(\partial C\setminus C^\aff)$,
we prove using Weyl group symmetry that their contribution  amounts to 0.
The second kind are those with images on $W(C^\aff)$.  Their contribution
to the equivariant Riemann-Roch also amounts to 0, provided we restrict
them as  functions  on $T$ to the subset $\explife$, and  each
 term has no pole on the subset. This is proved using affine Weyl group
symmetries.

 We shall refer to the above  phenomenon as cancellation, it is based
on  the fundamental formula of Section 6 and several identies proved in
Section 7. Also it requires detailed analysis
of the fixed points on the compactification locus.
 
The second ingredient starts with  the  surgery formula in Section 12, it enables us
to deal with the second kind of fixed points when  poles are present.
 Their contribution is encoded in the function $\RT: \explife \rightarrow \C$. The exact expression of $\RT$ is determined using transformation rule
of the affine Weyl group acting on the Riemann-Roch integrand, together
with symplectic cuts. The calculation is rather elaborate.

The construction of the space $X_N$, was first done in 1993 prior the symplectic
cuts. It shares certain  similarity with   the symplectic cuts,
except the cuts are made along the degenerated parts of the two
form $\omega|\mu^{-1}(\ft)$. 
The resulting space is only  symplectic outside  the inverse images
of $\partial C$.  
 The surgery formula for this kind of cuts is quite
different from that using symplectic cuts, as shown by Prop.~12.1.

The cancellation mentioned earlier has a consequence called    twin pair
construction. It says that for a generic  compact symplectic manifold with a
Hamiltonian  $G$-action, there is a different $G$-orbifold, with identical
equivariant Riemann-Roch.

\subsection{ What is ahead?}
As mentioned earlier, applications will be given in a forthcoming paper.

In a separate paper, the result here will be improved  so that all
symplectic $LG$-manifolds, with compact quotient $X/LG$,  will be covered. 
The present results assume the generic condition that $$\mu(X)\subset \lfg\times \{k\}\subset  \clg$$
is transversal to  $\ft$.

\subsection{   Organization of 
 the  paper}

In section 2, we discuss the construction of $X_N$. 
The paper [C1] emphasizes on the holomorphic aspects of $X_N$ while
here the construction builds around the symplectic structure.
 Several new results are presented here, including the existence
of  $LG$-invariant almost complex structures on 
a class of $LG$-manifolds. The existence of $T$-invariant almost complex structure
on $X_N$ is proved as well, which is not trivial considering that
the symplectic form on $X_N$ is degenerate.  

Section 3 contains the description of the fixed point sets, and their
stratification which is a must since $X_N$ is an orbifold. And the fixed point
formulas on orbifolds rely  not only on the fixed points but their
stratification as well. Section 4 includes the computations of weights
of the induced $T$-action on the normal bundles to the fixed point sets
and their lower strata, while Section 5 computes the curvatures of various
components of the normal bundles.

The root of the cancellation  
 is the the fundamental formula presented in  Section 6. 
 I proved this formula  in 1994.
The original proof was based on comparison of two different compactifications
of $G^\C$. Both have the same Riemann-Roch but have different  fixed points.
Hence one yields an identity, then by induction on the 
rank of the group, one proves this formula.
The result in Section 12 generalizes this `twin pair' construction.
 The present proof is 
simplified   by applying an identity which can be found in  [M].

Another important component in proving the cancellation is described in
Section 7.   Section 8 addresses a complication which occurs when  
$G\neq SU(n)$. For those groups, the affine alcove may not  be a  simple simplex,
with respect to the weight lattice.
 In  constructing  $X_N$, this fact introduces addition orbifold
singularities. That section describes the extra components of the isotropy
groups  associated with the orbifold singularities.

Section 9 provides a way of computing the  push-forward of certain  cohomology classes on 
a fibration whose fibers are homogeneous spaces like $K/T$. Several
integration formulas there  can be viewed as localization formulas for
families.  

In section 10, we discuss the relations between $T$-spaces and $G$-spaces,
characters of $T$-modules and $G$-modules. Several consequences of the
main theorem are proved.

Section 11 proves the main cancellation  which has a couple of consequences
discussed in Section 12, including the `twin pair'.
The proof  of the main result is completed in the last section after 
we find an expression for the remainder term.

\subsection{ Acknowledgment}
I am indebted to Prof.~Kac   who answered   my
questions with great patience and insight,
  to Prof.~Liu for many enjoyable and informative  conversations.
E. Meinrenken
and C. Woodward asked  an interesting question about the first version
of the results, I thank them for discussion.
Words of  encouragement  from Professors Cheng and  Stroock 
  have been valuable.
Last but not the least, I am grateful to my family for their support and
patience.

\section{Basic properties of the variety $X_N$}
 Many properties presented  here were proved in [C1], with the
exception of  the existence of a $LG$-invariant almost complex
structures on $X$,  and a $T$-invariant one on $ X_N$. In [C1], $X$ is assumed to be holomorphic, therefore it is not necessary.
We list the important properties for $\lfg, \X, X_N$ here.

\subsection{Basics of affine Lie algebra}
Let $\fg$ be the  Lie algebra of $G$  and $T$ be a maximal torus with $\ft$ as 
its Lie algebra.
The discussion works     without much modification for
semi-simple Lie groups.
Let    $LG$ be the loop group associated with $G$,
and $\lfg=\Lie LG$.

For functions on $S^1$, we use the $t\in[0,1]$ to parameterize them.
 The Fourier series components are $\{e^{\lp n t}\}$.
On $\lfg$, there is  the following well defined form in terms of the
invariant non-degenerate  form $(\cdot|\cdot)$ on $\fg^\C$:
$$(a|b)=\int_0^1(a(t)|b(t)) dt\in \R.$$
It induces a  symplectic form on $\lfg/\fg$:
$$B(a,b)=(a'|b)=\int_0^1(a'(t)|b(t))dt$$
where $a'=da/dt, t\in[0,1]$.
The form is degenerate when restricted to the constant $\fg$.
 We choose the form on $\fg$ with
the condition that
$ (\theta^\vv|\theta^\vv)=2$ or equivalently $(\theta|\theta)=2$, where
$\theta^\vv$ is the coroot corresponding to the highest root $\theta$
of $\fg$.
As pointed out in [PS, p. 46], the associated form $B$  defines  the smallest integral class on
$LG$.

The affine Lie algebra based on $\fg$, $\fg^\aff$, is defined as 
$$\fg^\aff=\lfg\oplus \R d\oplus \R K,$$ 
where $K$ is the central element and $d=d/dt$ is the differentiation.
The Lie bracket is 
$$[\xi+\lda d+cK,\eta+\ldaa d+c' K]
=[\xi,\eta]+\lda d\eta-\ldaa d\xi+B(\xi,\eta)K.$$
The central extension of $\lfg$, $\clg$, is given by
$\fg^\aff=\lfg\oplus \R K$.
  It is the Lie algebra of $\cLG$ which is a circle bundle over $LG$,
whose existence of $\cLG$ is proved in [PS].

The Lie algebra dual $\affgd$ is given by 
$$ \affgd=\lfg^*\oplus \R \delta\oplus \R \Lambda_0.$$ 
The bilinear form $(\cdot|\cdot)$ extends to $\affg$, $\affgd$, so that
it is invariant.  Its restriction to the
2-dim subspace $\R d\oplus \R K$,   is of the form 
$\begin{pmatrix}0 &1\\1&0\end{pmatrix}$.
See [K, Ch. 6] for more details. Continue to use $\nu$ for the map
$\affg\rightarrow\affgd$ defined by $(\cdot|\cdot)$.
 The bilinear form gives $$\nu(\lfg\oplus \R d)= \lfg^*\oplus
\R\Lambda_0.$$
The simple roots of $\affg$ consists of  simple roots of $\fg$ together with
 $\alpha_0:=\delta-\theta.$
Suppose that $\{E_i',F_i'\}$ are the Chevalley basis of $\fg$, let 
$E_0',F_0'\in \fg^\C_\theta$, so that 
$$[E_0',F_0']=-\theta^\vv,$$
where $-E_0$ is the Chevalley involution of $F_0'$, 
define $$e_0=z\otimes E_0',\quad  f_0=z^{-1}\otimes F_0'.$$
For $su(2)$, the pair is 
$$\begin{pmatrix}  0& z\\ 0&0 \end{pmatrix},\quad \begin{pmatrix}  0& 0\\ z^{-1}&0
\end{pmatrix}. $$
 The pair $e_0,f_0$ together with $ e_i=E_i',f_i=F_i'$, $ 1\le i\le l$ generate
$\lfg^\C$.

The   positive affine roots $\gamma\in \Delta_+(\fg^\aff)$ are
$$\{n\delta\pm\alpha| n\ge 1,\alpha\in \Delta_+(\fg)\}\cup
\Delta_+(\fg)\cup\{n\delta|n\ge 1 \}$$
and accordingly the basis of the root spaces  are $E_\gamma=z^n E'_{\pm\alpha},  E'_\alpha$ or $z^n h_\alpha$,
where $\{E'_\alpha\}$ is the root space basis of $\fg$.
The basis of the compact form are
$$\{x_\gamma=e_{\gamma}-f_{\gamma},  y_\gamma=i(e_{\gamma
}+f_{\gamma})\}$$
where $\gamma>0$.
The standard   complex structure $J^\lfg$ on $\lfg/\ft$ which inherits from
the map  $n^+\rightarrow \lfg/\ft\simeq \lfg^\C/(n^-+\ft^\C)$  has this
description:
 \begin{equation}\label{complex1}
J^\lfg(x_\gamma)=y_\gamma,\quad J^\lfg(y_\gamma)=-x_\gamma.
\end{equation}
By direct computation, for $h\in \ft$, we obtain
\begin{equation}
\begin{split}
 &x_\gamma'=2\pi ny_\gamma, \quad y_\gamma'=-2\pi nx_\gamma;\\
&[h,x_\gamma]=(-1)^{\sign(\alpha)} \alpha(h)/i y_\gamma=-(-1)^{\sign(\alpha)}
i\alpha(h) y_\gamma;\\
&[h,y_\gamma]=-(-1)^{\sign(\alpha)} \alpha(h)/i x_\gamma=(-1)^{\sign(\alpha)}
i\alpha(h) x_\gamma. 
\end{split}
\end{equation}
Hence for $\xi= \sum a x_\gamma+by_\gamma$,
\begin{equation}\label{kahlerlg}
\begin{split}
&(J\xi)'+[h,J\xi])=-\sum (2\pi n-(-1)^{\sign(\alpha)}i\alpha(h))(a x_\gamma
+by_\gamma);\\
&\left((J\xi)'+[h,J\xi])|\xi\right)=-\sum
(2\pi n+(-1)^{\sign(\alpha)}\alpha(h)/i)(|a|^2|x_\gamma|^2+|b|^2|y_\gamma|^2).
\end{split}
\end{equation}
Now suppose that $\frac{1}{\lp}h $ is in the affine alcove $ C$, then $$n+(-1)^{\sign(\alpha)}\frac{1}{\lp}\alpha(h)=<n\delta\pm\alpha,
\frac{d+h}{\lp}> \ge 0,$$ 
because all the positive roots are positive on $C+\frac{d}{\lp}$, (our
definition of $d$ differs from that in [K], hence the extra constant).
Therefore $\left((J\xi)'+[h,J\xi])|\xi\right)\le 0$.
Let $\Omega(\xi,\eta)= (\xi'+[h,\xi]|\eta)$, earlier calculation shows
that 
\begin{lemma} The form $\Omega(\cdot,J\cdot )$ is semi-positive definite.
\end{lemma}

\subsection{Adjoint and co-adjoint action}
Let  $\xi,\eta \in\lfg$, then  $[\xi+ad, \eta]=[\xi,\eta]+ad \eta$ which
integrates to $$g^{-1}(\xi+ad)g=\Ad_{g^{-1}}\xi+ag^{-1}dg+ad\in \affg.$$
The linear map $\nu:\affg\rightarrow\affgd$ defined by $(\cdot|\cdot)$
has  the effect that $\nu(\xi+ad)=\nu(\xi)+a\Lambda_0$, therefore
the adjoint action by $g$ is  given by
\begin{definition}\label{ad-action}
The adjoint action on $h+a\Lambda_0$ is 
$$ gh g^{-1}+ag \frac{d}{dt}g^{-1}+ad.$$
\end{definition}
The level of $h+a\Lambda_0$ is $a$, thus  the adjoint
action preserves the level.
 In case of level 1, the projection to the $\lfg$ part is exactly the gauge
transformation.

\subsection{Topology on $LG$}
Since the group $LG$ we are dealing with is the mapping space
  $Map(S^1,G)$,
there is the question as to  
   which norm is used for completion. The norms defined by
$B(\cdot,J^\lfg\cdot), \Omega(\cdot,J^\lfg\cdot)$ are  not strong enough.  The weakest norm which is
sensible geometrically, is
the  $H^1$-norm.
We can use other Hilbert metrics or Banach metrics as long as they are
stronger than $H^1$-norm.

{\bf Convention}:
The convention here is to
let $\lfg=H(S^1,\fg)$ be the space of maps completed under the norm of choice,
and $LG$ be the corresponding group. Set
$\lfg^0$ to be the set with one degree less of derivative, i.e.,  $$\lfg^0=H'(S^1,\fg)=\{h'|h\in\lfg\}\oplus \fg.$$
 It  is another   completion of $C^\infty(S^1,\fg)$,
 so that the map $$h\in \lfg\mapsto h'\in\lfg^0$$ is Fredholm and  bounded. 

 We will  simply use $\lfg$ for $\lfg^0$, just keep on mind that 
the target of $\mu$ consists of elements with  one degree less of 
differentiablity. 

\subsection{Loop group actions and the assumptions {\bf H1, H2}.}

Let $X$ be a Banach manifold with a differentiable action by
$LG$,
   $$\mu:X\rightarrow \lfg\oplus\R \simeq \lfg^*\oplus\R\Lambda_0 $$ be a moment map associated with a symplectic 2-form $\omega$ on $X$.
The isomorphism $\simeq $ is defined by the restriction of  $\nu$.
The moment map  is equivariant with respect to the $LG$-action on $X$
and the adjoint action on $\lfg^*\oplus\R\Lambda_0$,

{\it Remark:}   The Banach norm on
$TX$ does not have to be invariant; and the positive definite form
$\omega(\cdot,J\cdot)$, in general, defines a topology weaker than  the Banach norm, for any compatible almost complex structure on the tangent space.

\begin{definition}
 $\mu$ is of level $k\in \Z_+$
if $\mu(X)\subset \lfg\times \{k\}$.
\end{definition} 

{\it Remark:}  The topology on $\lfg\oplus\R$ as described makes
the co-adjoint action a bounded map.

The following assumptions will be made:

{\bf H1}:  $\mu$ is  proper with aforementioned topology.

{\bf H2}: $\mu$ is transversal to $\ft\times \{k\}$ in $\lfg\times \{k\}$.

The first one is essential and is equivalent to the compactness of $X/LG$,  and the second is  technical.

Assuming {\bf H1} and {\bf H2}, then $X_\ft=\mu^{-1}(\ft\times\{k\})$
 is a finite dimensional
submanifold. It is not symplectic, $\omega|X_\ft$ has serious degeneracy.
And it may not even be orientable.
Whenever the stabilizer of $\mu(p)/lp$ in $\lfg$,  $(\lfg)_{\mu(p)}$, has a  semi-simple part, 
$\omega|_{T_pX_\ft}$ is null on $(\lfg)_{\mu(p)}/T$.

\subsection{Toric bundle $\X$ over $LG/T$.}
It is a fundamental fact that the affine Weyl group $W^\aff$ acts on
$\ft\times \{k\}$, the quotient domain of the action is  
given by a simplex $k(C,1)$ where $C$ is the  affine alcove of $\fg$.
One can also consider the action by $W^\aff$ on the dual space
$\ft^*\times\{k\}$, the quotient  domain  is $kC^*$ with $C^*$ given by  
$$C^*=\{\lambda| \, (\alpha|\lambda )\geq 0, (\theta|\lambda)\leq 1\}=\{\lambda|<\alpha^\vv,\lambda>\geq 0, <\theta^\vv,\lambda>\leq 1\}$$
where $\alpha^\vv,\theta^\vv$ are the coroots. 
The above descriptions of $C,C^*$ are the duals of each other, through the map
$\nu:\ft\rightarrow \ft^*$.

When there is no confusion, we will not distinguish between $C,C^*$.

The simplex $C^*$ is not simple with respect to the weight lattice,
the edges do not form a base of the weight lattice of $\ft$. In fact
the edges are given by $$\Lambda_i/a_i^\vv,\quad i=1,...,l$$ where 
$l$ is the rank of $\fg$, $\{\Lambda_i\}$ is the set of fundamental 
weights of $\fg$ and $\{a_i^\vv\}$ are  the labels in the dual Dynkin diagram.
More on this can be found in [C1].

From the theory of toric varieties, we know that there is an orbifold toric 
variety $X_\fg$ and an orbifold line bundle $L_N$  associated with $C$. 
The pair is the quotient of $\C P^l,H$ by a finite group, where $H$ is the
hyperplane line bundle.
The details  are in [C1].

Given $X_\fg$, we can associate with it a  toric bundle over $LG/T$:
$$\X=LG\times_TX_\fg,\quad (gt,z)\simeq (g,tz).$$
The quotient is well defined since $T$ is compact and the action is free
of fixed points. The group $LG$ has the same topology as described earlier.

The construction of $\X$ given here  dated back to 1993.

There are many nice characteristics about $\X$, we list a few needed 
later.

On $\X$,  there is an action by $LG$ and $T$ respectively, the two actions
commute.
The actions by $LG,T$ satisfy the Hamiltonian conditions with moment maps given by
$$\muX:\X\rightarrow \lfg^0\times\{1\},\quad  \phi:\X\rightarrow \ft.$$ 
Let $\tphi=(\phi,1) \in \ft\times\{1\}$.
Then $$\muX([g,z])=\tAd_g(\tphi(z))$$ 
that the above is independent of the choice of $(g,z)$ in $[g,z]$ is evident. 
Let 
$T_{[I,z]}=\lfg/\ft\oplus T_zX_\fg$. The 2-form on $\X$ can be described as
 $$\omega_\X|_{[I,z]}( (\xi,a),(\eta,b))=(\xi'+[\phi(z),\xi]|\eta)+
\omega_{X_\fg}(a,b)\quad a,b\in T_zX_\fg \quad  \xi,\eta\in \lfg/\ft.  $$
 The form  is degenerate whenever $\phi(z)$ hits the boundary of the affine alcove
$C$. The null space is generated by $(LG)_\tphi^\sss/T$ where $(LG)_\tphi$ 
is the stabilizer of $(\phi,1)$.
  The complex structure $J$ 
is defined as:
$J|_{\lfg/\ft}=J^\lfg$ while on $T_zX_\fg$, it is given by
the  original one on $T_zX_\fg$.
The choice of this complex structure is due to the following reasons:

1). The action by $t\in T$ on the left is $$t(g,z)=(tg,z)\simeq (tgt^{-1},tz).$$

2). We want highest weight modules, rather than lowerst weight ones.

\subsection{The variety $X_N$}
Reverse the complex structure on $X$, so that $-\mu$ is the moment map,
and $-\omega$ is compactible with the complex structure $-J$.

The variety $X_N$ which is important to our study can now be described as
follows:
$$X_N=(\Psi^{-1}(0)/LG=\{(p,q)|\mu(p)=k\muX(q)\}/LG,$$
where  $\Psi=-\mu+k\muX:X\times \X\rightarrow \lfg^0$ is the moment map
associated with the  diagonal action by $LG$ on the product space, with
the 2-form $-\omega+\omega_\X$.

 The sign here is chosen so that no inversion of the complex structure 
on $\X$ is necessary, this way we still get the highest vectors in the end.

 Notice the level of $\Psi$ is 0.
Set $$Y_C=\mu^{-1}(kC\times \{k\}),$$
by {\bf H1,H2} it is a compact  manifold  with boundary.

The following shows that the toric variety can be  used to close the gash
which is the boundary of $Y_C$.
\begin{proposition} The space  $X_N$ is $$\{(p,q)\in Y_c\times
X_\fg|\, \mu(p)=k\tphi(q)\}/T.$$ 
\end{proposition}
The proof is  simple, since each pair $(p,q)$ with $\Psi(p,q)=0$ can
be conjugated to $(p',q')\in Y_c\times_T
X_\fg $ with $\mu(p')=k\muX(q')\in k(C,1)$; the pair is unique up to $T$.

From this description, considering the assumptions {\bf H1, H2}, it is
clear that $X_N$ is a compact orbifold.
 The claim that  $X_N$ is holomorphic whenever
$X$ is complex is  of more subtle nature, it is proved in [C1].

Because $T,LG$ commute on $\X$, the action by $T$ descends from the
product $X\times\X$ to $X_N$. So does the moment map $\phi$.
The form $-\omega+\omega_\X$ when restricted to $\Psi^{-1}(0)$ is invariant
under the $LG$-action, thus it descends down to a form on $X_N$. Denote it by
$\omega_N$. The pair $\phi$, $\omega_N$ satisfy the conditions for Hamiltonian 
action, though $\omega_N$ is degenerate.

When  $X$ is not complex, we will see there  is an $T$-invariant  almost complex structure
$J$ on $X_N$, such that $\omega_N(\cdot,J\cdot)$ is semi-positive definite.
If $\omega_N$ on  $X_N$ is symplectic, the existence of such a $J$ is well
known. For degenerate $\omega_N$, the existence is not automatic.

\subsection{The existence of an $LG$-invariant $J$ on $X$.}
The existence result only needs assuming {\bf H1}.

{\it Step 1: Existence of  a positive bilinear $LG$-invariant form on $X$. }

Let $\{v\}$ be the vertices of $ k(C,1)$. There are $l+1$ of them where 
$l$ is the rank of $\fg$.
For each $v$, let  $C^v$ be $k(C,1)$ after removing the face opposite to $v$.  
Let $W_v$ be the Weyl subgroup in $W^\aff$ generated by reflections with
respect to walls passing $v$. It is well known that 
$W_v$ is finite and is  the Weyl group of $(LG)_v$ which  stabilizes $v$ under the co-adjoint
action. Now set $$O_v=\cup_{w\in W_v}wC^v,$$
which is an open set. It is the star-shaped region with center $v$ if we view
the images of $C$ under $W^\aff$ as a triangulation of $\ft\times\{k\}$.
E.g. when $\fg=su(2)$, $C=[0,1]$ while $O_0=(-1,1), O_1=(0,2)$.

Clearly $\cup_vO_v$ is an open cover of $k(C,1)$, and  there exists
a partition of unity $\{\psi_v\}$  subordinate to the covering, and
 $\psi_v$ is invariant under $W_v$. 
Define 
 $$ S_v=\mu^{-1}\big(\tAd_{(LG)_v}(O_v)\big);\quad \calS_a=\mu^{-1}
(\tAd_{LG}O_a).$$ 
Using the conditions on $\mu$, we have:
$$\calS_a=LG(S_a)=LG\times_{(LG)_a} (LG)_v(S_v),$$
and $\{\calS_a\}$ is a $LG$-invariant open cover of $X$.
 Let $g_a$ be an $(LG)_a$-invariant Riemannian
metric on $S_a$. It exists, since $(LG)_a$ is compact and $S_a$ is of
finite dimension.

Next we define a  positive bilinear $LG$-invariant form $g_a^A$  on $\calS_a$.
Clearly 
\begin{equation}\begin{split}
T_p \calS_a=\lfg/(\lfg)_a\oplus
T_p S_a,
\end{split}
\end{equation}
on the second factor there is already a Riemannian metric $g_a$.
 On the first one, 
the choice for $J$ is 
$J^\lfg$.
 With this choice,
we have the positivity of 
$\Omega(\xi,J\xi)$ as shown earlier.

The form $(\cdot|\cdot)$ is the restriction  of the bilinear invariant 
form on $\fg^\aff=\lfg+\R \delta+\R d$, (see [K, Ch. 7]).
And $[J\eta,x]=J[\eta,x]$ which can be verified using the property that
if $$[E_\alpha,E_\beta= c_{\alpha,\beta} E_{\alpha+\beta},$$ then 
$$[E_{-\alpha},E_{-\beta}= -c_{\alpha,\beta} E_{-\alpha-\beta}.$$
  Thus $g_a^A$  is invariant under the 
conjugation on $\lfg/(\lfg)_a$ by $(LG)_\mu$.

\begin{lemma}\label{symdec} The decomposition $$T_p\calS_a=
\lfg/(\lfg)_a\oplus
T_p S_a$$ is orthogonal with respect to $g^A_a, \omega_\X$.
\end{lemma} 
{\it Pf:} It is orthogonal with respect to  $g^A_a$ by construction.
To check that for $\omega_\X$, let $u\in T_pS_a$ and $\xi\in \lfg/(\lfg)_a$, 
by the definition of moment map,
$$\omega_\X(u,\xi(p))=(D_u\mu|\xi)=0$$
where the last equality holds because $D_u\mu\in (\lfg)_a\perp \lfg/(\lfg)_a$
under $(\cdot|\cdot)$. 
 QED

We can use the group action by $LG$  to extend $g_a^A$ 
to a 
$LG$-invariant form denoted by the same  on $\calS_a$. This is made possible due to 

1).  the invariance
of $g_a$ under $(LG)_a$ on $S_a$;

2). The complement of $T_pS_a$ in $T_p\calS_a$ is $\lfg/(\lfg)_a$ whose 
positive bilinear form as in Eq. (\ref{posit}) is invariant under the
conjugation.

Now  for each $p$ in $\calS_a$, the image under $\mu$ of the orbit $LG(p)$ 
meets $ O_a$ in a $W_a$-orbit on which $\psi_a$ is constant, since $\psi_a$ 
is $W_a$-invariant. Therefore, we can extend $\psi_a$ to $\calS_a$. The
extension will still be denoted by $\psi_a$. Obviously, $\{\psi_a\}$ 
forms a partition  of unity on $X$ subordinate to $\{\calS_a\}$.

Now clearly  $g=\sum_i\psi_a g^A_a$ is an invariant positive bilinear form. 

{\it Remark:} We do not call this form a Riemannian metric in order
to avoid confusion, since the topology defined by $g$ is weaker than
that on $X$. 

{\it Step 2:  Existence of an $LG$-invariant almost complex structure}.

On finite dimensional space $M$, given an invariant symplectic form $\omega$ and 
a positive definite form $g$, there is an uniquely well defined almost complex 
structure. 
The exact construction of $j$ in terms of $g,\omega$ is like this:
the non-degenerate 2-form can be viewed as a linear map 
$$\omega:TM\rightarrow T^*M,$$
the form $g$ defines a map from $T^*M\rightarrow TM$ which is the inverse
of the map $\nu:TM\rightarrow T^*M$ defined by the metric.

Denote  $\nu^{-1} \cdot \omega$  by $f:TM\rightarrow TM$.
Then it is straightforward that $f^*=-f$, the dual is taken with respect to
$g$. Therefore, $-f f=(f^*f) $ is positive definite and 
has an uniquely defined positive definite square root, $(-f^2)^{1/2}$ which
is a linear map. Clearly $f,(-f^2)^{1/2}$ commute, so the following is the
desired $j$, 
$$j=f(-f^2)^{-1/2}.$$
This almost complex structure is clearly invariant.

{\it Two remarks:} 1). The positive form $\omega(\cdot,j\cdot)$ may not
coincide with $g$; 2). If $g$ is given by  $\omega(\cdot,J\cdot)$, then
$j=J$. 

Both of the above are easy to verify.

In the infinite dimensional situation,
we refrain from defining $J$ on $TX$ by directly using 
the above
 because the definition of $f^*f$ may cause problem.

In the particular situation we are facing, however, the map $f^*f$ will be shown
to be $I$ except on a finite dimensional space.
Let us examine the matter more closely.

Suppose $\mu=\mu(p)\in k(C,1)$ is covered by  $\{O_b\}$ but no other $O_a$. 
By definition of $O_a$, it is clear that the stabilizer of $\mu$ in $W^\aff$,
$W_\mu$, is 
$$W_\mu=\cap _{\{b| O_b\ni \mu\}} W_b.$$
When $\mu$ is in the interior, $W_\mu=\{I\}$ and all of $\{O_a\}$ cover it;
when $\mu=v$ is a vertex, only $O_v$ covers it.
Also $\mu$ is not in the support of $\psi_a$ if $\mu\notin O_a$.
Therefore at $p$, the positive form $$ g=\sum_{\{b|O_b\ni \mu\}}\psi_b g^A_b.$$

On $T_p\calS_b=\lfg/(\lfg)_b\oplus T_pS_b$, $g^A_b$ agrees with the form
\begin{equation}\label{lieform}
([\cdot,J\cdot]|\mu(p)) 
\end{equation}
on the first factor, for each $b$ with $p\in O_b$. Therefore
inside that tangent subspace, $g$ agrees with the expression in Eq.~\ref{lieform}
on $$\cap _{p\in O_b}\lfg/ (\lfg)_b.$$

By the previous lemma, the finite dimensional subspace $E_p$  generated by 
$$ \{T_pS_b\big|\mu( p)\in O_b\}$$
 is orthogonal to $\cap _{\{b| O_b\ni \mu(p)\}}\lfg/ (\lfg)_b$
 with respect to $\omega,g$. 
Hence, $\omega$ and $g$ are in diagonal  form in the decomposition
$$T_pX= \cap_{\{b| O_b\ni \mu(p)\}}\lfg/ (\lfg)_b \oplus E_p.$$ 
Thus the map $f=\nu^{-1}\omega$ is of diagonal form.
On the first factor, again by the previous lemma and the definition of $g_b^A$, 
$$g=([\cdot,J\cdot]|\mu),$$ therefore over that
subspace $$f=J=J^\lfg;\quad
 f(-f^2)^{-1/2}=f=J^\lfg.$$ 
On the finite dimensional subspace $E_p$, $J$ is uniquely defined by
the restriction of $f$. Hence it is invariant by $(LG)_\mu$.
  We can easily extend $J$ to the whole $X$ through invariance.
Thus we have just proved the following
\begin{proposition}
There is a $LG$-invariant almost complex structure on $X$ such that
$\omega(\cdot,J\cdot)=-(D_{J\cdot}\mu|\cdot)$ is positive definite.
\end{proposition}
{\it Remark:} If $X$ is complex, the above equality $J=-J^\lfg$ on $\cap _{\{b|
O_b\ni \mu(p)\}}\lfg/ (\lfg)_b$ is  not true in general. 

\subsection{The almost complex structure $J$ on $X_N$}

We will prove the existence of an almost complex structure on $X_N$ by
showing that there exists Hodge type decomposition  for
$T_{(p,q)}X\times\X$ where $\mu(p)=k\muX(q)$.
Such a decomposition is used in the finite dimensional situation to show the
existence of almost complex structure (or holomorphic structure if the
original manifold is) on the reduced space.

Suppose $\mu(p)=k\muX(q)\in k(C,1)$. 
It will be shown in the next section, as a consequence of {\bf H2}, that $\ft_p\cap \ft_q=0$ which 
 implies  $(\lfg)_p\cap(\lfg)_q=0$, since $(\lfg)_q=\ft_q$.

The following is easier to show than in the case when $X$ is holomorphic, for
the almost complex structure $J^X$ constructed earlier has very special
property.

From now on,
 let $X$ be equipped with the opposite of its usual complex structure.
And continue to use the same one on $\X$ so that $-\omega(\cdot,J\cdot)+
\omega_\X(\cdot,J\cdot)$ is
non-negative definite.

\begin{proposition}
The action by $\lfg$ on the product $T_{(p,q)}X\times\X, \Psi(p,q)=-\mu(p)+k\muX(q)=0$
has no kernel.
There is a $LG$-invariant (acting diagonally)
 decomposition of the tangent space:
$$T_{(p,q)}X\times\X=V_{(p,q)}\oplus \lfg_{(p,q)}\oplus  J\lfg_{(p,q)},$$
where $\lfg_{(p,q)}$ denotes the induced tangent vectors on the product,  $J$
is  on the product and 
$$V_{(p,q)}=\{(u,v)|D_{(u,v)}\Psi=D_{J(u,v)}\Psi=0\}.$$

The decomposition  is orthogonal with respect to $$h(\cdot,\cdot)=-\omega(\cdot,J\cdot)+\omega_\X(\cdot,J\cdot)$$ on the product space.
\end{proposition}
{\it Remark:} The above decomposition would have been obvious had $h$ defined a complete norm.

{\it Pf:} 1). Claim $D_{J\lfg}\Psi$ is onto $\lfg^0$. 

 In the same notations as before, there is a decomposition
of $T_pX$ as $\cap_b\lfg/(\lfg)_b\oplus E_p$ where $E_p$ is of finite
dimension. By construction,  $J^X=J^\lfg$  and $J^\X=-J^\lfg$. 
Therefore, for $\xi\in \cap_b\lfg/(\lfg)_b$, 
$$D_{J\xi}\Psi=-D_{J^X\xi}\mu(p)+D_{J\xi}\mu_\X(q)=-2D_{J^\lfg\xi}\mu(p)=2\tad_{\mu(p)}J^\lfg\xi\in \cap _b\lfg/(\lfg)_b,$$
where one  degree of differentiability is lost due to the adjoint action. 
Because $J^\lfg$ preserves $\lfg$, and $\tad_{ \mu(p)} (\cdot)$ with
$\mu(p)\in k(C,1)$ is an
isomorphism on $\cap_{ \{O_b|\mu(p)\in O_b\} } \lfg/(\lfg)_b$,  the above implies that $D\Psi$ is onto
$\cap _b\lfg^0/(\lfg)_b$. The orthogonal  complement of $\lfg^0/(\lfg)_b$ is the subspace generated by
$(\lfg)_b$ such that $  \mu\in O_b$. If $D\Psi $ is not onto that finite
dimensional complement,
there is a $\eta$, with
$$(D\Psi|\eta)=0.$$
In particular, $$0=D_{J\eta}\Psi(\eta)=
-\omega(J\eta,\eta)+\omega_\X(J^\X\eta,\eta),$$
which forces $\eta(p)=0$ and 
$\eta(q)$  in the null direction of $\omega_\X(\cdot,J^\X\cdot)$. 
The analysis on the null space of $\omega_\X(\cdot,J^\X\cdot)$ shows that
$\eta\in(\lfg)^{\sss}_\mu$, where $\mu=\mu(p)=k\mu_\X(q)$.
But {\bf H2} is equivalent to $$(\lfg)_p\cap[(\lfg)_\mu,(\lfg)_\mu]=0, $$
a contradiction.

2). Given a tangent vector in the product, $(u,v)$, 
from the claim there exists $\xi\in\lfg$  such that $D_{J\xi}\Psi=D_{J(u,v)}\Psi$. 
And choose $\eta\in \lfg$ to satisfy $D_{J\eta}\Psi=D_{(u,v)}\Psi$, 
then $$(u',v')=(u,v)-J\eta(p,q)-\xi(p,q)$$ 
satisfies
$$D_{(u',v')}\Psi=0;\quad D_{J(u',v')}\Psi=0$$
because $D_\tau\Psi=-\tad_\tau\mu(p)+\tad^*_\tau\mu_\X(q)=0, \forall
\tau\in\lfg$.
Therefore $$(u',v')\in \ker D\Psi\cap J\ker D\Psi=V_p.$$

3). The isomorphism between $V_p$ and the tangent space to the orbit
$ \ker D\Psi/\lfg$ holds due to the claim in Step 1. Since $V_p$ is invariant
under $J$,  there exists $J$ on $TX_N$. QED

\section{The Fixed Point Set on $X_N$}

\subsection{When a point is fixed by $T$}

Let $u\in X_N$, then $u=[p,q]$ with $(p,q)\in X\times\X$, the bracket denotes
the equivalence class under the diagonal action by $LG$. 

The point $q$ being in $\X$ is itself an equivalent class, $q=[h,z]$ with
$h\in LG$, $z\in X_\fg$.
The equivalent class
represented by $q$ is the orbit by the diagonal $T$ action on $LG\times
X_\fg$. Here the $T$ action is from the right on $LG$. 

Both equivalence relations can be kept tracked of by introducing the
following:
\begin{definition}\label{equiv}
For points in $X\times LG\times X_\fg$,
define  the equivalence relation:
\begin{equation}
(p,h,z)\simeq (gp,ghs,sz)
\end{equation}
 for $ g\in LG, s\in T$.
\end{definition}
A triple defines a point in $X_N$ iff $\mu_X(p)=k\tAd_g(\tphi(z))$ where
$\tphi=(\phi,1)$ and $\tAd$ is the adjoint action.

 On $X\times LG\times X_\fg$, the diagonal action by $LG$ on the first two
factors   and the diagonal $T$-action on the last two factors, commute with   $T$ acting on the last one.   And the $T$-action  acting on $X_\fg$ alone  preserves the moment
maps $\mu-\Phi$,  therefore it descends  to one  on $X_N$. 

We make a semi-canonical choice for $[p,q]$, so that $q=[I,z]$. Because
$0=\Psi(u)=\mu_X(p)-k\Phi(q)$, we have $\mu(p)\in k (\ft,1)$. 
The ambiguity is due to the  following:
$$(p,I,z)\simeq (s^{-1}p,I,sz),\quad \forall s\in T$$
as in  Definition \ref{equiv}.
Thus  $(s^{-1}p,I,sz),s\in T$  represents the same point as $(p,I,z)$  in $X_N$. 
\begin{lemma}
Let $\ft_z$ be the Lie algebra of the stabilizer $T_z\subset T$ of $z\in X_\fg$,
$\ft_p$ be that for the point  $p$ in $X$. 
Then a point $u=[p,I,z]$ is fixed by $T$, 
 iff that $\ft=\ft_z+ \ft_p$. 
\end{lemma}
{\it Pf:} The condition is clearly sufficient. Suppose $u$ is fixed by $g$, i.e.
$(p,I,gz)\simeq (p,I,z)$, that is  for some $h\in T$,
$$(p,I,gz)= (h^{-1}p,I,hz).$$
 Thus $$hp=p ,\quad h^{-1}gz=z$$for $h\in T$. So $h\in T_p$ and $h^{-1}g\in T_z$, and $g=h\cdot h^{-1}g\in$. 
 Therefore, any $g\in T$ can be decomposed into a product of
elements in $T_p,T_z$,   hence  $\ft=\ft_z+ \ft_p$. QED

Assume $(p,I,z)$ defines a point in $X_N$,  then we have 
$\mu(p)=k\tphi(z)=k(\phi(z),1)$ by the definition of $X_N$.
\begin{lemma}\label{localfree}
The following holds:

1.) $\ft_z=(\lfg)_{\tphi}^{\sss}\cap \ft$.

 2). Under the assumption that the image of $\mu_X$ is transversal to $(\ft,k)$ in $(\lfg,k)$, 
$\ft_z\cap \ft_p=0$. 
\end{lemma}
 {\it Pf:}
Let $C_\mu$ be the smallest wall of $\partial C$ containing $\phi(z)$, and let $V_\mu$ be the linear subspace in $\ft$  parallel to $C_\mu$.
Then $V^\perp_\mu=\ft_z$, which is a basic fact in toric variety theory, or
symplectic geometry.

The affine Lie algebra based on $\fg$, $\fg^\aff$, is   $\lfg+\R d+\R K$
 as in [K].
The dual is $\lfg+\R \delta+\R \Lambda_0$.
The simple roots  are 
$\{\alpha_0,\alpha_1,...,\alpha_l\}$ where
$\{\alpha_1,...,\alpha_l\}$ is the set of simple roots of $\fg$, and 
$\alpha_0=\delta-\theta$ which acts on $(h,l)\in \ft+\R \Lambda_0\subset
\lfg+\R \delta+\R \Lambda_0$
as $l-\alpha(h)$.
The boundary of the  alcove $C$ is defined  
by $\cap\alpha_i^{-1}(0)$. 

 The stabilizer of $\tphi$ is the same as that of $C_\mu$. And the   stabilizer
$(LG)_\mu$ of $C_\mu$ in $LG$ is generated by $\ft$ and $\lfg_{\alpha_i}$ such that
$\alpha_i(C_\mu)=0$. 
Clearly $\lfg_\mu^\sss=\sum_{\alpha(C_\mu)=0}\lfg_{\alpha}.$
We claim that $$(\xi|\tphi)=0,\quad \forall \xi\in (\lfg)_\mu^\sss.$$
This is clearly  true for $x_\alpha,y_\alpha$ in $\lfg_{\alpha}$, using the standard notation for $\lfg^\C_{\alpha}\simeq sl_2(\C)$. As  for 
the coroot $\alpha^\vv=[x_\alpha,y_\alpha]$,  it holds  as well because
$\alpha(\mu)=0$ and 
$$(\alpha^\vv|\mu(p))=(x_\alpha|[y_\alpha,\mu])=\alpha(\tphi)(x_\alpha|x_\alpha)=0.$$

Thus we have  $\lfg_\mu^\sss\cap \ft=\sum_{\alpha_i(C_\mu)=0} \R\alpha_i^\vv$, the right hand side are in $\ft_z$ from the orthogonal condition just
proved. We also know that the $V^\perp_\mu$ is exactly spanned by those
$\alpha^\vv$, since $C_\mu$ as a subface  is defined  this way.
Hence $\ft_z=\lfg_\mu^\sss\cap \ft$

2).
 The transversality condition  is equivalent to 
$$\lfg_p\cap [\lfg_{\mu},\lfg_{\mu}]=\{0\},$$
but $\lfg_\mu^{\sss}=[\lfg_{\mu},\lfg_{\mu}],$ and $\ft_p\subset \lfg_p$, therefore we have
$\ft_p\cap [\lfg_{\mu},\lfg_{\mu}]=\{0\}$ which implies that
$\ft_p\cap\ft_z=0$ since $\ft_z\subset [\lfg_{\mu},\lfg_{\mu}]$. QED

Let $T_p^0$ be the connected component of $I$ in  $T_p$. By the construction
of the toric variety $X_\fg$, it is easy to see that $ T_z$ is connected.
If $\ft_p\oplus \ft_z=\ft$, naturally $(t_p,t_z)\in T_p^0\times T_z\mapsto
t_pt_z\in  T$ is a covering map.
\begin{proposition}\label{prop1.1}
Let $\tF_p$ be the
connected component
of $T_p^0$-fixed point set in $X$  containing $p$, and $M_z$  be the 
connected component of $T_z$-fixed point set in $X_\fg$  containing $z$.

Then the connected component containing $(p,I,z)$ of $T$-fixed point set in $X_N$ is given by $(\tF_p\times \{I\}\times M_z)\cap \Psi^{-1}(0)$.
\end{proposition}
{\it Pf:} The inclusion of the set in the connected component of the fixed
point set by $T$ is clear. To see the other direction, suppose
$(q,I,w)$ is in the connected component of $(p,I,z)$ which is fixed by $T$,
and
suppose $q$
is close to $p$, $w$ is close to $z$.
The stabilizers in $T$ of  $q,w$ must be  subgroups of $T^0_p,T_z$, which is
well known. On the other hand,
we have shown  $\ft_q+\ft_w=\ft, \ft_q\cap \ft_w=0$. Therefore the inclusion
of $\ft_q\subset \ft_p$, $\ft_w\subset \ft_z$ are equalities instead.
Thus, $\ft_p,\ft_z$ vanish at $q,w$ respectively.

If $(q,I,w)$ is
connected in the fixed point set to $(p,I,z)$ through a 1-parameter curve
$(q_s,I,w_s)$, with $q_0=p, w_0=z$ and $q_1=q$, $w_1=z$,  then  for small $s$,
by the above argument, $q_s,w_s$ are fixed by $\ft_p,\ft_w$. So they are in
the
connected components of $p,z$ fixed by $T_p,T_z$ respectively. The above
argument shows also the set $\{s\}$ such  that $(q_s,I,w_s)$ is in the desired
product  is open.
Obviously the set  is also closed. Therefore, $(q,I,w)$ is in the product. QED

\subsection{More about the $T$-fixed point set on $X_N$}
From earlier discussion, we have learned that $[p,I,z]$ is a
fixed point of $T$ on $X_N$  iff $\ft_p\oplus\ft_z=\ft$.
Let $T_p^0$ denote the connected component of $T_p$.
To understand the structure of the fixed point set, we need the following:
\begin{definition}
Suppose $\ft_z\neq 0$, i.e. $\phi(z)\in \partial C$,
let $$\KK=\{g\in(LG)_\mu\Big| \Ad_gt=t,\forall t\in T^0_p\},\quad
\NN=\{g\in(LG)_\mu\Big| \Ad_g(T^0_p)=T^0_p\},$$
and $\fkk$, $\fnn$ be their   Lie algebras.
\end{definition}
\begin{lemma}
The  groups $(LG)_\mu$ and  $\KK$  are   compact and   connected.

Suppose $(\lfg)^{\sss}_\mu\cap \ft_p=0$ which is true under the assumption
that the image of $\mu$ is transversal to $\ft$, then $\fnn=\fkk$.

 The  weights of the $\ft_p$ action on $(\lfg)_\mu/\kk$  are non-trivial,
therefore $\KK$ is the largest connected group acting on $\tF_p$.
\end{lemma}
{\it Pf:} The assertion on $(LG)_\mu$ is well known. For $\KK$, the argument
is also standard but we include here anyway. Suppose
$\KK^0$ is the connected a component passing
$I$.
Let
$g\in \KK\setminus \KK^0$, 
 then $\Ad_g\KK^0=\KK^0$. The group $T$ is a maximal torus in
$(LG)_\mu$ and is contained in $\KK^0$. Therefore $\Ad_gT$ is a maximal  
torus.
By the uniqueness of maximal torus under the
adjoint action, there is a $h\in \KK^0$ such that
$\Ad_{hg}T=T$.
Obviously  $hg\in W(\KK^0)$,
so $hg$ is contained in the semi-simple part of $\KK^0$.
  Thus $g$ is connected to $I$ and  $\KK=\KK^0$.

The adjoint action by $T_p^0$ on $\nn\setminus \kk$ has those roots of $\nn$ which
are not roots of  $\kk$ as eigenvalues. Let $\alpha$ be one of them, then
because $\ft_p$ is normalized by $\NN$,  the reflection
element in the  Weyl group
$r_\alpha$ satisfies $r_\alpha(\ft_p)=\ft_p$. Using the definition of the
reflection, and the fact that $\alpha(\ft_p)\neq 0$,  one concludes  the  coroot $\alpha^\vv\in \ft_p$.  On the other hand
$\alpha^\vv\in \ft_z$, since $\ft_z$ is the Cartan subalgebra of the
$\fg_\mu^{\sss}$ which contains $\fnn^{\sss}$.  But $\ft_z\cap\ft_p=0$,
 hence it is impossible that
$\alpha^\vv\in\ft_p$. So all the roots of $\fnn$ are those of $\kk$.  
From this and  the fact that $\ft$ is the Cartan sub-algebra of both groups,
one concludes $\nn=\kk$.

 Now the assertion that the action by $\ft_p$ on $\lfg_\mu/\kk$ have non-trivial weights
follows   immediately.
QED

\subsection{ $Z$ and  $Z/T_z$ as fiber bundles}
\begin{proposition}
Let $X$ be the connected component containing $p$ in $\tF_p\cap \mu^{-1}
(\tphi)$.
Then the connected
component of T-fixed point set on $X_N$ passing
$[p,I,z]$ is given by $Z/T_z \times \{z\}
$.

If $\phi$ is in the interior of $C$, $T_z=I$. 

In case $\phi$ is on the boundary of $C$,  
$Z/T_z$ admits a projection to an orbifold $E$ with
fiber a finite  quotient of  $\KK/T$.
\end{proposition}
{\it Pf:}
Suppose $[q,I,w]$ is in a same connected component as $[p,I,z]$ in $X_N$ fixed
by $T$. It follows from Prop.~\ref{prop1.1} that $q\in \tF_p$ and $w\in M_z$.
A basic property of moment map restricted
to fixed point set dictates that
$$\mu(\tF_p\cap\mu^{-1}(\tft) )\subset \mu(p)+\ft^\perp_p,\quad \phi(M_z)\subset
\phi(z)+\ft_z^\perp,$$
and $\ft^\perp_p\cap \ft_z^\perp=0$ since  $\ft_p\oplus\ft_z=\ft$.
The condition  $\mu(q)=k\tphi(w)\in \tft$ then forces
$$\mu(q)=\mu(p)= k\tphi(z)=k\tphi(w).$$
As points in the toric variety $X_\fg$, the equality $\mu_\X(z)=\mu_\X(w)$
implies that
$w=tz$ for some $t\in T$.
Therefore  $[q,I,w]$ has a representative
of the form $(t^{-1}q,I,z)$ which defines a point in
the claimed set. Thus the  association of $[q,I,w]$ with a point 
in the quotient space  $(\tF_p\cap\mu^{-1}(\tphi)/T_z)\times\{z\}$ is
  1-1 and onto.

As for the last assertion, we already know that $\tF_p\cap\mu^{-1}(\tphi)$
admits the action of $\KK$, since it commutes with $T_p$ and preserves
$\mu^{-1}(\tphi)$.
Let $E=Z/\KK$. The fiber of the projection
$$\pi: \left( \tF_p\cap\mu^{-1}(\tphi)\right)/T\rightarrow \left(
\tF_p\cap\mu^{-1}(\tphi)\right)/KK=E$$ can be explicitly described. If
$[q]\in E$ with $q\in \tF_p\cap\mu^{-1}(\tphi)$, then the fiber
is
$\KK(q)=\KK/\KK_q.$  We know that $\KK_q$ has $T_p^0$ as  the
maximal connected subgroup,  since $\kk_q$ has no semi-simple part from 
the condition $\lfg_q\cap [\lfg_\mu,\lfg_\mu]=0$.
 Hence $T_p^0$ is a normal
subgroup and $\KK_q/T_p^0$ is a finite subgroup. Therefore the fiber
is a finite quotient of the homogeneous space $\KK/T$ by $K_q/T_p^0$.
QED

The following obviously holds:
\begin{equation}
Z=\tF_p\cap \mu^{-1}(\tphi)\longrightarrow \tF_p\cap \mu^{-1}(\tphi)/T
\longrightarrow \tF_p\cap \mu^{-1}(\tphi)/\KK=E.
\end{equation}
Since $T^0_p$ fixes points on $Z$, the action by $\KK$ is not effective on $Z$.
In fact, we have
\begin{lemma}
  Define $\ft^0$ to be a complement of $\ft_z\cap \kk^{\sss}$ in 
$\ft_z$, then $$ \kk=\kk^{\sss}\oplus \ft_0\oplus \ft_p, $$where $\kk^{\sss}$ is the semi-simple part of $\kk$.
\end{lemma}
{\it Pf:} It is known that $\ft\subset \kk$, so $\kk=\kk^{\sss}+\fh$ with 
$\fh\subset \ft$.  On the other hand, $[\kk^{\sss},\ft_p]=0$ by definition
of $\KK$, hence $\kk^{\sss}\perp \ft_p$. 
But $\kk^{\sss}\subset (\lfg)_\mu^{\sss}$ which has $\ft_z$ as its Cartan
subalgebra. Therefore $\kk^{\sss}\cap\ft \subset \ft_z$. 
We already knew that $\ft=\ft_p\oplus\ft_z$, hence
$$\kk=\kk^{\sss}\oplus \ft_p\oplus \ft_0$$
where $\ft_0\subset \ft_z$ and is perpendicular to $\ft\cap \kk^{\sss}\subset
\ft_z$. QED

Let $\KK'$ be the group with $\kk=\kk^{\sss}\oplus \ft_0$ as its Lie algebra.
This group  acts effectively on $Z$,
and $T_z$ is its Cartan subgroup.

So the fiber over $E$ is a finite quotient of $\KK'/T_z$.

\subsection{Connections of the orbifold fiberations}
Let  $\Lie \KK'=\ft_z+\fn$ be the Cartan decomposition.

Fix a $\KK'$-stable splitting of $TZ=T^HZ\oplus T^\perp Z$, where $T^\perp Z$
is tangent to the orbit by $\KK'$ on $Z$. Such a splitting exists because the
action by $\KK'$ is locally free, therefore  there is  automatically the bundle $T^\perp Z$. The $\KK'$-stable horizontal space can be obtained by an invariant
metric on $TZ$.  

Let  $A:TZ \rightarrow  T^\perp Z$ be the projection. 
There is a further decomposition:
$$T^\perp Z=T'Z+T''Z,$$
where $T'Z$ and $T''Z$ are  the vertical subbundle generated by vectors in $\ft_z$ and $\fn$ respectively.

\begin{lemma} 

1). If $(p,n)\in Z\times \fn\mapsto n(p)\in T_p''Z$, then the action
by $T_z$ on $T''Z$ induced from the action on $Z$,  maps $(p,n)$ to $(pt,
\Ad_t n),\forall t\in T_z$. 

2). $T''(Z/T_z)=Z\times _{T_z}\fn$ where the quotient is taken with
the previous $T_z$ action on $T''Z\times \fn$.
 And  $T(Z/T_z)=\pi^*TE\oplus T''(Z/T_z)$,
where $\pi :Z/T\rightarrow E=Z/\KK'$.
\end{lemma} 
{\it Pf:} Let $g_*$ denote the action on the tangent bundle.
 It is straightforward to verify that
\begin{equation}\label{vertmove}\xi(gp)=g_*(\Ad_{g^{-1}}\xi)(p)),
\forall g\in \KK', \xi\in \kk'.
\end{equation}
Hence $(gp,\xi)\mapsto (p,\Ad_{g^{-1}}\xi)$.
Or $(p,\xi)\mapsto (gp,\Ad_g\xi)$,
which implies the first assertion after applying it to $\xi=n\in \fn$ and
$g=t\in T_z$.

2). This is simply the decomposition of $T(Z/T_z)$ into horizontal and 
vertical parts. QED

Denote the map which  identifies the  vertical vectors  with $\kk'$ by $P$.
Combine with $A$,  we
obtain $$P A: TZ\rightarrow \kk'.$$
Furthermore $P=P'+P''$ with $P',  P''$ take values in $\ft_z,\fn$
respectively. Obviously $PA $ is a connection on $Z$.

Let $u,v$ be  vector fields in $T^HZ$ invariant under $\kk$,
 and let $\xi,\eta$ be vectors induced by two elements from $ \fn$.

The decomposition of the curvature given below will be useful later.
\begin{lemma}
 Suppose $T_z$ acts on $\C$ with character $\lambda$, then
 there is a connection $A_\lambda$ on the bundle $Z\times_{T_z}\C\rightarrow Z/T_z$
such  that its curvature is given by $\lambda(B)+\lambda(R)$;
 the 2-forms $B$ and $R$ are
$\ft_z$-valued, and  satisfy the following  $$B(u,v)= P'\cdot A ([u,v]), \quad 
R(\xi,\eta)= P'\cdot A([\xi,\eta]).$$
\end{lemma}
{\it Pf:} On the trivial bundle $Z\times \C$, the action by Lie algebra of
$\ft_z$ is given by ${\mathcal L}_t=(t(q),\lambda(t))$, for $t\in \ft_z$. The connection
defined by $\nabla=d+\lambda(P' A)$ satisfies 
$$\nabla_{t(q)}=d_{t(q)}+\lambda(P' A)(t(p))=d_{t(q)}+\lambda(t)={\mathcal
L}_t$$ by the 
definitions of $P,  A$. A connection,  satisfying the above condition,    descends
to the quotient bundle $Z\times_{T_z}\C$. The curvature is given by 
$d(\lambda(P' A))$.
For vector fields $u, \xi$,   by invariance of $u$, we know that
$$[u,\xi]=0 ,\quad A(u)=0, \quad PA(\xi(q))\in \fn$$ 
hence $P'A(\xi(q))=0$ and the differential of the 1-form, $d(\lambda(P' A))(u,\xi)=0$. Thus the curvature has no term mixing
$T^HZ$ and $ T''Z$, and  consists only of the vertical and the
 horizontal  part. 
The one coming from the horizontal vector fields $u,v$ is exactly
$$\lambda\cdot B(u,v)=\lambda\cdot P'A([u,v])=\lambda\cdot PA([u,v]),$$ 
since $A(u)=A(v)=0$ and $\lambda\cdot P''=0$. The one coming from the fiber, or the vertical directions,  is
$\lambda(PA([\xi(q),\eta(q)])=\lambda(P'A([\xi(q),\eta(q)])$. This follows from   $\lambda(\fn)=0$ and $P'A(\xi(q))=P'A(\eta(q))=0$, hence in calculating $dP'A(\xi,\eta)$ only the said
term remains.
 QED

Fix a base point $q$, one may ask how the forms $B, R$ transform along the 
fiber of $\pi$. 

\begin{lemma}
\begin{eqnarray*}\label{damn}
\lambda\cdot B(u,v)|_{gp}=\lambda (\Ad_gPA[u(q),v(q)]),\\
\lambda\cdot R(\xi(gq),\eta(gq))=\lambda (\Ad_g[\xi(q),\eta(q)]).
\end{eqnarray*}
\end{lemma}
{\it Pf:} Due to the invariance of $u,v$ and the fact that $g_*$ commutes with
$[,]$, 
we have  $[u(gq),v(gq)]=g_*([u(q),v(q)])$.
The map $g_*$ also commutes with $A$ by the invariance of $A$, therefore 
$$A[u(gq),v(gq)]=g_*(A([u(q),v(q)])).$$
Let $a=PA([u(q),v(q)])\in\kk'$, then its induced vector  $a(q)=A([u(q),v(q)])$ by the definition of 
$P$.  We know $g_*(a(q))=\Ad_g a(gq)$, 
thus
\begin{equation}\begin{split}
 PA[u(gq),v(gq)]&=Pg_*(A([u(q),v(q)]))\\
&=Pg_*(a(q))\\
&=P\Ad_g(a)(gq)\\
&=\Ad_g(a),
\end{split}\end{equation}
hence $\lambda(PA[u(gq),v(gq)])=\lambda (\Ad_gPA[u(q),v(q)])$.
For the same reason, we obtain
$$\lambda\cdot R(\xi(gq),\eta(gq))=\lambda \cdot\Ad_gPA([\xi(q),\eta(q)]). \QED
$$
This expression is crucial to  a calculation in Section 9, because 
$\lambda R$ behaves as the moment map on the fiber which is a coadjoint 
orbit.

\subsection{Stratification of the fixed point set}
Unlike the smooth case, the fixed point formula for orbifolds requires
the contributions of  lower strata of the fixed point set. Where do the
strata come from? They are present due to  the local isotropy groups on
the fixed point set $Z/T_z$. 

To describe them locally, let $p\in Z$,  $T_p$ may not be connected. Denote
 the connected component of $I$ by  $T^0_p$. 
Let $I_p=T_p\cap T_z$, $I^0_p=T^0_p\cap T_Z$. They are finite groups since
$\ft_p\cap \ft_z=0$. 
 When $T_z$ acts on $Z$, it
has $I_p$ as its stabilizer. All $q\in Z$  are fixed by $T_p^0$, so the
subgroup $I_p^0$ acts trivially on $Z$. The effective isotropy group on $Z$ is
$I_p/I_p^0$, though $I_p^0$ may have non-trivial action  on the normal
bundle and can not  be ignored.
Obviously the discussion is unnecessary if for all $p\in Z$, $T_p=T_p^0$.

For each $h\in I_p/I_p^0$, denote by $Z_h$ its fixed points in $Z$.
The  collection   $\{Z_h\}$ for $h\in I_p/I_p^0, \forall p\in Z$ form
stratification of $Z$. And their   quotient by $T_z$ in $F=Z/T_z$ contribute
to the fixed point formula computations, which is different from the 
smooth case.

Let $h\in I_p\setminus I_p^0$, let $Z_h=\{q\in Z|h(q)=q\}$.  
It is a submanifold.
  Clearly $T_z$ acts on it, and the points there are  fixed by $h, T_p^0$.
Let $\KK_h$ be the connected subgroup of $\KK$ which commutes with $h$, i.e.,
$\Ad_kh=h, k\in \KK_h$. Then  $\KK_h$ acts on $Z_h$, and it contains
$T$. 

In the language above, $Z$ itself can be thought of as $Z_h, h\in I^0_p$.

\subsection{ Lower stratum $Z_h$ and $Z_h/T_z$ as fiber bundles}
Obviously one has 
\begin{lemma}
The Lie algebra of $\KK_h$, $\kk_h$ is the maximal  subspace on which the action $\Ad_h|\kk$ is $I$.
\end{lemma}

This lemma  implies that $\kk/\kk_h$ induces a subspace
normal to  $T_pZ_h$ in $T_pZ$.

  As in the case of $Z/T_z$,  we can 
 realize $Z_h/T_z$ as a fiber bundle. 
As in the case of $\kk$, $\kk_h$
splits into $\kk'_h$ and $\ft_p$.
And the group $\KK'_h\simeq \KK_h/T_p^0$ acts effectively on $Z_h$.
 Associate with $Z_h$ the space  $E_h=Z_h/\KK'_h$.
Similar to $Z$, there is the following
sequence:
$$Z_h \rightarrow Z_h/T_z\rightarrow Z_h/\KK'_h=E_h$$
where the second projection yields a finite quotient of $\KK'_h/T_z$
 as the fiber.

\subsection{The action by a Weyl subgroup on the fixed points of $T$}
Suppose $(p,I,z)$ defines a fixed point of $T$ in $X_N$, that is to say
that $\ft_p,\ft_z$  generate $\ft$, and $\mu(p)=k\tphi(z)$.  

If $\phi(z)$ is in the interior of $C$, then $\ft_z=0$, and $\ft_p=\ft$ by the
above characterization of a fixed point. 

Suppose  $\phi(z)\in \partial C$, let $C_\phi$ be the smallest wall of $\partial
C$ containing $\phi$. Then $\lfg_A=(\lfg)_\mu^{\sss}$ 
commutes  with $\mu(p)=\phi(z)$, and it  is generated by
$\lfg_\alpha$ with $\alpha$ vanishing on the wall $C_\phi$. 

\begin{lemma}
The  subgroup $W^{\aff}\cap(LG)_\mu^{\sss}$  of $W^\aff$
is the
Weyl group of the finite dimensional semi-simple group $(LG)_\mu^{\sss}$.
It 
transforms one  fixed point set to another
in $\mu^{-1}(k\tphi)$ where $\tphi=\mu(p)/k$. 

For $w\in W(\KK')$, it preserves $Z$ but permutes among $\{Z_h\}$.
\end{lemma}
{\it Pf:} The group $(LG)_\mu$ is compact and connected, since it is the
stabilizer of $\mu(p)=(\phi(z),k)$  under the adjoint action. 
The semi-simple part is generated by $(\lfg)_\alpha$ where $\alpha$ vanishes
on $C_\phi$, therefore the Weyl group of $(LG)_\mu^{\sss}$ is generated by
the reflections with respect to those affine roots vanishing at $C_\phi$.  In
terms of the original affine Weyl group, it is exactly
$W_\mu=W^{\aff}\cap(LG)_\mu^{\sss}$.

 Suppose that $(p,I,z) $ defines  a fixed point, then 
$\mu(gp)=\tAd_{g}\mu(p)=\mu(p)$, for all $g\in (LG)_p^{\sss}$.
If $g$ is also in $W_\mu$, one has $g(\ft_z)=\ft_z$ because 
$\ft_z$ is the Cartan subalgebra of $(\lfg)_\mu^{\sss}$. Thus
$g(\ft_p)$ and $g(\ft_z)$ generate $\ft$.
To check that  
$(gp,I,z)$ is fixed by $T$,  let $ t\in\ft$, $$t=t_1+t_2, \quad
t_1\in g(\ft_p),\, t_2\in g(\ft_z)=\ft_z.$$
Then 
\begin{displaymath}\begin{split}
\exp(\lp t)(gp,I,z)&=(gp,I,\exp(\lp t)z)\\
&\simeq (\exp(\lp t_1)gp,I,\exp(\lp t_2)z)\\
&=(g\exp(\lp \Ad_{g^{-1}}t_1)p,I,z)\\
&=(gp,I,z).
\end{split}\end{displaymath}
In the above, we have used the fact that
$\exp(\lp\Ad_{g^{-1}}t_1)\in T_p$ if $t_1\in g(\ft_p)$.
This shows the  Weyl group of the semi-simple part of $\mu(p)$ acts 
on the fixed point set whose image  is $\mu(p)$.

Fix a lifting of $w$ to $\KK$, it certainly preserves $Z$, since the whole
$\KK$ does. If $Z_h$ is a stratum as described earlier, associated with
$h\in I_p/I^0_p$, for some $h$ and $p$, then the point $w(p)$ has $\Ad_wT_p$ as stabilizer,
and its isotropy group is $\Ad_wI_p$.
Since $T_z$ is the Cartan subgroup of $(LG)_\mu^{\sss}$,  and $w$ is in its Weyl
group, hence  $w$  preserves $T_z$ and  $$\Ad_wI_p=I_{w(p)},\quad \Ad_wI^0_p=I^0_{w(p)}.$$
Thus $w(Z_h)=Z_{w(h)}$.
QED

\section{Normal bundles to the fixed point  sets in $X_N$ and weights}
We will find out in this section the weights of the $T$-action on  the normal
space to  the   fixed point sets inside the compactification locus.

   The last section describes 
  the fixed  point sets
of $X_N$ in terms of data from $X$ and $\X$. Suppose $(p,I,z) $ defines a fixed point in $X_N$, and
$(p_s,h_s,z_s)$ is a curve with 
$$(p_0,h_0,z_0)=(p,I,z),\quad (p_0',h_0',z_0')=(a,\xi,x)\in T_{(p,h,z)}(X\times
LG\times X_\fg).$$  Let $t=t_p t_z$ so that
$t_p\in T_p$ and $t_z\in T_z$. 
Then the action by $T$ on the tangent vectors is given by
\begin{equation}\label{Tanaction}
t(a,\xi,x)=(a,\xi,t(x))\simeq (t_p(a),\Ad_{t_p}\xi,t_z(x))
\end{equation}
which is obtained by the definitions of the action of $T$ on $X_N$, and 
 the equivalence relation as in  Def. \ref{equiv}.

 By assumption, the equivalent class $[p,I,z]$ defines  a fixed point in $X_N$.
As before,   $T^0_p$
is  the connected component of $I$ in $T_p$.  
 Then $T^0_p$ and $T_z$ generate
$T$. Furthermore $T_p$ and $T_z$ act on the tangent spaces
$T_pX$ and $T_z\X$, respectively.
\subsection{Tangent space }
First let us describe  the tangent space to $X_N$
at $[p,I,z]$.  
\begin{proposition}\label{tangentsp}
Let $\mu=\mu(p)$ and $V_p$ be  defined as 
$$V_p=\{a\in T_pX|\exists b\in\lfg,    D_{(a+iJa)}\mu(p)=D_{(b+iJb)}\muX\}$$
where $\muX:\X=LG\times_T X_\fg\rightarrow\lfg\times\{k\}$. 
In the toric variety $X_\fg$, the point $z$ has $T_z$ as its stabilizer, let
the subspace tangent to $T^\C(z)\simeq (\ft_z^\perp)^\C$ be denoted by $H_z$,
then
the tangent space to $X_N$ at $[p,I,z]$ admits the following decomposition:
$$T_uX_N\simeq V_p\oplus (\lfg)_\mu/\ft\oplus H_z$$
where $(\lfg)_\mu$ is the Lie algebra of the stabilizer of the $\mu$ under 
the coadjoint action.

The tangent space $T_uX_N$ has a natural almost complex structure $J$ 
satisfying:
$$J|V_p\simeq J^X;\quad J|(\lfg)/\ft=J'=-J^\lfg; J|H_z=-J^{X_\fg}$$
where $J^X, J^\lfg, J^{X_\fg}$ are the almost complex structures  on the  
space $X, \lfg/\ft, X_\fg$ respectively. 
\end{proposition}
{\it Remark:} 
  The map $D\mu$ after extended to the complexified tangent space,
certainly is not holomorphic, thus the vectors in  $V_p$ are very special.

{\it Pf:} The tangent space of $T_uX_N$ is isomorphic to the space
$$\{(a,b)\in T_{(p,q)}X\times \X|D_{(a+iJa)}\mu(p)=D_{(b+iJ'b)}\muX(q)\},$$
as shown in Section 2 (or cf. [C1]). The subspace $H_z$ is contained there via
$b\mapsto (0,0,b)$, for $D_{(0,b+iJ'b)}\mu_\X(I,z)=0$.
Let $q=[I,z]\in\X$.
The subspace $\lfg_\mu/\ft$ is contained there via the inclusion
$\xi\in \lfg_\mu/\ft\mapsto (0,\xi,0)$, and it satisfies the equation
$$D_{(\xi+iJ'\xi,0)}\muX(I,z)=D_{(\xi-i\xi^J,0)}\muX(I,z)=[\xi-i\xi^J,\muX(I,z)]=0,$$
since $[\xi,\muX(q)]=[\xi^J,\muX(q)]=0$.  Also  we have used the fact that
$J'\xi(q)=-J^\lfg\xi(q)=-\xi^J(q)$ with $J^\lfg\xi=\xi^J\in \lfg_\mu/\ft$ being defined by the almost complex
structure 
on $\lfg_\mu/\ft$.

The subspace defined by $(\lfg)_\mu/\ft\oplus H_z$ corresponds exactly to
the subspace  $$U_0=\{(a,b)\in T_{(p,q)}X\times \X|a=0\}.$$
The complement of $U_0$ in $T_uX_N$ is exactly $V_p$.

There is a unique way to choose $b$ in the above definition of $V_p$: make it perpendicular to $H_z$, 
and choose it from $$\sum_{\alpha\notin \Delta(\lfg_\mu)}\lfg_\alpha $$
where the index means $\alpha$ is not a root of $\lfg_\mu$. 
This is possible since $D_b\muX=D_{J'b}\muX=0$ if $b\in H_z$ or 
$b$ is induced from $\lfg_\mu$. 

Then the map  $h: a\mapsto (a,b)$  satisfies $h\cdot J=J'\cdot h$. So the 
complex structure defined by $J^X$ on $V_p$ is mapped to $J=(J^X,J')$ restricted to the image of $h$. Thus the assertion is verified.
The  other two  idenities involving $J'$ follow strictly from the description of $J'$ on $T\X$.  QED


\subsection{Normal space to $Z$ and $Z_h$}
From the description of $V_p$ and $\lfg_\mu/\ft$, it is clear they inherit an almost
complex structure from $T_{[p,I,z]}X_N$. In the following  the weights
refer to the weights by  $\ft_p, h\in I_p\setminus I_p^0$ acting on 
the complex linear spaces.

Let $V^0_p$ and  $N_p$ be  the 0-weight and non-zero  weight subspaces of 
$V_p$ under the linear $t_p$-action respectively.
Then $V_p=N_p\oplus V^0_p$.

Similarly, denote by $\nno(Z_h,Z) \subset \{ V^0_p\}$  the maximal $h$-stable
subbundle 
subspace on which $\det(I-h)\neq 0$.

\begin{proposition}

1).  The normal subspace to the fixed point set $Z/T_Z$  is
$$N_p\oplus (\lfg)_\mu/\kk\oplus H_z.$$

2).  The normal bundle of $Z_h$ in $Z$ can be decomposed as
$$\nno(Z_h,Z)\oplus \kk/\kk_h.$$
The bundle $\nno(Z_h,Z)$  admits an action by $T_z$ which lifts the
action on $Z_h$.

On the normal subspaces in both situation, there is the induced almost
complex structure from $J$ on $T_uX_N$.
\end{proposition}
{\it Pf:} 
The first claim holds due to the fact that the normal space 
to the fixed point set is the non-zero weight space under the action by
the $\ft_p$ on the tangent space. From the description of the tangent space
and the definition of $N_p$, the assertion is evident. 

Part 2 follows from essentially the same reasoning, with a slight variation
since here we consider only the $h$-action.  The the normal space
is the maximal  $h$-stable subspace in $T_pZ$, on which 
1 is not an eigenvalue. 
There are two factors in the tangent to $T_pZ$, one is the $\kk/\ft$,
the other is $V^0_p$.
Then from the definition of $\nno(Z_h,Z)$ and Lemma 3.8, the decomposition
holds.

 The action by $T_z$ on $\nno(Z_h,Z)$
exists since $T_z$ commutes with $h$, so the non-zero weight space under $h$ in $V_p$ at $p$ is 
isomorphic to that in $V_{tp}$, $t\in T_z$.

The almost complex structure preserves those subspaces since 
$\ft_p$ and  $h$-actions  commute with $J$.  QED

\subsection{Description of weights}

From the last section, it is verified that $\ft=\ft_p\oplus\ft_z$ and $\ft_p\cap\ft_z=0$.
Let $t=t_p+t_z $,  $t\in \ft$ denote this  decomposition.

\begin{proposition}\label{asswei}
Let the weights of the $T_p$-action on $N_p$ be denoted by $\{\gamma\}$, the 
weights of the $T_z$-action on $H_z$   denoted by $\{\lambda\}$,  and the
positive roots of  $(\lfg)_\mu/\kk$ be denoted by $\{\beta\}$ respectively.

Then the weights by the   $T$-action on $T_{[p,I,z]}X_N$ are given by  the
corresponding three sets of weights:
$\{\tgamma\}$, $\{\tlambda\}$, $\{\tbeta\}$ such that
\begin{equation}\label{defwei}
\tgamma(t)=\gamma(t_p),\quad \tlambda(t)=-\lambda(t_z),\quad \tbeta(t)=-\beta(t_p).
\end{equation}
\end{proposition}
{\it Remark:} The negative sign for $\tlambda$ and $\tbeta$ reflects the 
complex structure on $X_N$ as described in Prop. 4.1.

{\it Pf:} The action by $t$ on the tangent space $T_{[p,I,z]}X_N$  is of the form:
$$t_p\cdot t_z (a,\xi,x)=(t_p(a),[t_p,\xi],t_z(x))$$
as shown by Eq. (\ref{Tanaction}). Now apply this to the three types of 
normal vectors as in  the last proposition, one has the assertion. 
The only thinking needed here is the observation that the embedding of
$V_p$ into $T_uX_N$ is $T_p$ equivariant.
The `-' sign reflects the choice of the complex structure on $H_z,
\lfg_\mu/\ft$. QED
\begin{figure}
\begin{center}
\setlength{\unitlength}{0.00041667in}
\begingroup\makeatletter\ifx\SetFigFont\undefined%
\gdef\SetFigFont#1#2#3#4#5{%
  \reset@font\fontsize{#1}{#2pt}%
  \fontfamily{#3}\fontseries{#4}\fontshape{#5}%
  \selectfont}%
\fi\endgroup%
{\renewcommand{\dashlinestretch}{30}
\begin{picture}(6012,3939)(0,-10)
\drawline(3000,1812)(4800,1212)
\drawline(4800,1212)(5100,612)
\drawline(3000,1812)(3900,1512)
\drawline(3776.671,1521.487)(3900.000,1512.000)(3795.645,1578.408)
\drawline(2100,2112)(1200,2412)
\drawline(1200,2412)(300,2712)(1500,3012)
\drawline(3000,1812)(3900,2412)
\drawline(3816.795,2320.474)(3900.000,2412.000)(3783.513,2370.397)
\drawline(4800,1212)(6000,1512)
\drawline(300,2712)(600,2112)
\dottedline{45}(3000,1812)(3300,1212)
\drawline(3219.502,1305.915)(3300.000,1212.000)(3273.167,1332.748)
\drawline(3000,1812)(3000,1212)
\drawline(2970.000,1332.000)(3000.000,1212.000)(3030.000,1332.000)
\dottedline{45}(3000,1812)(4200,2112)
\drawline(4090.859,2053.791)(4200.000,2112.000)(4076.307,2112.000)
\drawline(2100,2712)(2100,2712)
\drawline(2100,2712)(2100,312)
\drawline(1800,3312)(1800,3312)
\drawline(1800,3312)(1800,3312)
\drawline(3600,3912)(3600,3912)
\drawline(2100,2712)(3900,3912)
\drawline(2100,312)(2100,12)(3000,612)
\dottedline{45}(3000,1812)(3900,1212)
\drawline(3783.513,1253.603)(3900.000,1212.000)(3816.795,1303.526)
\dashline{60.000}(3000,1812)(2100,2112)
\drawline(2223.329,2102.513)(2100.000,2112.000)(2204.355,2045.592)
\put(0,2712){\makebox(0,0)[lb]{\smash{{{\SetFigFont{6}{7.2}{\rmdefault}{\mddefault}{\updefault}$Q$}}}}}
\put(3600,2412){\makebox(0,0)[lb]{\smash{{{\SetFigFont{6}{7.2}{\rmdefault}{\mddefault}{\updefault}$\lambda_1$}}}}}
\put(2700,1212){\makebox(0,0)[lb]{\smash{{{\SetFigFont{6}{7.2}{\rmdefault}{\mddefault}{\updefault}$\lda_2$}}}}}
\put(4275,2112){\makebox(0,0)[lb]{\smash{{{\SetFigFont{6}{7.2}{\rmdefault}{\mddefault}{\updefault}$\tilde{\lda}_1$}}}}}
\put(3900,1512){\makebox(0,0)[lb]{\smash{{{\SetFigFont{6}{7.2}{\rmdefault}{\mddefault}{\updefault}$\beta$}}}}}
\put(3375,1137){\makebox(0,0)[lb]{\smash{{{\SetFigFont{6}{7.2}{\rmdefault}{\mddefault}{\updefault}$\tilde{\lda}_2$}}}}}
\put(3975,1062){\makebox(0,0)[lb]{\smash{{{\SetFigFont{6}{7.2}{\rmdefault}{\mddefault}{\updefault}$\tilde{\beta}$}}}}}
\put(2400,612){\makebox(0,0)[lb]{\smash{{{\SetFigFont{6}{7.2}{\rmdefault}{\mddefault}{\updefault}$P$}}}}}
\put(2250,2262){\makebox(0,0)[lb]{\smash{{{\SetFigFont{6}{7.2}{\rmdefault}{\mddefault}{\updefault}$\gamma$}}}}}
\end{picture}
}
\end{center}
\caption{Weights at the fixed point set}
\end{figure}
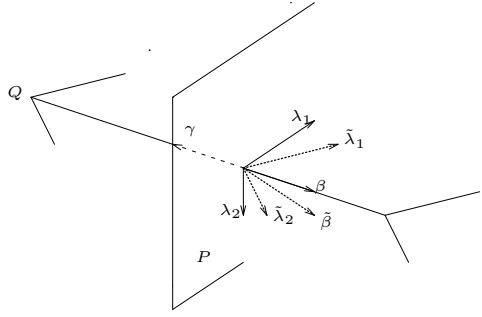

\subsection{Weights on normal space to the lower strata $Z_h$ in the fixed
point set}
As pointed earlier, the normal space to $Z_h$ consists of two parts: the
normal space to $Z$ and the normal space of $Z_h$ in $Z$. 
The contribution of $Z_h$ to the fixed point formula, in the orbifold sense,
depends on the weights of the action by $hT_p^0$ on the normal space. 
The calculation given earlier provides the answer for the action by $T$ on
the normal space to $Z$, here we need only determine 
 the
weights on
the second subspace  which is $\nno_p(Z_h,Z)\oplus (\kk/\kk_h).$
Recall that $h\in I_p/I_p^0=(T_p\cap T_z)/(T_p^0\cap T_z)$.
\begin{lemma}\label{stratawts}
Suppose the weights of the action by $h$ on $\nno(Z_h,Z)$, $(\kk/\kk_h)$
are given by $\{\theta_i\} $, $\{\beta_i\}$ respectively.
Let $(ht_p,h^{-1}t_z)\in hT_p^0\times T_z$ be a lifting of $t\in T$ in
$T_p\times T_z$,   $\{\ttheta\},\{\tbeta\}$ be the weights of  the action by 
$(ht_p, h^{-1}t_z)$
   on the two factors of the normal space. They are given by
$$\ttheta(ht_p\cdot h^{-1}t_z)=\theta(h),\quad \tbeta(ht_p\cdot h^{-1}t_z)=\beta(h).$$
\end{lemma}
{\it Pf:} Let $(p_s,k_s,z)$ be  a curve, with $p_s, k_s$ tangent to 
$Z,\KK'$ but normal to $Z_h,\KK'_h$ and $p_0\in Z_h, k_0\in \KK'_h$; 
and  $z$ is fixed by $T_z$.
The action by $t\in T$ on that curve after choosing the lifting $ht_p\cdot
h^{-1}t_z$, is simply 
$$ (ht_pp_s,\Ad_{ht_p}k_s,z), $$
which  is $(hp_s,\Ad_{h}k_s,z)$ since 
$t_p$ acts trivially on $Z$, and $\Ad_{t_p}k=k, \forall k\in\KK$ by the
definition of $\KK$.
Therefore, the weights are   of the forms described.
QED

\subsection{Three types of fixed points}

Suppose $(p,I,z)$ defines a fixed point in $X_N$,  we classify them according
to the following:

1). $\mu(p)/k=\tphi(z)$ is in the interior of $(C,1)$.

2). $\mu(p)/k=\tphi(z)$ is on $(\partial C,1)$ but not a vertex of $C$.

3). $\mu(p)/k=\tphi(z)$ is one of the vertices of the simplex $(C,1)$.

Considering  the decomposition of $\ft$ into $\ft_p,\ft_z$, the type 1) and 3)
correspond to the cases $\ft_z=0$ and $\ft_p=0$ respectively.

\subsection{The $(LG)_\mu/T$ factor}
In case of type 3), $\ft_z=\ft$, therefore $\ft_p=0$. Hence $\KK=(LG)_\mu$,
and $(LG)_\mu$ acts on the fixed point set, as shown for $\KK$. 
What's said earlier about the map $\pi:F\rightarrow E$ with fiber $\KK'/T_z$
now can be replaced by the factor $(LG)_\mu/T$. 

\subsection{More on $W^\aff$}
Suppose the fixed point is of  either type 2) or 3), then $\mu(p)$ has stabilizer
$(LG)_\mu$, under the co-adjoint action. The group has semi-simple part
generated by $(\lfg)_\alpha$, where $\alpha$ is an affine root 
vanishing
 at $\mu(p)/k$.

The following is a well known fact in affine algebra:
\begin{lemma}
The group generated by the reflections
$$r_\alpha=\{\alpha|\alpha(\mu)=0\}$$  is the Weyl  group  of  $T$ in $(LG)_\mu$, it acts on the type 2) fixed point set. Denote the Weyl group $W_\mu$ which 
is the subgroup of $W^\aff$ fixing the smallest wall containing $\mu$.
\end{lemma}

\subsection{Transformation weights by $W_\mu$. }
The discussion here is only meaningful for type 2) and 3) fixed points.

In the previous section, it was shown that the Weyl group $W_\mu$ permutes 
among the collection of fixed point sets with the same value for $\mu$.

Let $\{\tlambda'\}, \{\tbeta'\}, \{\tgamma'\}$ denote the three groups of  weights at
the point $(wp,I,z)$ as in Prop. \ref{asswei}.

Using the expressions in Eq.~(\ref{Tanaction}) and  Prop. \ref{asswei},  we conclude  that
all the non-trivial weights are those $\tlambda$'s in case of type 3) fixed
point, since $\ft_p=0$ and $\KK=LG_\mu$, so there is no  $\gamma,\beta$.
 Hence, along the $(LG)_\mu$-orbit of $(p,I,z)$ which is fixed by $T$ in
$X_N$, the weights are just those from $H_z$, $\{\lambda\}$.

In case of type 2) fixed point, we have the following relation between
the stabilizers:
$$T_{w(p)}=wT_p.$$
To derive the transformation rule for the weights, we need the following:
\begin{lemma}
$$V_{w(p)}=w_*(V_p).$$
where $w_*:T_pX\rightarrow T_{w(p)}$ is the isomorphism induced from the 
diffeomorphism $w$.
\end{lemma}
{\it Pf:} Let $a\in V_p$, i.e. $\exists \xi\in \lfg/\ft$ such that  $D_{a+iJa}\mu(p)=D_{\xi+iJ'\xi}k\tphi(z)$.
Then using the invariance of $J$ under $w$, we have
\begin{equation}\begin{split}
D_{w(a+iJa)}\mu(w(p))&=w(D_{a+iJa}\mu(p))\\
&=w(D_{\xi+iJ'\xi}k\tphi)\\
&=D_{w(\xi)+iJ'w(\xi))}\mu,
\end{split}
\end{equation}
 since $w$ preserves $\ft$, we have $w(\xi)\in\lfg/\ft$. 
So the assertion holds. QED

Write  $t$ in two different way:
$$t=t_p+t_z\in \ft_p\oplus \ft_z;\quad t=t_{wp}+t_z'\in \ft_{wp}\oplus \ft_z,$$  
how are $t_p, t_{wp}$ related?
\begin{lemma}
$wt_p=t_{wp}\in t_{wp}, \quad t_z'=
=t-wt +wt_z=t_p-wt_p+t_z\in \ft_z$.
\end{lemma}
{\it Pf:}
From the first equation, we have
$ wt=wt_p+wt_z$. Also  $ W_\mu=W(\lfg_mu)$ and $\ft_z=\lfg_\mu^\sss\cap\ft$,
we claim $wt-t\in t_z, \forall w\in W_\mu$. The claim holds easily for
reflections, hence it holds for all $w$.   Thus $$t=wt+(t-wt)=wt_p+(t-wt +wt_z)
\in \ft_{wp}\oplus \ft_z.$$ 
Hence $wt_p=t_{wp}$, $ t_z'=t-wt +wt_z$ by uniqueness of decomposition.
Another expression for $t_z'$ is 
$$t_z'=t-wt +wt_z=(t_p+t_z)-w(t_p+t_z)+wt_z=t_p-wt_p+t_z. \QED$$

Let $\{\gamma'\}$ be the weights of the $T_{w(p)}$-action on
$V_{w(p)}$, then 
 one has
$$\gamma'=w(\gamma).$$
Associated with $\gamma'$ is a weight $\tgamma'$ which is a character 
of  $T_{wp}^0\times T_z$, as shown in Prop. \ref{asswei}, defined by
$\tgamma'(t)=\gamma'(t_{w(p)}).$
Hence
\begin{equation}\label{gamma}
\begin{split}
\tgamma'(t)&=\gamma'(t_{w(p)}) =w(\gamma)(wt_p) =\gamma(t_p) =\tgamma(t).
\end{split}\end{equation}

 The first group of weights on $T_{(p,I,z)}X_N$, denoted by $\{\tgamma\}$ as  in Prop. \ref{asswei} 
are invariant under $W_\mu$,  see Fig 4.1. The group $W_\mu$ are generated by
reflections w.r.t. planes containing $\gamma$.

 Those weights $\{\tbeta'\}$ are defined by  $\tbeta'(t)=-\beta(t_{w(p)})$.
But $t_{w(p)}=w(t_p)$, therefore
\begin{equation}\label{beta}
\tbeta'(t)=-\beta(t_{w(p)})=-w(\beta)(t_p)=w(\tbeta)(t).
\end{equation}

 The weights $\{\tlambda'\}$, 
are defined  as 
\begin{equation}\label{lambda}
\tlambda'(t)=-\lambda(t'_z).
\end{equation}
The following transformation law for the weights is essential for future 
calculations. 

\begin{proposition}\label{transformwts}
Let $\spi:\ft\rightarrow \ft_z$ be the orthogonal projection, then
\begin{equation}
\tlambda(t)=\lambda(\spi  t)-\lambda(\spi t_p)
,\quad \tbeta(t)=\beta(\spi t_p).
\end{equation}

At the point $(wp,I,z)$ with $w\in W_\mu$, the weights are given by
\begin{equation}
\tgamma'=\tgamma,\quad \tlambda'(t)=\lambda(\spi   t)-\lambda(w\spi t_p),\quad
\tbeta'(t)=\beta(w\spi t_p).
\end{equation}

For $v\in W_\mu$, 
\begin{equation}
\begin{split}
&\tgamma'(vt)=\tgamma'(t)=\tgamma(t),\\
&\tlambda'(vt)=v\lambda(\spi t)- w\lambda(\spi t_p),\\
&\tbeta'(vt)=w\beta(\spi t_p).
\end{split}
\end{equation}
\end{proposition}
{\it Remark:} 1). The projection $\spi$ can be removed since
 $\beta,\lambda$ are linear functions on $\ft_z$ and  their
extension to $\ft$ by convention are compositions with $\spi$.
2). If $w\in W(\kk')$, then $wt_p=t_p$ since $[\kk',\ft_p]=0$,
and $\tbeta'=\tbeta$.


{\it Pf:} The proof is divided into 3 steps for each group of equations in the
above.

1). Since $t_z\in \ft_z$, and $t=t_z+t_p$, apply the map $\spi$, one has
$t_z=\spi t-\spi t_p$.  Hence
$$ \tlambda(t)=\lambda(t_z)=\lambda (\spi t)-\lambda(\spi t_p).$$ 
The root  $\beta$ vanishes on the orthogonal complement of  $\ft_z$, thus 
$$\tbeta(t)=\beta(t_p)=\beta(\spi t_p).$$

2). An useful observation can be made here:
\begin{equation}\label{useful}
\spi(w t)=w\spi(t)
\end{equation}
which is true for $w$ given by  reflection $r_\beta$ with respect to a root $\beta$ of
$(LG)_\mu$, for 
$$\spi(r_\beta t)=\spi( t_\beta-\beta(t)\beta^\vv)=\spi(t_\beta)-\beta(t)\beta^\vv$$
Since $\spi(\beta^\vv)=\beta^\vv$, also $\beta(t)=\beta(\spi
t)$. Combine the two one gets 
$$\spi(r_\beta t)=r_\beta(\spi t).$$
 Hence it  holds for any $w\in W_\mu$.

 It also has been shown that $\tgamma'=\tgamma$. As for $\tlambda'$ and $\tbeta'$, using $t'_z=t-t_{wp}=t-wt_p$, we obtain
\begin{displaymath}\begin{split} 
&\tlambda'(t)=\lambda(t_z')=\lambda(\spi t)-\lambda(\spi t_{w(p)}),\\
&\tbeta'(t)=\beta(t_{w(p)})=\beta(wt_p)=w\beta(t_p)=w\beta(\spi t_p).\\
\end{split}\end{displaymath}

3). Decompose $t=t_p+t_z=v^{-1}t_p+t_z+(t_p-v^{-1}t_p)$, we've shown that
$t_p-v^{-1}t\in \ft_z$, since $v\in W_\mu$. Therefore $t_z'=t_z+(t_p-v^{-1}t)\in \ft_z$ and $v(t)=t_p+v(t_z')$ with $v(t_z')\in \ft_z$, thus  $(vt)_p=t_p$
and 
 $$\tgamma'(vt)=\tgamma(vt)=\gamma((vt)_p)=\gamma(t_p)=\tgamma(t).$$

As shown in Step 2), by replacing $t$ there with $vt$,
$$\tlambda'(vt)=\lambda(\spi vt)-\lambda(w\spi (vt)_p),$$
but $\spi vt=v\spi t$, and it has just been shown $(vt)_p=t_p$, 
so we have 
\begin{equation}\begin{split}
\tlambda'(vt)&=\lambda(v\spi t)-\lambda(w\spi t_p)\\
&=v\lambda(\spi t)-w\lambda(\spi t_p).
\end{split}\end{equation}

For $\tbeta'$, following Step 2), one has
$$\tbeta'(vt)=w\beta(\spi (vt)_p),$$
but $(vt)_p=t_p$, therefore $$\tbeta'(vt)=w\beta(\spi t_p). \QED$$ 

Thus we have found how the weighs are related along $W_\mu$-orbit, and
how they change under $t\mapsto vt$. 

\subsection{Transformation of weights for $Z_h$}
For the $p\in Z_h$, the group $W(\KK)/W(\KK_h)$ acts on it, in addition to
 the transformation by $w$ in $ W((LG)_\mu)/W(\KK)$.
 The first group
preserve $Z$ but permutes among  $\{Z_h\}$. 
Since both $(LG)_\mu$ and $\KK$ contain $T$, there is the obvious exact
sequence 
$$W(\KK)/W(\KK_h)\rightarrow W((LG)_\mu)/W(\KK_h)\rightarrow W((LG)_\mu)/W(\KK).$$

 The strata $Z_{vh}/T_z\times \{z\}$ contains $[vp,I,z]$, where
$Z_{vh}=vZ_h$.  The normal
space to $Z_{vh}/T_z\times \{z\}$ acquires two more set of weights
$\{\theta^v_i\}$ and $\{\beta^v\}$ as shown earlier.

Let $t\in T$, with a fixed decomposition $t_zt_p$, the action on the normal
space at $[vp,I,z]$ is that of $vht_p\cdot (vh)^{-1}t_z$. 

The assertions in the next lemma are already verified in Lemma \ref{stratawts}
and Prop. \ref{transformwts}
\begin{lemma}
Let $v\in W(\KK)/W(\KK_h)$.
The weights $\{\ttheta^v_i\}$ and $\{\tbeta^v\}$ satisfy the following:
 $$\ttheta^v(vht_p\cdot (vh)^{-1}t_z)=\theta(h),\quad \tbeta^v(vht_p\cdot
(vh)^{-1}t_z)=\beta(vh).$$

If $v\in W((LG)_\mu)/W(\KK)$, $v$ does not preserves $T^0_p$. The lifting
of $t$  is given by $(vh)t_{vp}\cdot (vh)^{-1}t'_z$, where $t_{vp}=vt_p$.
With that lifting,
the weights $\{\ttheta^v\}, \{\tbeta^v\}$ are given by the same formula.

In both cases, the weights of the action by $((vh)t_{vp}, (vh)^{-1}t'_z)$ in the normal space of $vZ$ in $vF$ are simply  given by the evaluations of
the weights  given in Prop. \ref{transformwts}  on $((vh)t_{vp}, (vh)^{-1}t'_z)$
\end{lemma}

\subsection{A word about lifting of action by $t$}
For orbifolds, in order to evaluate the contribution by the fixed point set,
it is necessary to consider all the lifting of the action at a fixed point to the finite
smooth cover. Once there,   one needs to find the fixed point set of each
lifting and find its contribution. The strata $Z_h/T_z$ for $h\in
I_p\setminus I_p^0$ is one of those fixed point sets. What if  $h\in
I_p^0$? For such a $h$, $(ht_p,h^{-1}t_z)\in T_p\times T_z$  is a lifting of $t\in T$ as well, but $ht_p\in T_p^0$, and hence the consideration for that
lifting is already incorporated when we study the action by general
$(t_p,t_z)\in T_p^0\times T_z$.

\section{ Curvatures of various bundles}

In order to understand the contribution from the fixed points coming from 
compactification, i.e. those with images on $W(\partial C)$, we need to know 
the curvature of their normal bundles. And their transformation law  
by certain subgroups of $W^\aff$.
 
\subsection{A general fact}
Suppose $S$ acts on a manifold $N$ and $\C^n$, the action on $N$ is locally free and   the action on $\C^n$ is linear 
defined by $\lambda:\mathrm{Lie} S\rightarrow \fg l (\C^n) $. Let $A$ be a $S$-invariant  connection on $N$ with $\Lie S$ identified with the vertical
subspace in $TN$.
Denote the curvature of the connection by $F_A$ which is a $\Lie S$-values
horizontal two form.
\begin{proposition}\label{normalcurv}

1).  Then $(N\times \C^n)/S$ as vector bundle over $N/S$ has  a connection 
whose curvature is given by $ \lambda(F_A)$.

2). Suppose $V=\C^n/\Lambda$ where $\Lambda$ is a finite group acting on
the vector space linearly. If  $S$ acts on $N\times V$, such that the action
is locally free on $N$, and is a linear action on the orbifold line bundle $V$, i.e.
there is an extension of $S$ by the finite group $\Lambda$, $S'$ acting on
$\C^n$ linearly.  Then the orbifold line bundle $(N\times V)/S$ over $N/S$ has the same
curvature form as in 1). 
\end{proposition}
{\it Pf:} 1). On the trivial bundle $N\times \C^n$, defines the following
connection $d+\lambda\cdot A$, notice that $\lambda\cdot A:TN\rightarrow \fg
l(\C^n)$ is indeed  a  1-form. 
This connection has the feature that the curve $(g_t(p),g_t(v))$ is
horizontal for any 1-parameter subgroup $\{g_t\}$. Therefore, it descends to a connection on the quotient space
$N/S$. For the quotient connection constructed this way, the curvature is given by the descent of the curvature upstairs which is of the given form.

2). Let $S'$ act on $N$ through $S$. Then  apply the above argument to the extension group $S'$. QED  

\subsection{Action by $T_z$ on bundles}
The group $T_z$ has a local free action  on either  $Z$ or $Z_h$. The action
is obtained by restricting its action on $X$ to the given sets. Therefore it
acts on the various subbundles of the normal bundles.
Denote the action on normal bundle by $dt_z, \forall t_z\in T_z$.

\begin{lemma}

1). $dt_z: V_p\rightarrow V_{t_zp}$ is an isomorphism preserving the
decomposition into $V^0\oplus N$.

The  map $dt_z$ extends to an isomorphism 
\begin{equation}\label{normal}
dt_z:N_{p}\oplus (\lfg)_\tphi/\kk\oplus H_z \rightarrow
N_{t_zp}\oplus (\lfg)_\tphi/\kk\oplus H_z.
\end{equation}

2). The quotient by $T_z$ as described in  Eq.~(\ref{normal})  defines the normal orbifold line bundle to the fixed point set $F$ in
$X_N$.  
\end{lemma}
{\bf   Warning:}  The action by $t_z$  above should not be confused  with the
$T$-action on normal bundle to the fixed points studied in Section 3.
 The  action here defines the equivalent class, while  the action  
in Section 3  is on the set of equivalent classes.

{\it Pf:}  Let $(p_s,h_s,z_s)$ be a curve with $(p_0,h_0,z_0)=(p,I,z).$
Let $q=[I,z]\in \X$.
Using the defining equivalence relation,
 the diagonal action by $T_z$ on the first two factors is 
$$(p_s,h_s,z_s)\mapsto (t_zp_s,t_zh_s,z_s)\simeq (t_zp_s,t_zh_st^{-1}_z,t^{-1}_zz_s).$$
Differentiate the above at $s=0$ to obtain the diagonal action by $t_z$ on the
tangent
vectors:
\begin{equation}\label{5.2}
dt_z(p',h',z')\simeq(t_{z*}(p'),Ad_{t_z}h',t^{-1}_z(z')).\end{equation}
If $(p',h',z')\in V_p$, we have $z'=0$ and 
\begin{equation}
D_{p'+iJp'}\mu(p)=D_{h'+iJ'h'}k\mu_\X(q).
\end{equation}
Let $t_{z*}$ denote the induced action by $t_z$ on the tangent space.
Apply the action $\Ad_{t_z}$ on both sides, and combine that with the
properties of $d\mu$,  one has
\begin{eqnarray*}
D_{t_{z*}(p'+iJp')}\mu(t_zp)&=&\Ad_{t_z}(D_{(p'+iJp')}\mu(p))\\
&=&\Ad_{t_z}(D_{(h'+iJh')}\mu_\X(q))\\
&=&D_{t_{z*}(h'+iJh')}\mu_\X(q),
\end{eqnarray*}
which shows that $(t_{z*}(p'),t_{z*}(h'),0)\in V_{t_zp}$.

The calculation above also  shows for general $g\in LG$, the following holds:
\begin{equation}\label{Translaw}
g_*(p'+iJp',h'+iJh',0)=(g_*(p'+iJp'),\Ad_g(h'+iJh'),0).
\end{equation}

Since $T_p,T_z$ commute, one has $T_{t_zp}=T_p$, and clearly
$(p'+iJp',h'+iJh',0)$ has a non-zero weight under $T_p$ iff
$(t_{z*}(p'+iJp'), t_{z*}(h'+iJh'),0)$ has non-zero weight under
$T_{t_zp}=T_p$. Hence $N_p$ is $N_{t_zp}$ under $t_{z*}$.

The map $t_{z*}$ preserves $(\lfg)_\tphi/\kk$ since $\ft_z$ is the Cartan
subalgebra of the semi-simple part of $(\lfg)_\tphi$, and $T_z$ preserves $\kk$
under the adjoint action.
 The map $t_{z*}$ preserves $H_z$ because $T_z$ fixes $z$ and
$H_z$ is the subspace in $T_zX_\fg$ with non-zero weights. So it is stable
under the action by $T_z$.

2). The normal subspace to the fixed point set consists of non-zero weight
vectors in the   tangent space  $V_{p}\oplus (\lfg)_\tphi/\kk\oplus H_z$.
Therefore, it is given by $N_p\oplus (\lfg)_\tphi/\kk\oplus H_z$ before
mod out by the local free action of  $T_z$.
QED

\subsection{Curvatures of the normal bundles}
Next we calculate the curvatures of the three subbundles given in the
above. The transformation of the curvatures under the action by $W_\mu$ will
be given as well.

Following Eq. (\ref{5.2}), for $(p',h',z')$ in $(\lfg)_\tphi/\kk$, i.e.,
$p'=0,z'=0$, and $h'\in(\lfg)_\tphi/\kk$, the action by $T_z$ is
the adjoint action. And the action by $T_z$ on $H_z$ is the action by
isotropy group, since $T_z$ fixes $z$.

 Based on that Eq.~(5.2), one can write the curvature
in terms of the weights and the connection $A$ as follows:

\begin{proposition}\label{prop52}

1).The curvature of the bundles $$Z\times_{T_z}( (\lfg)_\tphi/\kk),
\quad Z\times_{T_z} H_z$$
are  given by $$-\oplus_{\beta\in \Delta_\mu/\Delta(\kk)} \beta\cdot dA, \quad -\oplus_{\lambda}
\lambda\cdot dA$$
respectively.  

2).
Let $\tnabla$ be a $T_z$-invariant connection on $\tN=\{N_p\}_{p\in Z}$ over
$\tF_p\cap\mu^{-1}(\phi)$, and $a_q(t)={\mathcal L}_t-\tnabla_t$, where
${\mathcal L}_t$ is the Lie derivative of $t$ acting on the bundle $\tN$, 
then the curvature of $N$ is given by
$$(\tnabla)^2+a_q(dA).$$
\end{proposition}

{\it Remark}: A simple but potentially important observation is that  $t_z$
acts
 in Eq.~(5.2) on the last component $z_s$ as $t^{-1}_z$. This introduces the
$-$ sign in  the next proposition.
Hence the weights are $\{-\lambda\}$ while the curvatures are given by
$\{\lda dA\}$. For the bundle defined by $\lfg_\mu/\ft$,
the weights are $\{-\beta\}=-\Delta(\lfg_\mu)_+$ while the curvatures are $-\{\beta dA\}$.

\subsection{Transformations of curvatures under $W_\mu$}

Recall that the fixed point sets $\tF_p\cap \mu^{-1}(\phi)$ transform
under the action by the subgroup in $W^\aff$, $W_\mu$, which is the Weyl group
of $(LG)_\mu$. We want to know first how the transformations act on the 
normal bundles and how the curvatures change.

If the following, for $w\in W_\mu$, we fix a lifting to $(LG)_\mu$ and it
will be denoted the same.

It follows directly from Eq. (\ref{Translaw}), at $(wp,I,z)$, the bundle
$V_{wp}=w^*V_p$.  Also, $T_{wp}=wT_z$, thus $v\in V_p$ is of non-zero weight
with respect to $T_p$ iff $w^*(v)$ is of the same weight with respect to
the action by $wT_p=T_{wp}$. Therefore, one has $N_{wp}=w^*N_p$. On the other
hand, $T_z$ is the Cartan subgroup of $(LG)^{\sss}_\mu$ which is the semi-simple
part of $(LG)_\mu$, whence $W_\mu$ is the normalizer of $T_z$. 
And conveniently we have $\{N_{wp}\}/T_{wp}=w(\{N_{p}\}/T_{p}).$
Hence we have proved  the first part of the following:
\begin{proposition}\label{transformcur}
1). The  pull-back to $\{N_p\}$ by $w$  of curvature form $(\tnabla)^2+a_q(dA)$  of the bundle $\{N_{wp}\}$ at
$wp$ is the same as that of $\{N_p\}$ at $p$ under the map $w^*$. 

2). The pull-back under $w$ of the  bundles $wZ\times _{T_z}(\lfg)_\mu/\kk$,
$wZ\times_{T_z} H_z$ have 
curvatures given respectively  by   
$$-\oplus_{\beta\in \Delta_\mu/\Delta(\kk)} w^{-1}(\beta)\cdot dA, \quad  \oplus_{\lambda}
w^{-1}(\lambda)\cdot dA.$$
\end{proposition}
{\it Pf: } The two expressions in Part 2 can be verified the same way.
Let's take the first one. For each subalgbera   $\lfg_\beta$,
one has the isomorphism:
$$w^*:\tF_{p}\cap \mu^{-1}(\phi))\times \Ad_{w^{-1}}(\fg_\beta)\rightarrow
\tF_{wp}\cap \mu^{-1}(\phi))\times \fg_\beta.$$
The action by $T_z$ on $\Ad_{w^{-1}}(\fg_\beta)$ has weight
$-w^{-1}(\beta)$. Therefore the curvature is given by
$-w^{-1}(\beta)\cdot dA$. QED

\subsection{More on the curvature of the bundle $N$}
Recall $\tN=\{N_p\}_{p\in Z}$,  and $N$ is the quotient $\tN/T_z$.
The group $\KK'$ acts on $\tN$ in a obvious manner, since $\KK'$ commutes
with $T_p$.
Let $\tnabla$ be a $\KK'$-invariant connection on the bundle $\tN$, and the
associated moment map $\epsilon$ be defined as $$<\epsilon(q),(\xi)>:={\mathcal
L}_\xi-\tnabla_\xi\in \hom(N_q).$$
Since we have $PA:TZ\rightarrow \kk'$, we can change $\nabla$ to
a new connection by adding a 1-form $$\nabla=\tnabla+<\epsilon(q),
PA(\cdot)>.$$
The invariant sections of $\tN$ now are horizontal with respect to $\nabla$.

Using invariant sections on $\tN$, one proves easily the following:
\begin{proposition}  The connection $\nabla$ descends 
  respectively to
connections on
$N\rightarrow F$ and  $\tN/\KK'\rightarrow Z/\KK'=E$. The curvature
2-form on $F$ is the pull-back of that over $E$. In particular the curvature is
trivial on the fiber $\pi:F\rightarrow E$.
\end{proposition}

\subsection{The curvature form  of bundles over $Z_h/T_z$}
Along the strata $Z_h/T_z$, there is the normal bundle in $Z$:
$\nno(Z_h,Z)\oplus \kk/\kk_h$ in addition to the restriction to $Z_h/T$  of the normal bundle of $Z/T_z\times \{z\}$ in $TX_N|_{Z_h/T_z\times \{z\}}$.

The curvatures and their transformations under $W_\mu$ of  those two
bundles  can be written down similar to what we have just done for $Z$.
Recall that $Z_h$ admits a locally free action by $\KK'_h$, and $Z_h/T_z$ is
a
fiber bundle with fiber given by a finite quotient of
$\KK'_h/T_z$.
On $Z_h$ fix a $\KK'_h$-invariant connection $A_h$, so that the tangent space
is
decomposed as $T^HZ_h\oplus T^\perp Z_h$, and $P_h:T^\perp Z_h\rightarrow
Z_h\times \kk'_h$ is the trivialization of the vertical bundle. The map $P_h$
performs the same function as that of $P$ for $Z$.
Fix a connection for the bundle $T^\perp Z$ which satisfies
$\nabla_\xi={\mathcal L}_\xi, \forall \xi\in\kk'_h$.
 Its existence follows the discussion on $Z$ with respect to
the group $\KK'$.

For convenience, we group the two components in the normal space to $Z_h$,
$(\lfg)_\mu/\kk$ and $\kk/\kk_h$ together as $(\lfg)_\mu/\kk_h$. They received
separate treatment earlier since the first one is  normal to $Z$
and the second one is tangent to $Z$ but normal to $Z_h$.

The proof of   the following is identical to that of
Prop. \ref{transformcur}, just replace $Z_h$ by $Z$.
\begin{proposition}
The curvature of $\nabla$ is the pull-back of  a form on $Z_h/\KK_h$.

The curvature of the bundle $Z_h\times _{T_z}((\lfg)_\mu/\kk_h)$ on $Z_h/T_z$
is
given by $-\sum_{\epsilon}\epsilon \cdot dP_hA_h$
where $\epsilon$ are the simple  roots of $(\lfg)_\mu$ not in $\kk_h$.
Furthermore, the curvature form $R_h=dP_hA_h $ can be decomposed into
$$dP_hA_h= B_h+R_h$$
where  $B_h(u,v)(pg)$  and $R_h(\xi,\eta)$ satisfy the same property as that
of $B, R$ in Lemma \ref{damn}. 

For $w\in W_\mu$, the
normal bundle at $wp\in Z_{wh}$ is given by
$$ wV_p^{01}\oplus wN_p\oplus (\lfg)_\mu/\kk_{wh}\oplus H_z,$$
where $V_p^{01}$ is the maximal $h$-stable subspace of $V_p^0$ on which
$\det(I-h)\neq 0$. 
The curvature of the first two components are the same as those at $p\in Z_h$.
The curvature of the last two components are given respectively by
$$- \oplus w^{-1}(\beta)\cdot dP_hA_h,\quad \oplus w^{-1}(\lambda)\cdot dP_hA_h.$$
\end{proposition}

\section{A fundamental formula}
In order to prove the desired fixed point formula, we need to understand  the
contribution from the  fixed points in $X_N$ of type 2) and 3). Those do not
come as fixed points of $X$ itself under the $T$-action, rather from certain compactification,
$X_N$, of the quotient  $X$ by the nilpotent subgroup  $LG^{\C+}$. For in the space
$X_N$, there is an open dense set which corresponds  to $X/LG^+$.  

We will prove in this section an important formula which enables us to calculate the
contribution from the type 2) and 3) fixed point set in the next section.

Suppose $\fk$ is the Lie algebra of the   semi-simple compact Lie group $K$,
with $\ft$ as its Cartan subalgebra, given a choice of Weyl chamber in $\ft$ of $\fk$,
 assume that $\{\lambda\}$ is the set of fundamental weights in $\ft^*$,
 $\{\alpha\}$ is the set positive roots, and $W$ is the Weyl group. 
Let $k$ denote the rank of $\fk$.

For a character of $\ft$, $l$, define the symbolic notation
\begin{equation}\label{z}
 z^l=\exp(2\pi i l(x)), \quad z_1^l=\exp(2\pi i l(y))
\end{equation}
where $x, y$ are $\ft$-valued forms of even order, e.g., $x=x'+x''$ with
$x'\in \ft$ and $x''$ a $\ft$-valued curvature two form. 
Denote  the number of simple roots by $m$, and the set of positive roots
by $\Delta^+(K)$, the set of fundamental weights by $\{\lambda\}$.

\begin{proposition}[Fundamental Formula]\label{fund}

\begin{equation}\label{pride}
\begin{split}
&\frac{1}
{ z^\rho  z_1^\rho \prod_{\alpha\in\Delta^+(K)}
(1-z^{-\alpha})(1- z_1^{-\alpha})}
\sum_{w,v\in W}
\frac{(-1)^{m+\sigma(w)+\sigma(v)}}
{ \prod_\lambda(1-z^{-w\lambda}z_1^{-v\lambda })}\\
&=\sum_{w,v\in W}
\frac{1}{\prod_\alpha(1-z^{-w\alpha})\prod_\lambda(1-z^{w\lambda}z_1^{v\lambda
})\prod_\alpha(1-z_1^{-v\alpha})}\\
&=0.
\end{split}
\end{equation}

\end{proposition}
{\it Pf:} Step 1): 
Apply the well known formula that
$$\sum_{w(\alpha)<0}w(\alpha)=w(\rho)-\rho$$
where $\rho$ is the half sum of the positive roots, 
to obtain  \begin{equation}\label{long}\begin{split}
&\frac{1}{\prod_\alpha(1-z^{-w\alpha})\prod_{\lambda}(1-z^{w\lambda}z_1^{v\lambda
})\prod_{\alpha}(1-z_1^{-v\alpha})}\\
&=\frac{(-1)^{\sigma(w)+\sigma(v)}}{z^{(-w(\rho)+\rho)}z_1^{(-v(\rho)+\rho)}\prod_\alpha(1-z^{-\alpha})\prod_{\lambda}(1-z^{w\lambda}z_1^{v\lambda
})\prod_{\alpha}(1-z_1^{-\alpha})}.
\end{split}
\end{equation}
Next observe that $\sum_{\lambda}\lambda=\rho$, or $\sum w(\lambda)=w(\rho)$. Likewise
for $v$, thus 
\begin{equation}\begin{split}
\prod_{\lambda}(1-z^{w\lambda}z_1^{v\lambda})
=(-1)^mz^{w\rho}z_1^{v\rho}\prod_{\lambda}(1-z^{-w\lambda}z_1^{-v\lambda}).
\end{split}
\end{equation}
Denote the left side of the equation (\ref{pride}) by $L.H.$,
after substituting the last expression into (\ref{long}) to get  
\begin{equation}
L.H.= \frac{1}
{z^\rho z_1^\rho\prod_\alpha(1-z^{-\alpha})\prod_\alpha(1-z_1^{-\alpha})}
\sum_{w,v\in W}(-1)^{m+\sigma(w)+\sigma(v)}
\frac{1}{\prod_\lambda(1-z^{-w\lambda}z_1^{-v\lambda })}.
\end{equation}

First we rename the index $v$ to $vw$, this is legitimate for obvious reason,
 then
\begin{equation}
\begin{split}
&\sum_{w,v\in W}(-1)^{\sigma(w)+\sigma(v)}\frac{1}
{\prod_\lambda(1-z^{w\lambda}z_1^{v\lambda })} \\
&=\sum_{w,v\in W}(-1)^{\sigma(w)+\sigma(vw)}
\frac{1}{\prod_{\lambda}(1-z^{w\lambda}z_1^{vw\lambda})}\\
&=\sum_{w,v\in W}(-1)^{\sigma(v)}
\frac{1}{\prod_{\lambda}(1-z^{w\lambda}z_1^{vw\lambda})}.
\end{split}
\end{equation}

Step 2): 
Next we shall verify that the last  sum  vanishes. 

It is well known that $\{\lambda\}$ spans the positive Weyl cone $\ft^+$,
and the set of cones of the form $w(\ft^+), \forall w\in W$ spans $\ft$, and
there is no overlapping in the interiors of those cones.

Let $\calS$ denote the set of all possible subsets of $\{\lambda\}$,
and $\forall S\in \calS$, let  $|S|=\#S$,  it satisfies $1\leq |S|\leq m$.
The maximum number is  $m$.

From [M], one obtains the following relation:
\begin{equation}\label{fund0}
\sum_{S\in\calS,w\in W}(-1)^{|S|}\frac{1}
{\prod_{\lambda\in S}\Big(1-\exp<{w\lambda},\lp x>\Big)}=1.
\end{equation}
From this we deduce that
\begin{equation}
\begin{split}
&(-1)^m\sum_{w\in W}\frac{1}{\prod_{\lambda
}\Big(1-\exp<{w\lambda},\lp x>\Big)}\\
&=1-\sum_{|S|<k,w\in W}(-1)^{|S|}\frac{1}{\prod_{\lambda\in
S}\Big(1-\exp<{w\lambda},\lp x>\Big)}.
\end{split}
\end{equation}
Let $x=t+v(s)$ with $z=e^{\lp t}, z_1=e^{\lp s}$, we have  
$$ \exp<w\lambda,\lp x>=z^{w\lambda}z_1^{vw\lambda}.$$
Summing over $v\in W$ with sign $(-1)^{\sigma(v)}$, we have
\begin{equation}
\begin{split}
\sum_{w,v\in W} \frac{(-1)^{k+\sigma(v)}}{\prod_{\lambda
}(1-z^{w\lambda}z_1^{vw\lambda})}
=\sum_{v\in W}(-1)^{\sigma(v)}-\sum_{ w,v\in W, |S|<k}
\frac{ (-1)^{|S|+\sigma(v)}}
{\prod_{\lambda\in S}(1-z^{w\lambda}z_1^{vw\lambda})}, \end{split}
\end{equation}
the first term in the last line is 0. When $|S|<k$, those weights in $S$
span a face of the Weyl chamber $\ft^+$. The face, or $S$, is fixed by a  non-trivial subgroup of $W$. The subgroup  is denoted by $W_S$ which is the Weyl
group of a subgroup of $K$.  And the subgroup $\Ad_wW_S$ fixes the set of
weights $w(S)$.

Let $W$ be divided into cosets by $\Ad_wW_S$ and let $[v_0]$ denote the coset, 
then for a fixed $S$ with $|S|<k$, one has
\begin{equation}
\begin{split}
&\sum_{w,v\in W}
\frac{ (-1)^{\sigma(v)}}{\prod_{\lambda\in S}\Big(1-z^{w\lambda}z_1^{vw\lambda}\Big)}\\
&=\sum_{[v_0] \in W/\Ad_w W_S,w\in W}
\sum_{v\in \Ad_wW_S}
\frac{ (-1)^{\sigma(v)}}{\prod_{\lambda\in S}\Big(1-z^{w\lambda}z_1^{vw\lambda}\Big)}\\
&=\sum_{[v_0] \in W/\Ad_w W_S, w\in W}
\sum_{v\in \Ad_wW_S}
\frac{(-1)^{\sigma(v)}}{\prod_{\lambda\in S}\Big(1-z^{w\lambda}z_1^{v_0w\lambda}\Big)}\\
&=0
\end{split}
\end{equation}
because $\sum_{v\in
\Ad_wW_S}(-1)^{\sigma(v)}=0$.
Thus we have the desired claim. QED

\section{Affine Weyl subgroups and the dual Coxeter numbers}
In order to apply the formula proved in the last section, 
we need to know more about those fixed points of type 2) and 3) occurring
along an affine wall of $C$. 

Some familiarity with the affine Kac-Moody algebra is required here, see [K]
for reference.  

 Let $(p,I,z)$ be a point which defines a fixed point in  $X_N$. 
Assume that $\mu(p)/k=\tphi(z)\in (\partial C,1)$, i.e.
that fixed point is of type 2) and 3). Let $\mu$ denote $\mu(p)$.

\noindent {\bf Conventions:}

 In the following, we treat the case when $\mu$ is of level 1 to simplify the
notation. The general case follows after replacing $\mu$ by $\mu/k$.

If $\mu\in \partial C$,  there are two possibilities:

a). $\mu\in C\setminus C^\aff$, i.e., $\mu$  is not on the wall defined by $\theta=1$, 

b). $\theta(\mu)=1$, or $\mu\in C^\aff$.  
 
In the first case, $W_\mu$ is generated by reflections of a subset of simple
roots, and $(LG)_\mu$ is generated by all the $\fg_{\beta}$ such that   $\beta$ is a 
simple root vanishing at $\mu$. In particular, $(LG)_\mu$ is a subgroup of
$G$ and $W_\mu$ is a subgroup of the regular Weyl group $W$.

In case of b), the affine root $\alpha_0=\delta-\theta$ plays an important role. And the reflection
with respect to $\alpha_0=0$ or $\theta=1$ is the composition of 
the reflection defined by $\theta=0$ and a translation element. 

To understand
this and the group $(LG)_\mu$ more thoroughly, we need a few things from the theory of affine Lie
algebra. 

\subsection{A few facts on affine Lie algebra}

\begin{proposition} Let $\Delta^+$ be the set of simple roots of $\fg$ with
respect to the cone spanned by the alcove $C$, $\{\Lambda_i\}$ be the set of 
fundamental weights of $\fg$.

a). If $\mu$ is not on the affine wall, then 
$$\Delta^+_\mu=\{\beta|\beta (\mu)=0, \beta \in \Delta^+\}$$
 is a set of simple roots for
$(\lfg)_\mu$. The fundamental weights of $(\lfg)_\mu$ is the orthogonal
projection with respect to the form $(\cdot |\cdot)$ of $\{\Lambda_\beta\}$ 
to the linear span of $\Delta^+_\mu$.

b). If $\mu$ is on the affine wall defined by $\theta=1$, 
let $$\Delta^0_\mu=\{\beta|\beta (\mu)=0, \beta\in \Delta^+\},$$
$$\Delta^+_\mu=\{\delta-\theta\}\cup \Delta_\mu^0$$
is the set of simple roots.
The fundamental weights are given by the orthogonal projection of 
$$-\mu,\quad \Lambda_\beta-a_i^\vv\mu, \, \beta\in \Delta^+_\mu$$ 
to the linear span of $\Delta^+_\mu$.

In both cases, the group  $(LG)_\mu$ is connected and its Weyl group is the 
subgroup in $W^\aff$ generated by reflections using elements in
$\Delta^+_\mu$.
 \end{proposition}
{\it Remark:} The fundamental weights of $(\lfg)_\mu$ may not be weights of
$\fg$. Nevertheless, they are in the rational span of $\Delta^+$.

{\it Pf:} For  a), the  Lie algebra $(\lfg)_\mu$ is a subalgebra of $\fg$,
generated by $\fg_\beta$, $\beta\in \Delta^+_\mu$. The assertions  are  well known facts  about finite dimensional semi-simple algebras.

For part b), first we verify the assertion about the fundamental weights.
\begin{equation}\begin{split}
&2(\Lambda_\beta-a^\vv_\beta\mu|\alpha)/(\alpha|\alpha)=
2(\Lambda_\beta|\alpha)/(\alpha|\alpha)=\delta_{\beta\alpha},\quad \alpha\in
\Delta^0_\mu;\\
&(-\mu|\alpha)=0,\quad \alpha\in\Delta^0_\mu;\\
&(-\mu|-\theta)=1;\\
&(\Lambda_\beta-a^\vv_\beta\mu|\theta)=a^\vv_\beta-a^\vv_\beta=0,
\end{split}
\end{equation}
the last equality is due to the fact that $(\Lambda_\beta|\theta)=a^\vv_\beta$.
Thus, the given set is dual to $\Delta^+_\mu$, and therefore its orthogonal projection is the set of fundamental weights of $(\lfg)_\mu$.

The best way to see $\Delta^+_\mu$ of  part b) is a set of simple roots is
to use the theory of
affine Lie algebras, from which we learn that $$
\Delta_\mu=\{\alpha_0:=\delta-\theta\}\cup \Delta^+$$ form a set of simple roots for the 
algebra $\affg$. Now it is known that any subset of simple roots is a
set of simple roots for the subalgebra generated by the subset. Therefore,
$$\Delta_\mu^+=\{\alpha_0\}\cup  \Delta^0_\mu$$  is a  set of simple roots. 
On the other hand, $$(\delta|\delta)=0,\quad (\delta| \alpha)=0,  \ \forall \alpha\in \Delta,$$ see
[K,Ch.6]. Thus the inner product of a pair in $\Delta'$ is the same as that of the
corresponding pair in $\Delta^+_\mu$. In particular,     the Dynkin diagrams formed by $\Delta'$ and $\Delta_\mu^+$ are the
same. Therefore, $\Delta_\mu$ form a set of simple roots for the subalgebra
$(\lfg)_\mu$.

From the characterization of the simple roots of $(LG)_\mu$, it is clear that
its Weyl group is generated by reflections using $\Delta_\mu^+$.
 In case $\mu$ is
 on an affine wall, the reflections are with respect to
$\{\alpha=0|\alpha\in\Delta^0_\mu\}$ and $\theta=1$, as desired. 

The connectedness is based on the well known argument in Lie theory. Decompose
$(LG)_\mu=\cup  K_i$ into connected component, and  $K_0$ contains 
$I$. For $g\in K_i$, $\Ad_{g}K_0=K_0$. Multiplying $g$ with an
element in $K_0$ if necessary, we can assume that $\Ad_{g}T=T$ where $T$ is
the maximal torus. Therefore, $g$ is in the Weyl group of $K_0$, therefore
$g\in K_0$. QED

\subsection{The half sum of positive roots of $(LG)_\mu$}
Let $\Delta_\mu$ be as before, and $\{\lambda_i\}$ be the  set of
fundamental weights of $(\lfg)_\mu$. We have seen that $\lambda_i$ is
given by the orthogonal projection of $\Lambda_\beta-a^\vv_\beta\mu$.  

Let $\rho_\mu$ be the half sum of
positive roots of $(\lfg)_\mu$.
 For finite dimensional  Lie algebra, it is well known
that
\begin{equation}
\rho_\mu=\sum_i\lambda_i.
\end{equation}

The following is as important as the fundamental formula:

\begin{proposition}\label{pride2}
Let $\rho$ be the half sum of positive roots of $\fg$, $\rho_\aff=\rho+h^\vv
\Lambda_0$, where $h^\vv$ is the dual Coxeter number defined by $h^\vv-1=\sum_{i=1,...,l}a^\vv_i$,
then the following holds:

1). If $<\theta,\phi><1$:
$$w(\rho_\mu)-\rho_\mu=w(\rho)-\rho \  \mod \Z \delta, \quad  \forall w\in W_\mu$$
and 
\begin{equation}\label{half}
e^{2\pi i(<vk\phi -\rho+v\rho-v\rho_\mu,t>)}=e^{2\pi
i<k\phi+\rho_\mu,t>}.
\end{equation}

2). If $<\theta,\phi>=1$ is on the affine wall of $C$,
$$r(\rho_\mu)-\rho_\mu=r(\rho)-\rho \ \mod \Z\delta  $$
where $r$ is a reflection defined by a simple root of $\fg$.
Let $r_{\theta}$
be  the reflection with respect to the affine wall $\phi(\theta^\vv)=1$, 
$$r_{\theta}(\rho_\mu)-\rho_\mu=r_\theta(\rho)-\rho+h^\vv \theta,$$
where $\nu:\affg\rightarrow \affgd$ is the map induced by the bilinear form
$(\cdot|\cdot)$ on $\affg$.

And when $t$ is restricted to the lattice $\frac{M^*}{k+h^\vv}$, with $M^*$
being the dual of the long root lattice,
$$e^{2\pi i(<vk\phi +\rho-v\rho+v\rho_\mu,t>)}=e^{2\pi
i<k\phi+\rho_\mu,t>}, \forall v\in W_\mu^0$$
where $W_\mu^0$ is the subgroup of $W$ generated by the reflections
with respect to $\alpha\in \Delta^0_\mu$ and $\theta$. It is isomorphic
to $W_\mu$. 
\end{proposition}
{\it Pf:}
1). It is enough to verify that for reflection $r_\alpha$,
$\alpha\in\Delta^+_\mu$. In this case, $\alpha$ is also a simple root of $\fg$.
Thus, $$r(\rho_\mu)-\rho_\mu=-\alpha=r(\rho)-\rho.$$
In this case,  all the simple roots in $\Delta^+_\mu$ vanish at  $\tphi$ by
definition of $(\lfg)_\mu$. Therefore, $v\phi=\phi$, the eq. (\ref{half})
holds as a consequence of the identity just proved.

2). If $r$ is generated by a simple root of $\fg$, which is the case if
$\alpha\in \Delta^0_\mu,$ then the previous argument works.
If $r=r_\theta$, $-\theta$ is a simple root of $(\lfg)_\mu$ but not of $\fg$,
we obtain
 $$r_\theta(\rho_\mu)-\rho_\mu=\theta  \quad  \mod \Z\delta .$$
On the other hand $$r_\theta(\rho)-\rho=-<\rho,\theta^\vv>\theta=-(h^\vv-1)
\theta,$$
where the last equation is from the definition of $h^\vv$:
 $$h^\vv=1+\sum_i a_i^\vv=1+ <\rho,\theta^\vv>.$$
Now $$r_\theta k\phi=k\phi-k<\phi,\theta^\vv>\theta=k\phi-k\theta,$$
thus we obtain
\begin{equation}
\begin{split}
&e^{2\pi i<r_\theta k\phi -\rho+r_\theta\rho-r_\theta\rho_\mu,t>}\\
&\ =e^{2\pi i<k\phi -(k +h^\vv)\theta-\rho_\mu,t>}\\
\end{split}\end{equation}
which equals $e^{2\pi i<k\phi -\rho_\mu,t>}$ on the lattice since 
$$<(k +h^\vv)\theta,t>\in \Z,\quad \forall   t \in \latk.$$
The above holds now for the  generators in $W_\mu^0$, therefore it holds
on the lattice $\latk$ for all $W_\mu^0$.

The isomorphism between $W_\mu$ and $W_\mu^0$, the only difference among
the generators is the form has a reflection w.r.t. $\theta=1$ while the latter
has one w.r.t. $\theta=0$. 
QED

\section{Further orbifold complications } 

In order to apply fixed point principle to the space $X_N$, we  need  better
understanding of the action by $T$ on the normal bundles to the fixed point set
of all three types. 

\subsection{Weights on the toric variety $X_\fg$}

 We divide the discussion on the normal bundles according to the types of the
fixed point sets.

We continue to use the convention from the last section.
\subsection{When $\mu$ is on  $(\partial C,1)$}

This is the more interesting case. Recall that $X_\fg$ is constructed as a
global orbifold toric variety. First we investigate what is the stabilizer
and the weights at a point on $X_\fg$. 

From the last section, we have
learned that the fundamental weights of $(LG)_\mu$ is given  by orthogonal
projections  of the following vectors: 
$$\Lambda_\beta-a_\beta^\vv\mu,\quad \beta\in \Delta_\mu$$
together with $\-\mu$ if $\alpha_0=\delta-\theta\in \Delta_\mu$. 
 In the above, $\Lambda_\beta$ is a fundamental weight of $\fg$ as well since
$\beta$ is a simple root of $\fg$. Because $X_\fg$ is an orbifold, the polytope
$C$ is not a simple convex polytope (see [O] for definition)  with respect to the weight lattice $M$ 
generated by $\{\Lambda_\beta\}$,   but rather it is a simple polytope in 
the larger lattice  $M'$ generated by $\{\Lambda_\beta/a^\vv_\beta\}$. The
larger  
lattice defines a unique covering of $T$, $T'$ so that $\pi T'\rightarrow T$
has the quotient $M'/M$  as its kernel.

The dual lattice of $M$ is given by the coroot lattice  $$N=\sum\Z\alpha^\vv_i.$$ The dual
lattice of $M'$ is given by the  sublattice  $$ N'=\sum \Z(a_i^\vv\alpha^\vv_i).$$
So $T=\R^l/N$ and $T'=\R^l/N'$. 

The  polytope $C$ is integral with respect to $M'$, and it  is actually  a simple
simplex. Thus it  defines  smooth toric variety $X'$ with respect to the group $(T')^\C$, $X'$ is in fact the projective
space $\C P^l$, and $X_\fg=X'/\ker \pi$, see [Od, p96].

\begin{proposition}\label{orbiwts} 
On $X'\simeq \C P^l$,  assume  $z'\in \phi^{-1}(\partial C)$. 
The stabilizer of $z'$, $T'_z$ and the  weights of the action by $T'_z$ on the normal bundle is given by
 the following:

1). If $\mu=\phi(z')$ is not on the affine wall, 
the stabilizer is given by
$$\sum _{\beta\in\Delta_\mu} \R\beta^\vv/(\sum
_{\beta\in\Delta_\mu} \Z a_\beta^\vv\beta^\vv).$$
 The weights are given by the orthogonal projection to $\nu(\ft_z')$ of $$ 
\Lambda_\beta/a_\beta^\vv, \forall \beta\in \Delta_\mu.$$  

2). If $\mu$ is on the affine wall, recall $\Delta_\mu=\{-\theta\}\cup
\Delta^0_\mu$ and $\Delta_\mu^0=\{\beta\in \Delta_+ | \beta(\mu)=0\}$. 
The coroot $\theta^\vv\in N'$ by  definition is $\sum_ia^\vv_i\alpha_i^\vv$. 

In particular $\sum_\beta\Z a^\vv_\beta\beta^\vv$ is a sublattice of $N'$, where
$a^\vv_\theta=1$.
The stabilizer again is given by $$\sum _{\beta\in\Delta_\mu} \R\beta^\vv/(\sum
_{\beta\in\Delta_\mu} \Z a^\vv_\beta\beta^\vv)\subset T'.$$

 The weights of the stabilizer  are the orthogonal projection  to $\nu(\ft_z')$ of  $$\{-\mu\}\cup 
\{\Lambda_\beta/a_\beta^\vv-\mu\big|\beta\in \Delta_\mu^0\}.$$
\end{proposition}

{\it Pf:} The simplex $C$ defines a polarization and a moment map $\phi$ which
has $C$ as its image. Let $\partial C_\mu$ be the smallest face passing $\mu$
and $V_\mu$ be the smallest linear space containing $\partial C_\mu-\mu$.  If $\phi(z')=\mu\in\partial C$, the image of $d\phi(z')$
is exactly $V_\mu$,  and the subspace perpendicular to $V_\mu\subset \ft^*$ is the image under $\nu$ of the stabilizer $\ft'_z$.

On the other hand, the face $\partial C_\mu$ is defined by $\cap _{\beta\in
\Delta_\mu}\beta^{-1}(0)$, if $\mu$ is not on the affine wall. And given by
$$ \theta ^{-1}(1)\cap  _{\beta\in
\Delta^0_\mu}\beta^{-1}(0),$$
if $\mu$ is on the affine wall. In either cases, 
$$V_\mu^\perp=\nu(\sum_{\beta\in\Delta_\mu}\R \beta^\vv).$$
Therefore, the Lie algebra of the stabilizer $\ft'_z$ is of the desired form.
Clearly the lattice $\sum
_{\beta\in\Delta_\mu} \Z a^\vv_\beta\beta^\vv$ is in $\ft'_z$, and in fact it is 
$\ft'_z\cap N'$. Thus the stabilizer of $z'$ in $T'_z$  is  
 $\ft_z'/(\ft'_z\cap N')$ whose explicit form is given by the proposition.

To understand the claim on the weights of the action by the stabilizer on the
normal bundle, we recall first that each point on the toric variety $X'$,  the
neighborhood is constructed as follows:

 Let $A$ be the tangent cone of the
simplex $C$ at the point $\phi(z')$, $A$ is a convex cone. Take the semi-group
$\sigma_\mu=N'\cap A$. Because $C$ is a simple simplex, it is easy to see that
$$\sigma_\mu=\sum_i\Z_{\geq 0}\eta_i+\sum_j\Z \xi_j,$$ where  $\{\eta_i\}\cup
\{\xi_j\}$  is a base of the
lattice $N'$. Since $\sum_i\Z\xi_i$ is a sublattice, it is given by
the lattice points in the maximal linear subspace contained in $A$, $V_\mu$.
Then the action of $T'$ on a  neighborhood of $z'$ in the toric variety $X'$ is given by
$$t(z_1,...,z_m,w_1,...,w_n)=(t^{\eta_1}z_1,...,t^{\eta_m }z_m,t^{\xi_n}w_1,...,t^{\xi_n}w_n),$$
the point $z'$ correspond to a point with $w_j=0,z_i\neq 0$   which also
defines the fixed point set of the subgroup $T'_z$. 
The above facts about toric varieties can be found in Ch. 1.2 and Ch. 2.4 in [Od].
From this it is easy to read
off the weights by the action of $T_z'$ near $w_i=0$. They are given by the
restriction of those weights $\xi_i$ to $T'_z$. Or the projections to
$\nu^{-1}\ft_z'$ of $\{\xi_j\}$. 

What are those weights $\{\eta_i\}\cup \{\xi_j\}$?  First select a lattice point
$\Lambda_i/a_i^\vv$  on the
affine subspace spanned by the smallest face $\partial C_\mu$ passing $\mu$, 
then $$\{\eta_k\} =\{\Lambda_k/a_k^\vv-\Lambda_i/a_i^\vv\in V_\mu\},\quad
\{\xi_j\}=\{\Lambda_j/a_j^\vv-\Lambda_i/a_i^\vv\notin V_\mu\}.$$
On the other hand, if one replace $\Lambda_i/a_i^\vv$ by a point on $\partial
C_\mu$, such as $\mu$ itself, the projection to the orthogonal complement of $V_\mu$ does
not change. Therefore, the weights of the action by $T'_z$ can also be given
by the orthogonal projections of $$\{\Lambda_j/a_j^\vv-\mu\big|\Lambda_j/a_j^\vv-
\mu\notin V_\mu\}.$$ 
In the above construction,  $\{\Lambda_j/a_j^\vv=0$ is allowed to account for the
weight which is the projection of $-\mu$. 
QED

\begin{definition}\label{defofn}
Let $\epsilon^\vv$ be an element of the coroot lattice so that it is given by
$ (1/n)\sum_{ \alpha_i\notin \Delta_\mu^0}a_i^\vv \alpha^\vv_i$, where $n\geq
1$ and no fraction of $\epsilon$ is in the coroot lattice. 
\end{definition}
Obviously, for 
$su(l+1)$, there is only one choice  $n=1$.   

Notice the definition of $n$ depends on $\Delta^0_\mu$.
\begin{corollary}\label{toriciso}

1). When $<\phi,\theta>< 1$ where $\phi=\phi(z)\in\ft^*$, the group
$T_z'/T_z$ is given by the finite group  $$\exp(\sum_{\alpha_i\in \Delta_\mu 
}b_i\alpha_i^\vv),\quad 0\leq b_i< a_i^\vv.$$

2). When $<\phi,\theta>=1$ where $\phi=\phi(z)$,  the subgroup $T'_z/T_z$ is given by $$\exp(b\epsilon^\vv+\sum_{ \alpha_i\in \Delta^0_\mu }b_i\alpha_i^\vv),\quad 0\leq b< n, \quad 0\leq b_i< a_i^\vv .$$

\end{corollary}
{\it Pf:} The group in question satisfies  $T'_z/T_z\simeq N/N'$. From the expression of $N, N'$, we
can identify the elements in the kernel of the map $T'_z\rightarrow T_z$
easily.

 The second assertion is based on the simple observation that
$$n\epsilon ^\vv=\theta^\vv \quad \mod \  \sum_{\beta\in \Delta^0_\mu}\Z
a^\vv_\beta\beta^\vv.$$

Now we claim  $\{n\epsilon^\vv\}\cup \{a^\vv_\beta\beta^\vv\}_{\beta\in\Delta^0_\mu}$ spans $N'$ as well. 
To verify the claim, let $m$ be in the sub-lattice  $N\cap(\sum_{\alpha_i\in\Delta^0_\mu}\R\alpha_i^\vv+\R\theta^v), $
then $m=\sum_{\alpha_i\in\Delta^0_\mu} r_i\alpha^\vv_i+r\theta^\vv,$
let $\Lambda_i$ act on both sides, by assumption $\Lambda_i(m)\in \Z$,
therefore $m_i=r_i+ra_i^\vv=\Lambda_i(m)\in \Z$. 
So $m=\sum_{\alpha_i\in\Delta^0_\mu}(m_i-ra^\vv_i)\alpha^\vv_i+r\theta^\vv.$
Replace $m$ by
 $m'=-\sum_{\alpha_i\in\Delta^0_\mu}ra^\vv_i\alpha^\vv_i+r\theta^\vv=r\sum_{\alpha_j\notin\Delta^0_\mu}a^\vv_j\alpha_j^\vv$,
where in the last equation the definition of $\theta^\vv$ is used.
Clearly $m'$ is in the same sub-lattice as $m$,
and $m'$ is a integer multiple of $\epsilon^\vv$ because $ra^\vv_j\in\Z$. This proves the claim.

From there, it is easy to identify what $N/N'\simeq T'_z/T_z$ is.
QED.

\subsection{The transformation of $T'_z/T_z$}
\begin{lemma}\label{bepsilon}
Under the action by $W_\mu$, the isotropy group $T'_z/T_z$ transforms
into itself.

For $ \phi$  with $<\phi,\theta>=1$ and $\mu=(\phi,1)$,  $$e^{2\pi i <\phi,w(b\epsilon^\vv+\sum_{ \alpha_i\in
\Delta^0_\mu }b_i\alpha_i^\vv) >}=e^{2\pi i<\phi,wb\epsilon^\vv>}=e^{2\pi
i<\phi,b\epsilon^\vv>}.$$

\end{lemma}
{\it Pf:} The first one is easy to verify using reflections defined by simple
roots of $(LG)_\mu$, to be more explicit:
$$r_\beta(b\epsilon^\vv+\sum_{ \alpha_i\in
\Delta^0_\mu }b_i\alpha_i^\vv)=b\epsilon^\vv+\sum_{ \alpha_i\in
\Delta^0_\mu }b_i\alpha_i^\vv-<b\epsilon^\vv+\sum_{ \alpha_i\in
\Delta^0_\mu }b_i\alpha_i^\vv,\beta>\beta^\vv$$ 
where $\beta$ is either $\theta$ or in $\Delta^0_\mu$. In either case,
because the coefficient of $\beta^\vv$ above is in $\Z$, after mod out
the lattice defining $T'_z$, it is clear that the element above is
in $T'_z/T_z$.

For the second part, first observe that $e^{2\pi i<\phi,\beta^\vv>}=1$
if $\beta\in \Delta_\mu^0$ or $\beta=\theta$. So the reflections does
not change the value, hence it is invariant under $W_\mu$.
Or $$e^{2\pi i<\phi,w(b\epsilon^\vv+\sum_{ \alpha_i\in
\Delta^0_\mu }b_i\alpha_i^\vv) >}=e^{2\pi i<\phi,b\epsilon^\vv+\sum_{
\alpha_i\in
\Delta^0_\mu }b_i\alpha_i^\vv >}=e^{2\pi
i<\phi,b\epsilon^\vv>}. \quad \mathrm{QED}$$

\subsection{The relations among the groups $T_z,T_z'$ and the maximal torus in
$(LG)_\mu^{\sss}$.}
We have just studied the relation between $T_z,T_z'$. A third Abelian
group is the maximal torus $S$ of $(LG)_\mu^{\sss}$, the three share the same
Lie algebra. The difference is the defining lattice.
\begin{lemma}\label{Sgroup}
The three groups are related as:
$$T_z'\rightarrow S\rightarrow T_z$$
where each arrow is a covering. 
When $\theta(\phi)\neq 1$,
$S\simeq T_z$. When $\theta(\phi)=1,$ $S/T_z\simeq \Z/n\Z$ where $n$ is
defined before.
\end{lemma}
{\it Pf:} If $\theta(\phi)\neq 1$, the coroot lattice of $(LG)_\mu^{\sss}$
is given by $\{\alpha^\vv\}_{\alpha\in\Delta^+_\mu}$, where
 each 
$\alpha^\vv$ is also a coroot of $\fg$. So the lattice defining $T_z$ is the
same as that defines $S$.

If $<\theta,\phi>=1$, the coroots are $\{-\theta^\vv\}\cup
\{\alpha^\vv\}_{\alpha\in\Delta^0_\mu}$ which form the  lattice defining $S$.  
On the other hand, the lattice defining $T_z$ is generated by
$\epsilon^\vv$ and $\{\alpha^\vv\}_{\alpha\in\Delta^0_\mu}$, therefore the
claim is verified. QED

\section{A couple of integration formulas}
One of the key steps in prove the cancellation formula 
is the following evaluation of certain 
differential forms on  a space which is a fiber bundle.
To be more precise, let $Z\rightarrow F\rightarrow E$ 
be an sequence so that $F=Z/S$ and $E=Z/K$ where $S$ is a the maximal torus
of $K$ and $K$ admits a local free action on $Z$.

We shall use the notations introduced in Sect 1 on the connection $dPA$ defining
the vertical and horizontal parts of $F$. Let $B,R$ be the vertical and
horizontal part of $ dPA $ respectively.

There is one exception  here,  $S$ is used instead of $T_z$ and $K$ is used instead of
$\KK''$.

\begin{proposition}\label{integration}
1).  Let $\epsilon$ be a weight of $S$, then the bundle $Z\times _S \C$ with
$(ps,v)\simeq (p,s^\epsilon v)$ defines a line bundle on $Z/S=F$ with 
curvatures given by $<\epsilon,dPA>$. The Chern class is given by
$<\epsilon, B+R>$.

2). The following holds: $\Td(TF)=\Td(T^HF)\cdot \Td(T''F),  \Td(T^HF)=\pi^*\Td(E) $ and 
$$
\Td(T''F)=\prod_{\tau\in {\Delta }_+}\frac{-\op<\tau, B+R>}{(1-e^{\op<\tau,B+R>})}$$
where ${\Delta}_+$ is the set of positive roots of $\kk$.

3).  Localization for a family:
 $$\int_{\pi^{-1}([p])} \Td(TF)e^{\op<\epsilon,B+R>}
=\Td(E)\sum_{u\in W(K)}\frac{e^{\op<u\epsilon,B>}}{\prod_{\tau\in {\Delta}_+}(1-e^{-\op<\tau, uB>})}.$$
 (This  equation and the next  should be viewed as  identities   about differential forms.)

4).\begin{equation}\begin{split}
& \int_{\pi^{-1}([p])} \frac{\Td(TF)}{\prod_\epsilon(1-e^{\lp<\epsilon, t+\fp
R+\fp B>)})}\\
&
=\Td(E)\sum_{u\in W(K)}
\frac{1}{\prod_{\tau\in  {\Delta}_+}(1-e^{<\tau, \op uB>})(1-e^{\lp <\epsilon,
 t+\fp uB>})}.
\end{split}\end{equation}
\end{proposition}
{\it Pf:} The first part repeats Prop.~5.1. Since $TF=T^HF\times T''F$, the
Todd class satisfies the equality $$\Td(F)=\Td(TF)=\Td(T^HF)\times \Td(T''F).$$ 
On the other hand, as discussion in Sect.1, $T''F\simeq Z\times _{S}\fn$,
where $\fn=\fs^\perp\subset \fk$. Hence $$T''F=\oplus_{\tau\in {\Delta}_+} Z\times_S\C,$$
where  $-\tau$ is the
character of $S$ acting on $\C$, the sign reflects the choice of the complex
structure on $X\times\X$. According to Part 1), the Chern class
of the line bundle is represented by  $-<\tau,B+R>$ and the Todd class of $T''F$ is the
product of the Todd class of the line bundle corresponding to each positive
root  $\tau$. Thus the expression is verified.

To see 3), we employ the localization of equivariant cohomology class
\begin{equation}
\int_{K/S}\prod_{\tau\in{\Delta}_+}\frac{-\lp<\tau,
\fp R+ X>}{(1-e^{\lp<\tau, X+\fp R>})} e^{\lp<\epsilon,X+\fp R>}=\sum_{u\in
W(K)}\frac{e^{2\pi i<\epsilon,uX>}}{\prod(1-e^{2\pi i<\tau,uX>})},
\end{equation}
the above is  an identity for  analytic  functions in $ X, \forall  X\in\fs$.
 If we plug the $\fs$-valued  2-form $B$ ,
instead of $X$, we end up with an equality of forms. This is exactly the claim
of Part 3). 

To see the claim of Part 4, inserting $0<r<1$ and introducing the multi-index
$\vn\in \Z^m$ with $m=\#\{\epsilon\}$, and
$\vepsilon=(\epsilon_1,...\epsilon_m)$. Then 
\begin{equation}
\begin{split}
&\int_{\pi^{-1}([p])} \frac{\Td(TF)}{\prod_\epsilon(1-re^{\lp<\epsilon,t+\fp R
+\fp B>})}\\
&=\Td(T^HF)\sum_\vn r^{| \vn|}e^{2\pi i\vn\cdot \vepsilon t}\int_{\pi^{-1}([p])} \Td(T''F)e^{<\vn\cdot\vepsilon,\op (R+B)>}\\
&=\Td(T^HF)\sum_\vn r^{|\vn|}e^{\vn\cdot\vepsilon t}\sum_{u\in
W(K)}\frac{e^{<\vn\cdot\vepsilon,\op uB>}}{\prod(1-e^{\op<\tau,uB>})}\\
&=\Td(T^HF)\sum_{u\in W(K)}\frac{1}{\prod(1-e^{<\tau,\fp uB>})\prod_\epsilon(1-re^{\lp<\epsilon,t+\fp uB>})},
\end{split}
\end{equation}
once the formula is established for $0<r<1$, we can take limit $r\rightarrow
1$.
QED

\section{$T$, $G$-spaces and  the consequences of the  main result}
Before we tackle the technical difficulty, the main  cancellation, 
let's first see a couple of 
consequences of the main theorem. 

There are two groups of results presented here, one is in the general
case and the other is the holomorphic case.

\subsection{Passing from $T$-modules to $G$-modules.} 
If $V_T$ is a $T$-module with  weight vectors in $\ft_+^*$, then one can
apply the holomorphic induction to get a $G$-module from it. Namely, let 
$$V_T=\oplus_{\lambda\in \ft_+^*}m_\lambda\C v_\lambda$$
where $v_\lambda$ is a weight vector with weight $\lambda$ and $m_\lambda\in
\Z_+$ the multiplicity of the weight $\lambda$, then
define  $$V_G=\oplus_{\lambda\in \ft_+^*}m_\lambda V_\lambda$$
where $V_\lambda$ is the unique highest weight $G$-module with $\lambda$ as
the 
highest weight.

The above is just the usual holomorphic induction.

Apply this construction to $H^0(X_N,L_N)$ to get a $G$-module which is denoted by
$V_G(X_N,L_N)$. 

The long passage we have taken
that started from the $\tLG$-module $H^0(X,L)$ to $H^0(X_N,L_N)$, and then
$V_G(X_N,L_N)$ has some advantage.  It is relatively easy to 
construct $X_N$ which should be viewed as the compactified quotient of 
$X$ by the maximal Borel $B^+\subset LG^\C$. We could have taken the quotient
by  $B^+_I\subset B^+$ which is $I$ at a fixed point on the circle. It would
be more difficult to find its compactification,
and will be done on another 
occasion. Although the space $X_N\times _TG$ which will be discussed more in the next section can be thought of as poor man's version of $X/B^+_I$.  

\subsection{The character functions of $T$ and $G$ modules.}
The character functions of the  modules $V_T,V_G$  are related. 
Let $$\chi_G(g)=\tr (t|_{V_G}), \quad \chi_T(t)=\tr(t|_{V_T}).$$
It is known that $\chi_G(g)$ is a class function, so its values are determined
by the restriction  to $t\in T$.
 As an easy
application of the famed Weyl character formula, one has the following:
\begin{equation}\label{TGcharacter}
\chi_G=\sum_{w\in
W}w\cdot\frac{\chi_T}{\prod_{\alpha\in\Delta^+}(1-e^{-\alpha})}.
\end{equation}

\subsection{Relating $T$, $G$-spaces}
The description above has  a generalization.
 Suppose that $P$ is a 
compact symplectic manifold (or orbifold) and $V$ is a complex line bundle on it, and the pair admits an
action by $T$, so that the action by $T$ on $P$ is Hamiltonian.  
Assume the above data fit together in the sense of geometric quantization,
c.f. [GS]. Then one can define the following $T$-equivariant Riemann-Roch:
\begin{equation}\label{RR}
\RR^T_P(t)=\int_P\tTd(TP)\tCh(V)(t),
\end{equation}
the above is also the equivariant index of a $\mathit{spin}^\C$-complex.
Via the fixed point formulas of Atiyah-Bott-Segal-Singer, 
the above can be written as contribution as $\sum_{F} \FC_F(t)$
where $\{F\}$ is  set of connected components of fixed points, $\FC_F(t)$ is  
an integral on $F$ involving $\Td, \Ch$ and equivariant classes of the normal
bundle of $F$ in $P$. The exact expression will be given later.

By assumption on $P$, there is a moment map $\phi:P\rightarrow \ft$. 
Suppose that $\phi(P)\subset \ft^+$ which is the positive Weyl chamber of
$\fg$,  we can associate a $G$-space with $P$, $P\times _TG$, a $G$-bundle
with $V$ over $P\times _TG$, $\pi^*(V)/T$ where $\pi:P\times G\rightarrow P$ is the projection,
the action by $t\in T$ on $\pi^*(V)$ is given by the following:
\begin{equation}
(p,g,v)\in \pi^*(V) \mapsto  (tp,tg, t^{\phi(p)}v).
\end{equation}
The following proposition is an easy exercise, as an application of the  Atiyah-Bott-Segal-Singer fixed
point formula on both $T$ and $P\times_TG$:
\begin{proposition}\label{TGRR}
\begin{equation}
 \RR^G_{P\times _TG}(t)=\sum_{w\in W}\frac{\RR^T_P(wt)}{\prod_{\alpha\in
\Delta^+}(1-(wt)^{-\alpha})}
\end{equation}
where the left side  represents  the $G$-equivariant Riemann-Roch number associated with
$(P\times_TG,\pi^*(V)/T)$.
\end{proposition}

\subsection{Contribution from fixed point sets}
Given a connected  component $F$ of $T$-fixed point set  in
$X_\ft=\mu^{-1}(\ft\times\{k\})$, define
\begin{equation}
\FC_F(t)=\int_F\frac{\tTd(F)\tCh(L|_F)}{\det(1-t^{-1}e^{-\Omega})|\no^\C(F,X_\ft)},
\end{equation}
where $\Td(F), \Ch(L_F,t)$ are Todd class and $T$-equivariant Chern class
of $TF, L|_F$  respectively; the denominator is the standard equivariant
class in the finite dimensional fixed point formula of the complex normal
bundle
of $F$ in $X_\ft$.
The existence of the complex structure on the normal bundle has been shown 
in Section 1.

\subsection{Coefficients of modular transformations.}
Use the notations introduced in Section 1, and let $\chi_a$ be the character
function of the highest weight $G$-module defined by  $a\in P_+$.
The following is in [K, Ch. 13]:
\begin{proposition}\label{modular}
\begin{equation}
\frac{(-1)^l}{\big|\frac{M^*}{(k+h^\vv)M}\big|}\sum_{\lambda \in
P^k_+}\big(\chi
_b \cdot \chi_{\bar{a}}\cdot D^2\big)( e^{2\pi
i\nu^{-1}(\frac{\lambda+\rho}{k+h^\vv})})=\delta_{b,a},
\end{equation}
where $\bar{a}$ is the  weight whose highest weight module is contragredient
to the one defined by $a$. (or $\bar{a}=w_L(-a)$ with  $w_L$ being the longest
element in $W$).
\end{proposition}

\subsection{Consequence of the main theorem}
 Assume the main theorem, how can we determine the
function $\RR(Y)$ from its values on the subset $\explife$?
We will prove Cor. 1.1 here.

\begin{lemma}\label{range}

1). $$\RR^T(X_N,L_N)(t)=\sum_{a\in P^k_+} m_a t^a $$
where $P^k_+$ is the set of weights in $kC$.

2).  $m_a=\RR(\calM_a, L_a)$ the Riemann-Roch number of  the pair
$(\calM_a, L_a)$ where $\calM_a=\phi^{-1}(a)/T$ and $L_a$ is the induced
orbifold line bundle on $\calM_a$.

3). $\RR^G(Y,L_Y)=\sum_{a\in P^k_+}m_a \chi_a$ where $\chi_a$ is the
character function of the highest weight $G$-module  defined by $a$.
\end{lemma}
The first assertion follows from an important fact that $m_a=0$ if
$a$ is outside the image of $k\phi$. We already know that the image lies
in $kC$, so $a$ must be in $kC$ if $m_a\neq 0$. That $a$ has to be a weight,
because $\RR^T(X_N,L_N)(t)$ is a function on $T$.

Part 2 follows the affirmative  solution (the Abelian case) to a conjecture by
Guillemin-Sternberg.

The last part uses part 1 together with the Weyl character formula:
\begin{equation}\begin{split}
\RR^G(Y,L_Y)&=\sum_{w\in
W}w\frac{\RR^T(X_N,L_N)}{\prod_{\alpha\in\Delta^+}(1-e^{-\alpha})}\\
&=\sum_{w\in
W}w\frac{\sum_{a\in P^k_+} m_a
t^a}{\prod_{\alpha\in\Delta^+}(1-e^{-\alpha})}\\
&=\sum_{a\in P^k_+} m_a\chi_a. \QED
\end{split}\end{equation}

{\it Pf. of Cor. 1.1:}  Let $\RR(Y)=\sum_{a\in P_+}m_a\chi_a,$  
then from the discussion of the $T$ and $G$-spaces, we know that
$$\RR(X_N)(t)=\sum_{a\in P_+}m_a t^a.$$
Thus only those $\{a\}$ inside $kC$ will occur, since 
$\phi(X_N)\subset kC$.
So one can write $\RR(Y)$ as $\sum_{a\in P^k_+}m_a\chi_a$. Multiplying $D^2\chi_{\bar{a}}$
on both sides and sum over $$\tau\in \explife$$ to get
\begin{equation}\label{98}
\begin{split}
m_a&=\frac{(-1)^l}{\big|\frac{M^*}{(k+h^\vv)M}\big|}\sum_{\tau} (\RR(Y) \cdot
D^2\cdot\chi_{\bar{a}})(\tau)\\
&=\frac{(-1)^l}{\big|\frac{M^*}{(k+h^\vv)M}\big|}\sum_{\tau}
 \chi_{\bar{a}}(\tau)\cdot D^2(\tau)\Big(\sum_{\{F|\mu(F)\in kW(C^\inte)\}}
\FC_F(\tau) +\RT(\tau)\Big). 
\end{split}\end{equation}
QED

\subsection{Application to the holomorphic case.}
Assume $X$ is holomorphic, and it satisfies the conditions that 
$\mu $ is both transversal and proper.

The result in this subsection does not rely on the main theorem of [C1].
First we  write down the  character formula/fixed point formula  of the module      
$V_G(X_N,L_N)$ in terms of the fixed points on $X$. 
\begin{theorem}\label{holomom}
Assume
$H^0(X_N,L_N)=\RR(X_N,L_N)$,  which holds if the higher cohomology groups
vanish.

1). $H^0(X_N,L_N)\simeq \sum _{a\in P^k_+} m_a \C_a$ where $C_a$ is the
$T$-module on which $T$ acts with weight $a$.

 2). Let $\chi_G(e^{it}), t\in\ft$ be the 
character of the $G$-module $V_G(X_N,L_N)$, then
\begin{equation}
\chi_G(\tau)= (\sum_{F}\FC_F+\RT)(\tau), \quad \tau\in \explife
\end{equation}
Furthermore, the multiplicity  is given by
of an irreducible component with the highest weight $a$ is 
given by the same expression as $m_a$ in  Eq. (\ref{98}).
\end{theorem}
{\it Pf:}
\begin{equation}
\begin{split}
\chi_T(X_N,L_N)(t)&:=\tr(t|H^0(X_N,L_N))\\
&=\sum_i(-1)^i\tr(t|H^i(X_N,L_N))\\
& =\RR(X_N,L_N)\\
&=\sum_{a\in P^k_+}m_at^a.
\end{split}
\end{equation}
The first part follows.

Using the earlier result relating $\chi_T(X_N, L_N)$ and $\chi_G$, we have
\begin{equation}
\begin{split}
\chi_G&=\sum_{w\in W}w\frac{\chi_T(X_N, L_N)}{\prod_{\alpha\in
\Delta^+}(1-e^{-\alpha})}\\
&=\sum_{w\in W}w\frac{\RR(X_N,L_N)}{\prod_{\alpha\in
\Delta^+}(1-e^{-\alpha})}\\
&=\RR(Y),
\end{split}
\end{equation}
where the last equality follows from Prop. \ref{TGRR}. 
From there, the assertion is shown to be true by the main theorem and Cor.
1.1.
QED 

\subsection{The character function of the $\cLG$-modules}

Now  we can derive a character formula for the representations of $\tLG$ on
$H^0(X,L)$,
  assuming the main theorem of [C1].

 The derivation here  requires certain formal
manipulations. The formal aspect to the approach here can be traced back
to the derivation of Weyl-Kac character formula.

Let $\chi_\tLG(X,L)(t)$ be the trace of $t$ acting on the part of $H^0(X,L)$
generated by the highest weigh modules of level $k$. The qualifier 
here about the trace is included because we do not know at this 
point whether $H^0(X,L)$ is generated by the highest weight modules. 
If $H^0(X,L)$ is a representation of  finite energy, then automatically it
has the desired quality, see [PS].    

\begin{theorem}\label{ABSScousin}
Define the regular and reduced Weyl-Kac denominators as follows
  $$  \tDwk=\prod_{\talpha>0}(1-e^{-\talpha}), \quad \tDwko=\prod_{\talpha>0,\talpha\notin \Delta^+}(1-e^{-\talpha})$$
where $\talpha$ in the first function  runs through all the  positive roots of $\clg$, while the second does not contain the positive roots   $\fg$, $\Delta^+$. 
Under the same assumptions as in Thm \ref{holomom}, the character
$\chi_\tLG(X,L)(t)$ is given by
\begin{equation}\begin{split}
&\chi_\cLG(X,L)(t)\\
&=\sum_{F\in \mathcal{F}} \int_{F}\frac{\Td(
F)\Ch(L|_{F},t)}
{ \det(I-t^{-1}e^{-\Omega})|_{\no^\C(F,X)}}+\sum_{w\in W^\aff/W}
w\frac{\RT(\tau)}{\tDwko}.
\end{split}
\end{equation}
where $\no^\C(F,X)$ is the complex normal bundle of $F$ in $X$.

\end{theorem}
{\it Remark:} It is known that Weyl-Kac formula has a formal flavor to it.
The infinite sum and product in the denominator above reflects that. 

{\it Pf:} It is clear that $\tDwk=D\cdot \tDwko$ with 
$D$ as in Section 1.1.

Given $a\in P^k_+$, there is a unique highest weight $\tLG$-module at level
  $k$, $\tV_a$. The character 
$\tr(t|\tV_a)$ is provided by the Weyl-Kac formula as
$$\tchi_a(t)=\sum_{w\in W^\aff}w\frac{e^a}{\tDwk}(t).$$
The resemblance to the characters of the $G$-modules is obvious.

From the representation theory of $\fg_\aff$,
we learn that each $T$-module with weights in $P^k_+$ induces such a $\tLG$-module, and
these are all the irreducible highest weight module of $\tLG$ at level $k$.

Now the part of $H^0(X,L)$ generated by the highest weigh modules of level $k$
is simply 
$$\oplus _{a\in P^k_+}m_a\tV_a,$$
from our identification of   all the highest weight vectors in $H^0(X,L)$.
Thus we conclude for $t\in \explife$:
\begin{equation}\begin{split}
\chi_\cLG(X,L)(t)
&=\sum_{w\in W^\aff}w\frac{ \RR(X_N,L_N)}{\tDwk}(t)\\
&=\sum_{w\in W^\aff/W}w\frac{\RR(Y)}{\tDwko}(t)\\
&=\sum_{w\in W^\aff/W}w\frac{\sum_F\FC_F+\RT }{\tDwko}(t).
\end{split}
\end{equation}
Next we  treat the first term. Since each  connected component of fixed point
sets in $X_\ft$ can be written as $wF$ for some $w\in W_\aff$, and
$F\in\calF^1$, it suffices to show 
\begin{equation}\label{dito}
\begin{split}
w\frac{\sum_{F\in \calF^1}\FC_F}{\tDwko}
& =\int_{wF}\frac{\Td(w F)\Ch(L|_{wF},t)}
{ \det(I-t^{-1}e^{-\Omega})|_{\no^\C(F,X)}}.
\end{split}
\end{equation}
The map $w:X \rightarrow  X$ preserves the complex structure,
naturally  $w$ induces isomorphism between  $TF,TwF$. Also
the symplectic form which is defines  $c_1(L|_{wF}), c_1(L|_F)$ up to a
constant, and is invariant under $LG$. Hence the two Chern classes 
are equal under pull-back. The equivariant Chern classes are related 
by $$w\Ch(L|_{F},t)=e^{2\pi i\mu (wt)+c_1(L|_F)}=w_*e^{2\pi i\mu
(wt)+c_1(L|_{wF})}$$
where $w_*$  pulls back  forms. And $\Td(TF)=w_*\Td(TwF)$. The only tricky
part is the identification of the classes associated with the normal bundles.

Let $D$ be in the denominator of $\FC_F$, i.e., 
$$D(t)=\det (1-t^{-1}e^{-\Omega})$$
which can be written in terms of Chern roots $\{x_i\}$,  the weights $
\{\theta_i\} $ and the roots of $\fg$  as
$$\prod_{i}(1-t^{-\theta_i}e^{-x_i})\prod_{\alpha\in \Delta^+}(1-e ^{-\alpha}),$$
the index here is finite, since $\no(F,X_\ft)$ is of finite dimension.

The manipulation of $w\tDwk$ is similar to the compact case, and we obtain
\begin{equation}\begin{split}
&\det (1-t^{-1}e^{-\Omega}) w (\tDwko)\\
&=w \prod_{i}(1-t^{-\theta_i}e^{-x_i})\prod_{\alpha\in\Delta^+(\fg_\aff)}(1-t^{-\alpha})\\
&=\prod_{i}(1-t^{-w\theta_i}e^{-x_i})\prod_{\alpha\in\Delta^+(\fg_\aff)}(1-t^{-
w\alpha})\\
&=\det (1-t^{-1}e^{-\Omega})|_{\nor(F,X)},
\end{split}
\end{equation}
where we have used the fact that $\nor_x(F,X)\simeq T_x\mu^{-1}(\ft)\oplus \lfg/\ft$.
We observe that 
1).  $w^*(x_i^w)=x_i$ where $x_i^*$ is
the Chern root of the corresponding line bundle at $wF$.
2). Each $(\lfg)_\alpha$ in
$(\lfg/\ft)$ induces  a trivial bundle over $F$ by group action, hence it
has no curvature. 
 3). The weights on  $wF$ on $\no(wF,X)$ is given by $w\{\theta_i\}\cup
w\Delta^+(\fg_\aff)$. 

After observing the above, immediately we obtain
$$R.H.= w_*\det(1-t^{-1}e^{-\Omega})|\no(wF,X).$$
Thus we complete the proof. QED

\section{The proof of the main cancellation }
\begin{proposition}\label{fantasy}
Let $F_h$ denote a connected component of fixed point sets (here $h$ may be
$I$). Suppose $\mu=\mu(F_h)\in k(\partial C,1)$ and is preserved by $W^\aff_\mu$.
Let $w$ be a lifting in $W_\mu^\aff$  of $[w]\in  W_\mu/W(\KK_h)$, $wF_h=F_{wh} $
(or denoted by $F^w_h$)
is another component of fixed point sets with the same value under $\mu$.

1). If $\mu(F_h)=k(\phi,1)$ with $\phi\in \partial C$, but $\phi\notin
C^\aff=\{\theta=1\}$,
then
$$\sum_{v\in W}v\frac{\sum_{[w]\in
W_\mu/W(\KK_h)}\FC_{wF_h}}{\prod_{\alpha\in
\Delta^+}(1-e^{-
\alpha})}(t)=0.$$
2). If $\phi\in \partial C\cap \{\theta=1\}$ and $vFC_{wF_h} $ has no pole
on $\explife$, then
$$\sum_{v\in W}v\frac{ \sum_{[w]\in
W_\mu/W(\KK_h)}\FC_{wF_h}}{\prod_{\alpha\in\Delta^+}
(1-e^{-2\pi i\alpha})}=0 \on \frac{M^*}{k+\cheh}.$$

\end{proposition}

 {\it Pf:}

{\it Step 1: Lifting action by $T$ to a covering group $T_p\times T'_Z$.}

A  key ingredient in Atiyah-Bott-Segal-Singer's fixed point formula
is the contribution from  the normal bundles which appear in the denominators
of an integral on the fixed point set.
Ignoring for a moment the complications from
orbifolds, the fixed point formula requires evaluating:
\begin{equation}\label{fixfor}
\FC_F(g)=\int_F\frac{\Td(TF)\Ch(L|_F,g)}{\det(1-g^{-1}e^{-\Omega})},
\end{equation}
where $\Omega$ is  $i/2\pi$ times the curvature operator of the complex normal bundle.

In what follows, we will  evaluate the integral above for $wF, wF_h$ and
for $g=t,vt$ where $w,v$ are elements in Weyl subgroups specified later.
Furthermore, we need to include the consideration that $\{F\}$  are
orbifolds, if $F$ is of type 2, 3, i.e., $\phi(F)\in \partial C$.

Suppose $v,w\in W_\mu$,  
 then the decomposition of $t=(t_{wp}, t^w_z)
\in \ft_{wp}\oplus \ft_z$ has the following property as
 shown in
Section 4:
$$ (vt)_{wp}=wt_p,\quad (vt)^w_z=vt_z+(vt_p-t_p)\in \ft_z.$$

Use  $t_z$ to denote $vt_z+(vt_p-t_p)$, there should be no confusion.
 And denote by
$e^{2\pi it'_z}, e^{2\pi i s'}\in T_z'$  certain  lifting of $e^{2\pi
it_z}\in T_z$, $ e^{2\pi is}\in T_z\cap T_p$  respectively.

1). If $\mu$ is on the affine wall,
 the lifting of $ve^{2\pi it}\in T$ to $T_{wp}\times T'_z$  are given by
\begin{equation}\label{lifting1}
(e^{2\pi (wt_p+ws)}, e^{2\pi i(t_z'-s'+\sum
b_i\alpha^\vv_i+b\epsilon^\vv)})\in T_{wp}\times T_z' \ \mathrm{with} \   e^{2\pi i(ws)}\in wI_p^0=I_{wp}^0, 
\end{equation}
where $0\leq b_i<a^\vv_i, 0\leq b<n$ with $n,\epsilon$ the same as in Corollary8.1.

2). If $\mu$ is not on the affine wall, then the lifting are  given by
\begin{equation}\label{lifting2}
(e^{2\pi (wt_p+ws)}, e^{2\pi i(t_z'-s'+\sum
b_i\alpha^\vv_i)})\in T_p\times T_z \ \mathrm{with} \   e^{2\pi i (w s)}
\in wI_p^0=I_{wp}^0.
\end{equation}
\begin{lemma}\label{equichern} Let $g$ be a lifting as above, and $v\in
W_\mu^0$ which is the subgroup of $W$ isomorphic to $W_\mu^\aff$.
For $<\phi,\theta>=1$,
$$g^\phi=e^{2\pi i<\phi,vt+b\epsilon>}=e^{2\pi i<v \phi,t+b\epsilon>};$$
for $<\phi,\theta>\neq 1$,
$$g^\phi=e^{2\pi i<\phi,vt>}.$$ 
\end{lemma}
{\it Pf:} First of all, $\phi$ is a weight on $T_p$, since $p$ is fixed by
$T_p$, and $\phi$ is the character of $L$ at $p$. Now   use   Lemma \ref{Sgroup}, Part 2 of Corollary \ref{toriciso}
and the fact that $\phi$ is a weight on $S$ in the notation there, we have the
assertion. QED

{\it Step 2: Defining the denominators.}

Let $g$ denote a lifting just described.
Recall $F_h$ is a strata associated with $F$ and is fixed by $T$ and $h$ in
the isotropy group of $F$.
From now, we use  $F_h$ instead of $F$ and treat $F$ as the special case
$h=I$. 
 We have treated
$F$ and $F_h$ separately so far, this leads to repetitions on occasions.
We will point out  whenever the difference
requires additional attention.

The following are various factors which will appear in the denominator
 as in Eq.(\ref{fixfor}) along $wF_h=F^w_h$.
The expressions are obtained using the weights and curvatures computations
done in Sections 3-5.

{\bf Remark}: The  signs in front of the weights below are determined 
by the following consideration:

a). $X_N$ is the reduced space of the product $X\times\X$, for the
map $\mu_X-\muX$. The weakly symplectic form is $\omega-\omega_\X$. Thus one chooses   the original $J$ on $X$, and $-J$ on  $T\X$ to make the form
semi-positive definite.
Hence a negative sign for the weights $\{\beta\},\{\lambda\}$ on the  the normal subbundle $\{\lfg_\mu/\ft\oplus H_z\}$ from $T\X$ 

b). The expression in the denominator involves $g^{-1}e^{-\Omega}$, therefore
another negative sign. 

c). The difference in the signs  for the term   $dA_h$ 
in various determinants below   was referred
to in the Remark after Prop. 5.2. 

Now we can write down the expressions for the denominators:

1). If $\mu$ is on the affine wall, let $\vb=(b,b_1,...,b_\II)$ with
$\II=\#\Delta^0_\mu$ the number of simple roots of $\fg$ which vanish
at $\mu/k$.

\begin{equation}\label{denomdef}
\begin{split}
&\D_0^w(vt)=\det(1-g^{-1}e^{-\Omega})|_{N_{wp}\oplus \nno(wZ_h,wZ)}
= \prod(1-e^{-2\pi i \gamma(t_p+s+1/4\pi^2 \nabla^2)}),\\
&\tD^w(vt,\vb)=\det(1-g^{-1}e^{-\Omega})|_{H_z}\\
&\  =(1-  e^{-2\pi
iw\mu(t_z'-s'+b\epsilon^\vv+\sum_ib_i\alpha^\vv_i-1/4\pi^2 dA_h)})\\
&\quad\quad \times \prod(1- e^{2\pi i w\frac{\Lambda_i-a_i^\vv\mu}{a^\vv_i}
(t_z'-s'+b\epsilon^\vv+\sum_ib_i\alpha^\vv_i - 1/4\pi^2 dA_h)}),\\
&\D^w(vt,b)=
(1-  e^{-2\pi iw\mu(t_z-s+b\epsilon^\vv -1/4\pi^2 dA_h)})\\
&\quad\quad  \times
\prod(1- e^{2\pi i w(\Lambda_i-a_i^\vv\mu)(t_z-s+b\epsilon^\vv- 1/4\pi^2 dA_h)}),\\
&\D_b^w(vt)=\det(1-g^{-1}e^{-\Omega})|_{(\lfg)_\mu/\kk_{wh}}
 =\prod_{\beta\in \Delta^+\setminus \Delta^+(\KK_{wh})}(1-e^{2\pi
iw\beta(t_p
+s+ 1/4\pi^2 dA_h}),\\
&\D_a^w(vt)=\prod_{\alpha\in\Delta^+}(1-e^{-2\pi i\alpha(vt)}).
\end{split}\end{equation}

2). If $\mu$ is off the affine wall, let $\vb=(b_1,...,b_\II)$,
the one   difference is:
  \begin{equation}\label{denomdef2}
\begin{split}
&\tD^w(vt,\vb)=\det(1-g^{-1}e^{-\Omega})|_{H_z}\\
&\ = \prod(1- e^{2\pi
 i w\frac{\Lambda_i-a_i^\vv\mu}{a^\vv_i}(t_z'-s'+\sum_ib_i\alpha^\vv_i
-1/4\pi^2 dA_h}),\\
&\D^w(vt)= 
\prod(1- e^{2\pi i w(\Lambda_i-a_i^\vv\mu)(t_z-s+b\epsilon^\vv-  1/4\pi^2  dA_h}).\\
\end{split}\end{equation}

Obviously in the above definition, when $h=I$, $\no(wZ,wZ_h)$ is 
trivial since $Z=Z_h$. 

{\it Step 3: Expressing $\FC_{F^w_h}(vt)$.}

Continue to let $g\in T_p\times T_z'$  where the first component of $g$ passes 
$w(hT^0_p)$ in $w(T_p)$, as a lifting of $vt\in T$. This form of  lifting appeared in Eq. (\ref{lifting1}),
(\ref{lifting2}).
In the notations just introduced, we have  along the normal bundle of $F^w_h$
\begin{equation}\label{wholedet}
\det(1-g^{-1}e^{-\Omega})=\D_0^w(g)\D_b^w(g)\tD^w(g)
\end{equation}
by  definitions of the factors on the right and the structure of the normal
bundle.

The line bundles  $L|_{F_h}, L|_{F^w_h}$ (orbifold bundles actually) are
related through the map $w:F_h\rightarrow wF_h=F_{wh}=F^w_h$ as
$$w^*(L|_{F^w_h})=L|_{F_h}, \quad w(\phi)=\phi.$$
This  relation holds for the pair   $Z_{wh}, Z_h$ obviously, and
it holds on the quotient $T_z$ because $w\in W_\mu$ normalizes $T_z$.
In terms of $F_h$,  following Lemma \ref{equichern}, one obtains 

1). when  $\mu$ is on the affine wall and $w\in W_\mu^\aff$,
\begin{equation}
\label{Ch}
\Ch(L|_{F^w_h},g)=g^{w\phi} e^{\omega}=g^{\phi} e^{\omega}=e^{2\pi i<kv\phi,t+b \epsilon^\vv>}e^{\omega},
\end{equation}
where $\omega$ is the symplectic form on $X_N$ restricted to $F_h$.

2). when $\mu$ is off the affine wall,
\begin{equation}
\label{Ch2}
\Ch(L|_{F^w_h},g)
=g^{w\phi} e^{\omega}=g^\phi e^{\omega}
=e^{2\pi i<kv\phi,t>}e^{\omega}.
\end{equation}

Depending on whether $\mu$ is on or off the affine wall, one obtains now 
an expression of $\FC_{F^w_h}$ in terms of an integral on $F_h$ as 
follows:
\begin{equation}\begin{split}
&\phi(\theta^\vv)= 1:\\
&\FC_{F^w_h}(vt)
=\frac{1}{|I_p^0||T'_z/T_z|}\sum_g\int_{F_h}\frac{\Td(TF_h)e^{2\pi i<kv\phi, t+
b\epsilon^\vv>}e^{\omega}}{\D_0^w(g)\D_b^w(g)\tD^w(g)};\\
&\phi(\theta^\vv)\neq 1:\\
&\FC_{F^w_h}(vt)=\frac{1}{|I_p^0||T'_z/T_z|}\sum_g\int_{F_h}\frac{\Td(TF_h)
e^{2 \pi i<kv\phi,t>}e^{\omega}}{\D_0^w(g)\D_b^w(g)\tD^w(g)}.\\
\end{split}
\end{equation}
In the above the summation is over  all the possible lifting of $vt$
in $w(hT_p^0)\times T_z'$ which is the isotropy group of the strata
$F_{wh}$.

{\it Step 4:  Summation  over $T'_z/T_z$.}

For an  upcoming calculation, it is crucial to replace the fractional weights
$$(\Lambda_i-a_i^\vv\mu)/a^\vv_i, \quad \alpha_i\in \Delta^0_\mu$$
by the integral weights on $T_z$, $$\Lambda_i-a_i^\vv\mu, \quad
\alpha_i\in \Delta^0_\mu.$$ 
This amounts to replacing $\tD^w$ by $\D^w$ and is an important step. The
basic observation is for $\mu$ on the wall, 
$$|T'_z/T_z|=n\prod_{\alpha_i\in\Delta^0_\mu}  a_i^\vv;$$
for $\mu$ off the wall,
$$|T'_z/T_z|=\prod_{\alpha_i\in\Delta^0_\mu}  a_i^\vv;$$
and more importantly:
\begin{equation}\begin{split}
\sum_{b_i}\frac{1}{\tD^w(g)}=\big(\prod_{\alpha_i\in\Delta^0_\mu}
a_i^\vv\big)\frac{1}{\D^w(g)}.
\end{split}\end{equation}
The last equation is based on two observations:

1).  Using geometric series expansion, and
the fact that  $\sum_{b_i}e^{2\pi i<a, \sum_ib_i\alpha^\vv_i>}$ is either
$0$ or $\prod_{\alpha_i\in\Delta^0_\mu}  a_i^\vv$, for all the possible weight
$a$, depending on $a$ is a weight on $T_z'$ but not on $T_z$, or $a$ is a
weight on $T_z$.

2).  Another  observation used in the above is
\begin{equation}\begin{split}
 & e^{-2\pi i<w\mu,t_z'-s'+b\epsilon^\vv+\sum_ib_i\alpha^\vv_i-\fp dA_h>}
=  e^{-2\pi i<w\mu,t_z-s+b\epsilon^\vv-\fp dA_h>},\\
 & e^{2\pi
 i <w(\Lambda_i-a_i^\vv\mu),t_z'-s'+b\epsilon^\vv+\sum_ib_i\alpha^\vv_i
-\fp dA_h>}=
 e^{2\pi
 i <w(\Lambda_i-a_i^\vv\mu),t_z-s+b\epsilon^\vv
-\fp dA_h>},
\end{split}\end{equation}
due to $\phi(\alpha^\vv_i)=0;
\Lambda_i(\epsilon^\vv),\Lambda_i(\alpha^\vv_i)\in
 \Z$;
also $-\mu, \Lambda_i-a_i^\vv\mu$ are  fundamental weights of $T_z$, their
values
at $e^{2\pi i(t_z'-s')}\in T'_z$ are  the same at the    projections $e^{2\pi
i(t_z-s)}\in T_z$.

Thus one has:
\begin{equation}\begin{split}
&\phi(\theta^\vv)= 1:\\
&\FC_{F^w_h}(vt)
=\frac{1}{n|I_p^0|}\int_{F_h}\frac{w^*\Td(TF_{wh})
e^{2\pi i<kv\phi,t+ \epsilon^\vv>}
e^{\omega}}{\D_0^w(vt)\D_b^w(vt)\D^w(vt,b)};\\
&\phi(\theta^\vv)\neq 1:\\
&\FC_{F^w_h}(vt)
=\frac{1}{|I_p^0|}
 \sum_g\int_{F_h}\frac{\Td(TF_h)e^{2 \pi i<kv\phi,t>}e^{\omega}}
{\D_0^w(vt)\D_b^w(vt)\D^w(vt)}.\\
\end{split}
\end{equation}

{\it Step 5: An identity of Todd classes.}
Given $w\in W_\mu$, there is the  map $$w: F_h=Z_h/T_z\rightarrow
wZ_h/T_z=Z_{wh
}/T_z.$$
\begin{lemma}\label{Todd}
\begin{equation}\begin{split}
&w^*(\Td(T^HF_{wh}))=\Td(T^HF_h),\\
&w^*(\Td(T''F_{wh}))=\prod_{\tau'\in\Delta^+(\KK_{wh})}\frac{-i/2\pi <w\tau',dA_h>}{1-e^{ i/2\pi <\tau',dA_h>}}.
\end{split}\end{equation}
\end{lemma}
{\it Remark:} The sign convention reflects again the the complex structure
chosen on $\lfg/\ft$.

{\it Pf:} If $\pi:F_h\rightarrow E_h$ as in Section 3,
then $\Td(T^HF_{wh})=\pi^*\Td(E)$. The same holds for $F_{wh}, E_{wh}$. It
is
easy to see that
$w:Z_h\rightarrow Z_{wh}$ is an  equivariant isomorphism  with respect to the
action
by $K_h, \Ad_w\KK_h=\KK_{wh}$ on $Z_h,Z_{wh}$ respectively.
Therefore, $\Td(E_h)=w^*\Td(E_{wh})$ whose pull-back to $F_h,F_{wh}$ yield the
first equation.

On $T''Z_{wh}$, after decomposing it  into  a sum   of line bundles
according to the roots $\{\tau'\}$, the curvature  is given
by $$-\oplus <\tau',dA_{wh}>=-\oplus <w\tau',dA_h> .$$
The second identity follows that. QED

{\it Step 6: A sufficient condition for the validity of the  cancellation.}

First of all, recall $$w^*\Td(TF_{wh})=\Td(T^HF_h)w^*\Td(T''F_{wh})$$ where 
$\Td(T^HF_h)=\pi^*\Td(E_h)$, with 
$$\pi: F_h=Z_h/T_z\rightarrow E_h=Z_h/\KK'_h.$$
Also it was shown in Prop. 4.4
 $D_0^w(vt)=D_0(t)$.
 We also know both forms below  $$
w^*\Td(T^HF_{wh})=\Td(T^HF_{h}), \quad  e^{\omega}$$ are  the pull-back of  forms on
$E_h$, because they are null in the fiber direction of the map  $\pi$. 
 Also clearly that $D_a^w(vt)$ is  constant with respect to the integration
variable. 
Therefore, 
in order to evaluate
$$\sum_g\int_{F_h}\frac{w^*\Td(TF_{wh})e^{2\pi
i<kv\phi,t+b\epsilon^\vv>}e^{\omega}}{\D_0^w(vt)\D_b^w(vt)\D^w(vt,b)}$$
(for $<\phi,\theta>= 1$) or 
$$\sum_g\int_{F_h}\frac{w^*\Td(TF_{wh})e^{2\pi
i<kv\phi,t>}e^{\omega}}{\D_0^w(vt)\D_b^w(vt)\D^w(vt)}$$
(if $<\phi,\theta>\neq 1$),
one can pull $D_0^w(vt), \Td(T^HF_{h})$ and $ e^{\omega }$ out,  when integrating along $\pi^{-1}([p])$;
those factors except $D_a^w(vt)$ are independent of $v\in W^0_\mu,w\in
W_\mu$. 
Thus to prove the proposition, which involving evaluating a  sum over
$W_\mu$ of the above integrals,  it suffices to prove
\begin{equation}\begin{split}
\label{equivalent}
&\phi(\theta^\vv)=1:\\
&\sum_{w\in W_\mu/W(\KK_h);v\in W_\mu^0} \int_{\pi^{-1}([p])}\frac{w^*\Td''(TF_{wh})e^{2\pi
ik\phi(t_p+s+b\epsilon^\vv)}}{\D_a^w(vt)\D_b^w(vt)\D^w(vt,b)}
=0, \on \latk ;\\
&\phi(\theta^\vv)\neq 1:\\
&\sum_{w\in W_\mu/W(\KK_h);v\in W_\mu^0} 
 \quad \int_{\pi^{-1}([p])}\frac{w^*\Td(T''F_{wh})e^{2\pi
ik\phi(t_p+s+b\epsilon^\vv)}}{\D_a^w(vt)\D_b^w(vt)\D^w(vt)}
=0.
\end{split}\end{equation}
In the above, a lifting of each $w\in W_\mu/W(\KK_h)$ to $W_\mu$ is fixed and
denoted by the same.

{\it Step 7: Turn the integrals along the fiber to a sum via equivariant cohomology.}

Recall from Section 3 $ dA_h=B_h+R_h$, where $B_h$ is the horizontal part of the
curvature, $R_h$ is the vertical part tangent to $\pi^{-1}([p])$.
Recall also that $\pi^{-1}([p])$ is   a finite quotient of  $\KK'_h/T_z$. Hence
the integrals can be pulled to $\KK'_h/T_z$.

We will continue to use the same notations for the pull-back to $\KK'_h/T_z$ of various curvature forms on $\pi^{-1}([p])$.
In Eq.(\ref{equivalent}), the term $\D_a^w(vt)$ is a constant on
$\KK'_h/T_z$. 
The integral of the rest were calculated using equivariant cohomology.  As
a straightforward application of 
 Formula 4 in Prop.~\ref{integration},   the answers are
\begin{equation}\begin{split}
&\phi(\theta^\vv)=1:\\
&\int_{\KK'_h/T_z}\frac{w^*\Td(T''F_{wh})e^{2\pi
ik\phi(t_p+s+b\epsilon^\vv)}}{\D_a^w(vt)\D_b^w(vt)\D^w(vt,b)}\\
&=\frac{e^{2\pi
ik\phi(t_p+s+b\epsilon^\vv)}}{\D_a^w(vt)}\sum_{u\in W(\KK_h)}
\frac{ 1}{\prod_{\tau'\in \Delta^+(\KK_{wh})} (1-e^{2\pi i<w\tau',
\fp uB>}) \rd_b^w \rd^w},\\
\end{split}\end{equation}
where 
\begin{equation}\begin{split}
\rd_b^w 
&=\prod_{\beta}(1-e^{2\pi iw\beta(t_p+s+ \fp uB_h)}),\\
\rd^w&=
(1-  e^{-2\pi iw\mu(t_z-s+b\epsilon^\vv +\fp uB_h)})
\prod(1- e^{2\pi i w(\Lambda_i-a_i^\vv\mu)(t_z-s+b\epsilon^\vv- \fp uB_h}),
\end{split}\end{equation}
in the above the $\{\beta\}$ are in $\Delta^+_\mu\setminus \Delta^+(\KK_{wh})$.
(In applying Prop.~\ref{integration} above, we treat factors
$\D_b^w(vt)\D^w(vt,b)$ collectively as
$\prod_\epsilon(1-e^{\lp<\epsilon, it+\fp R_h+\fp B_h>})$; and use $w\tau'$ instead of  
$\tau$ there.)

For the other case $\mu(\theta^\vv)\neq 1$, the only difference is in
the definition of $\rd_b^w$, 
\begin{equation}
\rd^w=
\prod_{i\in\II}(1- e^{2\pi i (w\Lambda_i)(t_z-s+b\epsilon^\vv-\fp uB_h}).
\end{equation}

{\it Step 8: The first  lucky break.}

We will now simplify the expression just obtained to be in a position
to  apply the fundamental formula of Section 6. 

The first break comes when we  group the factor involving the positive roots
$\{\tau'\}=\Delta^+(\KK_{wh})$, with the one involving $\{\beta\}=\Delta^+_\mu\setminus \Delta^+(\KK_{wh})$, we realize that the first factor can be  made look just  like the
second one. To be more specific: 
\begin{equation}\begin{split}
&\prod_{\tau'\in \Delta^+(\KK_{wh})} (1-e^{ \op<w\tau', uB_h>})\\
&=\prod_{\tau'\in \Delta^+(\KK_{wh})} (1-e^{2\pi i<w\tau', t_p+s+\fp uB_h>}).
\end{split}\end{equation}
There are two reasons for the above equation: 1). $e^{2\pi i u(t_p+s)}=e^{2\pi i(t_p+s)}\in hT^0_p$, since $u\in W(\KK_h)$ and $\KK_h$ commutes with $hT^0_p$ by its definition;
2). $e^{2\pi iw\tau'(t_p+s)}=1$, the adjoint action by $hT^0_p$ on $\kk_h$
is $I$ since the two commute by the definition of $\kk_h$, and $w\tau'$ is a root of $\KK_h$. 

In this form,  the product is clearly in the same species  as 
$\rd_b^w$. 

{\it Step 9:  The second break and the finale.}
The denominator
needs to be written in a form so we can apply Prop.\ref{fund}.
The denominator acquires an extra factor
after integration as shown  in the previous two steps. So the total is 
$$\calD^w_v=\D_a^w(vt)\D_b^w(vt)\D^w(vt,b)\prod_{\tau'\in \Delta^+(\KK_{wh})} (1-e^{2\pi
i<w\tau', t_p+s+\fp uB_h>})$$
if $\phi(\theta^\vv)=1$. For $\phi(\theta^\vv)\neq 1$, $\calD^w_v$ is defined
the same way except $\D^w(vt,b)$ is replaced by $\D^w(vt)$ as defined in
Step 2 above.

According to the previous step, $\D_b^w(vt)$ which has $\beta$ running over 
$\Delta^+_\mu\setminus \Delta^+(\KK_{wh})$  can be written together with
$\prod_{\tau'\in \Delta^+(\KK_{wh})} (1 -e^{2\pi
i<w\tau', t_p+s+ \fp uB_h>})$ as
\begin{equation}
\prod_{\beta\in \Delta^+_\mu} (1 -e^{2\pi i<w\beta, t_p+s+\fp  uB_h>}),
\end{equation}
which can be further written as 
\begin{equation}
\prod_{\beta\in \Delta^+_\mu} (1 -e^{2\pi i<w\beta, t_p+s+b\epsilon^\vv+ \fp uB_h>}),
\end{equation}
for $\beta$ is  a weight of $T_z$, therefore is trivial on the  $T'_z/T_z$ in
which $e^{2\pi ib\epsilon^\vv}$ lies.

Now the full expression of $\calD^w_v$ is given by
\begin{equation}\begin{split}
\calD^w_v
& = \prod_{\beta\in \Delta^+_\mu} (1 -e^{2\pi i<w\beta, t_p+s+b\epsilon^\vv+ \fp uB_h>})
\prod_{\alpha\in \Delta^+} (1 -e^{-2\pi i<v\alpha, t>})\\
&\quad \times 
\prod_{\lambda}(1-  e^{2\pi i(<v\lambda,t>-<w\lambda, t_p+s+b\epsilon^\vv +\fp
uB_h)})
\end{split}\end{equation}
where $\lambda$ runs through the fundamental weights of $(\lfg)_\mu$. If
$\phi(\theta^\vv)\neq 1$, just remove the term $b\epsilon^\vv$ in the above.
Also
\begin{equation}\begin{split}
&e^{2\pi i<w\beta, t_p+s+b\epsilon^\vv>}=e^{2\pi i<w\beta,
u(t_p+s+b\epsilon^\vv)>},\\
& e^{-2\pi i<w\lambda, t_p+s+b\epsilon^\vv>}=e^{-2\pi i<w\lambda,u( t_p+s+b\epsilon^\vv)>},
\end{split}\end{equation}
for $u\in W(\KK_h)$, this property can be deduced from
Lemma \ref{bepsilon} and the fact that  $u|_{hT^0_p}=I$, for $u\in W(\KK_h)$ commutes with
$hT^0_p$ as $\KK_h$ does.

As mentioned earlier, $w$ is a lifting of $W_\mu/W(\KK_h)$ in $W_\mu$,
and $u\in W(\KK_h)$. The product $wu$ on $\fg$ runs through the whole $W_\mu$.
Or $uw$ on $\fg^*$ goes through $W_\mu$.
Denote the combined by $w\in W_\mu$.

Let $y=t_p+s+b\epsilon^\vv+\fp B_h$.
Finally, for $v\in W^0_\mu,w\in W_\mu$, and $m$ being the rank of $\kk_h$ we
have the following after pulling out some of  the exponential terms:
\begin{equation}\label{sumden}\begin{split}
\calD^w_v
&\ = \prod_{\beta\in \Delta^+_\mu} (1 -e^{2\pi i<w\beta, y>})
\prod_{\alpha\in \Delta^+} (1 -e^{-2\pi i<v\alpha, t>})\\
&\quad \times
\prod_{\lambda}(1-  e^{2\pi i(<v\lambda,t>-<w\lambda, y>) })\\
&\ =(-1)^{m+\sigma(w)+\sigma(v)}e^{2\pi i(<w\rho_\mu-\rho_\mu
-w\sum\lambda, y>+<\rho-v\rho+v\sum\lambda,t>}\\
& \quad \times
 \prod_{\beta\in\Delta^+_\mu}(1-e^{-2\pi i<\beta, y>})
\prod_{\alpha\in \Delta^+} (1 -e^{-2\pi i<\alpha, t>})
\prod_{\lambda}(1-  e^{-2\pi i(<v\lambda,t>-<w\lambda, y>) }),
\end{split}\end{equation}
in the exponent of the first term above, there is  the 
cancellation  
$$\rho_\mu-\sum\lambda=0. $$

  The final calculation is
\begin{equation}\begin{split}
&\sum_{v\in W^0_\mu,w\in W_\mu}\frac{e^{2\pi i<vk\phi,t+b\epsilon^\vv>}}{\calD^w_v}\\
&\ =\frac{e^{2\pi i<-\rho_\mu,y>}}{\prod_{\alpha\in \Delta^+} (1 -e^{2\pi i<v\alpha, t>})
\prod_{\beta\in\Delta^+_\mu}(1-e^{-2\pi i<\beta, y>})}\\
&\quad \times \sum_{v\in W^0_\mu,w\in W_\mu}
\frac{(-1)^{m+\sigma(w)+\sigma(v)}
e^{2\pi i(<vk\phi, t+b\epsilon^\vv> -<\rho-v\rho+v\rho_\mu,t>)}}{\prod_{\lambda}(1-  e^{-2\pi i(<v\lambda,t>-<w\lambda, y>})},\\
\end{split}\end{equation}
after observing $e^{2\pi i<vk\phi,b\epsilon^\vv>}=e^{2\pi
i<k\phi,b\epsilon^\vv>}$, hence it can be taken outside the summation.
It now suffices to show the vanishing of 
\begin{equation}\begin{split}
\sum_{v\in W^0_\mu,w\in W_\mu}
\frac{(-1)^{m+\sigma(w)+\sigma(v)}
e^{2\pi i(<vk\phi +\rho-v\rho+v\rho_\mu,t>)}}
{\prod_\lambda( 1-  e^{-2\pi i(<v\lambda,t>-<w\lambda, y>)})}.\\
\end{split}\end{equation}
The numerator, when $\phi(\theta^\vv)=1$ and on  the lattice $\latk$
as shown in Prop. \ref{pride2} agrees with
$$e^{2\pi i<k\phi+\rho_\mu,t>}$$
which is independent of $v,w$ and can be pulled outside the summation.
When $\phi(\theta^\vv)\neq 1$, the numerator equals  
$$e^{2\pi i<k\phi+\rho_\mu,t>}$$ everywhere.
The lattice is invariant under $W$.
Therefore, over this lattice, the vanishing of the  above sum is 
implied by 
\begin{equation}
\sum_{v\in W^0_\mu,w\in W_\mu}
\frac{(-1)^{\sigma(w)+\sigma(v)} }
{\prod_\lambda( 1-  e^{-2\pi i(<v\lambda,t>-<w\lambda, y>)})}=0,
\end{equation}
whose validity is shown by Prop. \ref{fund}. 
Thus we have completed the proof of Part 1, 2 of  Prop. \ref{fantasy}.

\section{ Twins and   a new  surgery formula}
\subsection{A consequences from the proof: Twin pairs of compact  $G$-manifolds}

Here we give an easy application for finite dimensional symplectic $G$-manifolds, or more
generally symplectic $G$-orbifolds.

Let $M$ be a symplectic manifold (or orbifold)  with a Hamiltonian $G$-action,   $f$ be  the moment map. 
Suppose that $M$ is compact and $f^{-1}(\ft)$ is smooth, i.e. the image
of $f$ is transversal to $\ft$, then we can construct a $T$-orbifold $M_N$, just
as we have done for the $LG$-space $X$. 
Specifically, $M_N$ is constructed as follows: Let $k\in \Z_+$ so that
$f(M)\cap \ft_+\subset kC$. Let $X_\fg$ be the same toric variety as in
Section 2, with $\Phi$ as its moment map. 
Then $$M_N= \{(p,q)\in f^{-1}(\ft_+)\times X_\fg |f(p)=k\Phi(q)\}/T.$$
Another way to see it is to  start with
  $f^{-1}(\ft_+)$ which is a manifold with  non-smooth  boundary $f^{-1}(\partial \ft_+)$. Each face $Q$ on $\partial \ft_+$ has  a semi-simple 
stabilizer  $G_Q\subset G$ under the adjoint action. The normal subspace to $Q$ 
is the Lie algebra of a maximal torus of $G_Q$, $T_Q$.  Then $$M_N=\cup_Q f^{-1}(Q^\inte)/T_Q.$$  
The construction looks similar to symplectic cuts, but the two
have  quite different properties when one consider surgery formulas for
equivariant Riemann-Roch.

 The map $\phi(u) =f(p)$ with $u=[p,q]$, is well defined, since $f$ is
invariant under $T$. 
The variety $M_N$ enjoys the same property as $X_N$ in Section 2.3.
The $T$-fixed points on $M_N$ comes from either fixed points on $M$, or as a
result of surgery  along 
$$ \big(w, \phi^{-1}(\partial \ft_+)\big),\quad  w\in W.$$ The latter were studied in Section 3.

Similar to $X_N$,   $M_N$ has degeneracy as 
a symplectic orbifold and  has  a $T$-invariant  almost complex structure.

The  space $$M_G=G\times _T M_N$$ is a new $G$-space. In general, $M_G$ is
 different from $M$, when $f(M)\cap\ft$ intersects  the boundary of the Weyl
chambers.   Though the restriction of the image of the moment maps of the two
spaces coincide.

\begin{figure}\label{122}
\begin{center}
\setlength{\unitlength}{0.00041667in}
\begingroup\makeatletter\ifx\SetFigFont\undefined%
\gdef\SetFigFont#1#2#3#4#5{%
  \reset@font\fontsize{#1}{#2pt}%
  \fontfamily{#3}\fontseries{#4}\fontshape{#5}%
  \selectfont}%
\fi\endgroup%
{\renewcommand{\dashlinestretch}{30}
\begin{picture}(7074,3194)(0,-10)
\drawline(12,860)(387,35)
\drawline(2412,35)(387,35)
\drawline(987,2585)(1887,2585)
\drawline(1887,2585)(2862,860)
\drawline(1437,1160)(1437,1160)
\drawline(1437,1160)(1437,1160)
\drawline(1437,1160)(1437,1160)
\drawline(428,12)(428,12)
\drawline(428,12)(428,12)
\drawline(428,12)(428,12)
\drawline(987,2585)(12,860)
\drawline(2862,860)(2412,35)
\drawline(7062,860)(6687,35)
\drawline(4662,35)(6687,35)
\drawline(6087,2585)(5187,2585)
\drawline(5187,2585)(4212,860)
\drawline(5637,1160)(5637,1160)
\drawline(5637,1160)(5637,1160)
\drawline(5637,1160)(5637,1160)
\drawline(5637,1160)(5637,1160)
\drawline(5648,1141)(5648,1141)
\drawline(4212,860)(4662,35)
\drawline(5637,2585)(5637,35)
\drawline(2637,485)(2637,485)
\dottedline{45}(1437,2585)(1437,2585)
\drawline(1437,2585)(1437,2585)
\drawline(6087,2585)(7062,860)
\drawline(6837,410)(4737,1760)
\dottedline{90}(537,1760)(2637,485)
\dottedline{90}(2337,1760)(237,485)
\dottedline{90}(1437,2585)(1437,35)
\drawline(4437,485)(6612,1685)
\put(537,3035){\makebox(0,0)[lb]{\smash{{{\SetFigFont{6}{7.2}{\rmdefault}{\mddefault}{\updefault}$f(M)$}}}}}
\put(5037,3035){\makebox(0,0)[lb]{\smash{{{\SetFigFont{6}{7.2}{\rmdefault}{\mddefault}{\updefault}$f(M_G)$}}}}}
\put(4887,1685){\makebox(0,0)[lb]{\smash{{{\SetFigFont{6}{7.2}{\rmdefault}{\mddefault}{\updefault}$Q_1$}}}}}
\put(5412,2285){\makebox(0,0)[lb]{\smash{{{\SetFigFont{6}{7.2}{\rmdefault}{\mddefault}{\updefault}$Q$}}}}}
\end{picture}
}
\end{center}
\caption{The extra cuts on $f(M_G)$}
\end{figure}
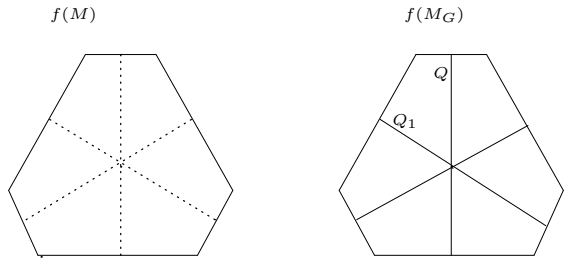

The figure~\ref{122} illustrates the cuts and the intersections of the images
of the two spaces with $\ft$.

As already mentioned, the new fixed points all have  their images on the
boundary of the Weyl chambers. As shown in the proof in  the last section, the sum of contributions of  
fixed points with the same image vanishes. 
Therefore, we have the following interesting consequence for the pair:
\begin{corollary}\label{twins}

The two $G$-spaces $M$ and $M_G$ have the same $G$-equivariant Riemann-Roch
 numbers.
\end{corollary}

\subsection{Another consequence: A new surgery formula }

Let $M$ be a symplectic  $G$-orbifold satisfying the transversality condition,
i.e. $f(M)$ is transversal to $\ft\subset\fg$,  and 
$V$ be   a  smooth symplectic subvariety of $M$. Suppose  the $T$-action preserves $V$,   with $V\subset f^{-1}(\ft)$. 

\begin{definition}
1). $(\tU, G_U)$ is an orbifold chart 
of $U$ if $G_U$ is an finite group acting with no non-trivial kernel on $\tU$
so that    
 $$\pi :\tU\rightarrow U=\tU/G_U,\quad \pi(\tilp)=G_U(\tilp)\in U.$$  

 2). $I_V\subset T$ is  the isotropy group of $V$ if for an open set $U$ with $U\cap V\nonempty$, and   
 $I_V\subset G_U$ is the stabilizer of $\pi^{-1}(V\cap U)$.  
For  $\tau\in T$  which fixes $V$, let $\tau _V$ denotes all the local
liftings
of $\tau$.
\end{definition}
From [C], we know that up to isomorphism, the group $I_V$ is independent of the chart, and the open set $U$. So it is well defined over $V$. 

{\it Remark:} The local liftings may not be extended over the entire $V$, since
there could be global  monodromy. On the other hand, the following
characteristic class over $V$ can be defined, as an average over $\tau I_V$: 
$$\frac{1}{|I_V|}\sum_{t\in \tau I_V}\int_V\frac{\tTd(V)\tCh(L_V)}{\det_{\no(V,M)}(1-t^{-1}e^{-\Omega})}$$
where $\Omega$ is the
 $i/2\pi$-curvature operator of the normal bundle of $V$ in 
$M$.
We want to find another expression for it.

 The expression is in terms of the subvarieties in  $M_G$ and its lower
strata.

For each face $Q$ in $\{w\ft_+ \big| w\in W\}$, $Q$ can be written as $wQ'$ with
$Q'$ as a face of $\ft_+$. 
Thus each $Q$  has a stabilizer in $W$, denoted by $W_Q$.
 Let
\begin{equation}\label{NQ}
\begin{split}
&M_Q= \{(p,q)\in f^{-1}(Q)\times X_\fg |f(p)=k w\Phi(q)\}/T;\\
&V_Q= \{(p,q)\in V\times X_\fg | f(p)\in Q, \quad f(p)=k w\Phi(q)\}/T.
\end{split}
\end{equation}

Both $M_Q$ and $V_Q$  are orbifolds, as a consequence of the transversality
of $f(M)$ to $\ft$,  and both have $T$-invariant almost complex structures.
Therefore one can define Todd class and makes sense of the equivariant
normal bundles of $V_Q$ in $M_Q$. 

The readers are warned that
there is no natural almost complex structure on the normal bundle of $V$ in $M_\ft$.

\begin{lemma}\label{gporder}
1). Let $M_Q^G=G\times M_Q.$
Then $M_Q, M_Q^G$ and   $V_Q$ are subvarieties in $M_G=G\times_T M_N$.

2). 
If $F$ is a subvariety of $V_Q$, then the isotropy groups, $I(F,V_Q),
I(F,M_G)$
of $F$ in $V_Q, M_G$ respectively satisfies the following exact short sequence:
$$1\rightarrow I(V_Q,M_G)\rightarrow I(F,M_G)\rightarrow I(F,V_Q)\rightarrow
1$$
where the second term  is the isotropy group of $V_Q$ in $M_G$. In particular:
$$| I(F,M_G)|=|I(V_Q,M_G)|\cdot |I(F,V_Q)|.$$
\end{lemma}
{\it Pf:}  Part 1  is obvious.
To see Part 2), let $(\tilU,G_U)$ be an orbifold chart meeting $F$.  If $g\in  I(F,M_G)$,   is $I$ on $\pi^{-1}(M_G\cap U)$, it remains  so on the smaller
set $\pi^{-1}(V_Q\cap U)$.
Hence $g\in I(F,M_G)$.  The next map in the sequence  is defined by
 restriction. 
If the restriction of $g$  is trivial, it is $I$ on $V_Q$ hence  $g\in I(V_Q,M_G)$.

The relation for the size of the group now is an immediate consequence. 
QED.

Suppose that $Q$ is a sub-face of $Q_1$, 
then the above construction shows that $V_Q$ is a subvariety of 
$V_{Q_1}$.  Obviously 
$$V_Q=V_{Q_1}\cap M_Q.$$

Let $\{\Delta\}$ be all the top dimensional faces, i.e. they are $\{w\ft_+\}$, 
if $Q\subset \Delta$, the normal bundle of $V_Q$ in $V_\Delta$ admits a 
$T$-action since $T$ acts on both.

The following is a new  surgery formula and will be used in the
next section to find the remainder term $\RT(\tau)$ in the fixed point formula:
\begin{proposition}\label{cutup}
Suppose $\tau\in T$ fixes $V$.
Let $\Lambda^{\max}\no(V_Q,V_\Delta)$ be the determinant
line bundle of  $\no(V_Q,V_\Delta)$, then
\begin{equation}\begin{split}
 &\frac{1}{|I_V|}\sum_{s\in \tau I_V} \int_V \frac{\tTd(V)\tCh(L_V )(t)}{\det_{\no(V,M) }(1-s^{-1}e^{-\Omega})}\\
&=\sum_\Delta\sum_{Q\subset \Delta}
\frac{1}{|W_Q||I_{V_Q}|}\sum_{s\in \tau I_{V_Q}} \int_{V_Q}\frac{\tTd(V_Q)\tCh(L_{V_Q}\oplus \Lambda^{\max}\no(V_Q,V_\Delta)(t))}{
\det_{\no(V_Q,M_Q^G)}(1-s^{-1}e^{-\Omega})}
\end{split}\end{equation}
where $W_Q, I_{V_Q},I_V$ are stabilizer of $Q$ in $W$ and  isotropy groups
of $V_Q, V$ in $M_G,M$ respectively.
\end{proposition}
{\it Pf:} 
 We prove the statement  by first representing both sides in terms of the $T$-fixed
points of $V, V_Q$.
Since $T$ acts on $V$ and its normal bundle, we can apply the localization
of equivariant cohomology classes on $V$ to get
\begin{equation}\begin{split}
 &\frac{1}{|I_V|}\int_V \frac{\tTd(V)\tCh(L_V)}{\det_{\no(V,M)}(1-t^{-1}e^{-\Omega})}\\
&=\sum_{F\subset V}\frac{1}{|I_V||I(F,V)|}\sum_{t\in \tau I(F,M_G)}\int_F \frac{\tTd(F)\tCh(L_F)}{\det_{\no(F,V)\oplus \no(V,M)}(1-t^{-1}e^{-\Omega}) |_F }
\end{split}\end{equation}
where the denominator can be combined as $\det_{\no(F,M)}(1-t^{-1}e^{-\Omega})$
obviously.

As for the right hand side,  for the same reason as above one  can express it as
$$\sum_\Delta\sum_{Q\subset \Delta} \sum_{F\subset V_Q}\frac{1}{m}
\sum_{t\in \tau I(F,M_G)}\int_F\frac{\tTd(F)\tCh(L_ F\oplus \Lambda^{\max}\no(V_Q,V_\Delta)|_F)}{\det_{\no(F,M_Q^G)}(1-t^{-1}e^{-\Omega})}$$
where $m=|W_Q||I_{V_Q}||I(F,V_Q)|$.

There are two kinds of fixed point sets on $\{V_Q\}$: those already
on   $V$,  
and those as a result of  the cuts along the boundary of Weyl chambers. 
The first kind have their isotropy groups come with the property that 
$I(F,M_G)=I(F,M)$ since the cut does not pass through $F$.
The second kind have  images under $\phi$ on the 
boundary of $wC$ for some $w$.

Also we may apply Lemma~\ref{gporder} 2) to the orders of the isotropy groups,
so that $ |I_{V_Q}||I (F,V_Q)|=|I(F,M_G)|$ independent of $V_Q$, if $\phi(F)$
is on the boundary of the Weyl chambers.

 Next we describe the  weights on the various bundles. 
Without loss of generality, assume $Q$  is a face of  $\ft_+$.
Let $\fg_\phi$ be the stabilizer of $\phi(F)\in Q\subset \fg$.
We claim that the total contribution of all the terms involving $F$ vanishes,
if $\phi(F)$ is on the boundary of the Weyl chambers.

Section 3 analyzed the formation of fixed point sets, here we continue to use the notations.

Let $\fg_\phi^\sss$ be the  semi-simple part of $\fg_\phi$,  
$\ft_\phi= \fg_\phi^\sss\cap \ft$  is a Cartan subalgebra of $\fg_\phi^\sss$.
Let $\Lambda_\phi$ be the set of  fundamental weights of $\fg_\phi^{\sss}$, and  
\begin{equation}\label{Qthing}
\Lambda_Q=\{\tlambda\in  \Lambda_\phi|\tlambda \in Q\},\quad  \Lambda_Q^\perp=\Lambda_\phi\setminus \Lambda_Q.
\end{equation}
Here $\ft_\phi$ is identical to $\ft_z$ in Section 3 and Section 4, which
is true since both subalgbras are the orthogonal complement of  the smallest 
face   containing $\phi(F)$.

As shown in Prop.~4.3, 
$\Lambda_\phi$ induces the following weights  along  $F$,
$$s^\lambda=(s_\phi)^\tlambda$$
for each lifting of $s$ to  $(s_p,s_\phi)\in T_p\times T_\phi$.
Or in terms of Lie algebra notation:
$$\exp2\pi i<\lambda,t>=\exp2\pi i<\tlambda,t_\phi>=\exp2\pi i<\tlambda,t-t_p>.$$
The set  $\{\lambda| \tlambda\in \Lambda_Q^\perp\}$ are all the weights 
on $ \no(V_Q,V_\Delta)|_F$.

We remark that for each fixed point set $F$ with image on the boundary of
$\ft_+^*$,  there is a submanifold  $Z\subset M$ with $Z/T_\phi$ where
$T_\phi=T\cap G_\phi^\sss$, just as in Section 3.
The action by $T_z$ is locally free, hence there is an associated 
connection $A$ on it.  Thus the equivariant Chern class of
$\Lambda^{\max}\no(V_Q,V_\Delta)|_F$
is given by
 $$\exp  \sum_{\tlambda\in \Lambda_Q^\perp} -2\pi i<\lambda,t-\fp dA >=\exp  \sum_{\tlambda\in \Lambda_Q^\perp} -2\pi i<\tlambda,t-t_p-\fp dA>$$
where $\exp 2\pi i (t-t_p)=s$.

If $\Delta=w\ft_+$, then the above is replaced by
 $$\exp  \sum_{\tlambda\in \Lambda_Q^\perp} -2\pi i<\lambda,w t-wt_p-1/4\pi^2wdA>.$$

Let $\{\gamma\}$, $\{\beta\}$, $\{\alpha\}$  be the same as in Prop.~4.3, where
 $X_N$ is replaced   by $M_N$,   then the weights to $F\subset V_Q$ in  $M_Q^G$ are
given by the same expressions as in Prop.~4.4, except only those $\lambda$ 
with  $\tlambda\in \Lambda_Q$ contribute, because those in $\Lambda_Q^\perp$ are normal 
to $M_Q$ or to $M_Q^G$.

As for the $\det_{\no(F,M_Q^G)}(1-s^{-1} e^{-\Omega})$, similar to the
expressions given in the last section,  we can express it 
in terms of the weights $\{\gamma\}$, $\Lambda_Q$,  $\{\beta\}$, $\{\alpha\}$:
$$ \prod_{\lambda\in\Lambda_Q}D_0^w(ws)D_a^w(ws)D_b^w(ws)$$
where $D_0,D_a,D_b$ are the same as in Eq. (\ref{denomdef}) of the last section.

There might be a non-Abelian $\KK\varsupsetneq$ commuting with $T_\phi$,  hence it acts on
$Z$.
 As was shown in Step 7 and 
Step 8 in the last section, the presence of a non-Abelian $\KK$ after
integrating the integrand along $\KK/T_\phi$ leaves little trace behind,
except replacing the two form  $\fp dA$ by  $\fp B$ which is a two form on
$Z/\KK$.  Therefore we deal
directly with the final expression of the denominator. 

So the fixed points contribution of $F\subset M_Q^G$ is given by
\begin{equation}\begin{split}
&\sum_{\Delta, Q\subset \Delta}\frac{1}{|W_Q|}\sum_{F\subset Q}\frac{1}{|I(F,M_G)|}\sum_{s\in \tau I(F,M_G)}\int_F \frac{\tTd(F)\tCh(L_F
\oplus\Lambda^{\max}\no(V_Q,V_\Delta)|_F)}{ \det_{\no(F,M_Q^G)}(1-s^{-1}e^{-\Omega})}\\
&=\sum_{w\in W, Q\subset w\ft_+}\sum_{F\subset Q}\frac{1}{|W_Q||I(F,M_G
)|}\sum_{t\in \tau I(F,M_G)}\int_F
\frac{\tTd( F)\tCh(L_F)}
{d}\\
&d= e^{\sum_{\lambda\in \Lambda_Q^\perp}2\pi i<\lambda,wt-wy >}
\prod_{\lambda\in \Lambda_Q}(1-e^{\lp<\lambda,wt-wy>}) D_0^w(ws ) D_b^w(ws)D_a^w(ws),
\end{split}
\end{equation}
where $y=t_p+\fp B$.

Recall that $D_0^w(ws)=D_0(s)$ as shown in Prop.~4.4, 
and use the expressions from Step 2 in the last section to obtain:
\begin{equation}
\begin{split}
D_b^w(ws)&=(-1)^{\sigma(w)}e^{\lp\sum_{w\beta<0} \lp<w\beta,y>}
\prod_{\beta\in\Delta_+(\fg_\phi)}(1-e^{\lp <\beta,y>})\\
&=(-1)^{\sigma(w)}e^{\lp< w\rho_\phi-\rho_\phi,y>}D_b(s);\\
D_a^w(ws)&=(-1)^{\sigma(w)}e^{-\sum_{w\alpha<0} \lp<w\alpha,t>}
\prod_{\beta\in\Delta_+(\fg_\phi)}(1-e^{-\lp <\alpha,t>})\\
&=(-1)^{\sigma(w)}e^{\lp<\rho_\phi-w\rho_\phi,t>}D_a(s).
\\
\end{split}
\end{equation}
This is the same procedure as in Eq.~(\ref{sumden}). 
 From Lie theory, we have $$\sum_{\lambda\in \Lambda_Q} w\lda+\sum_{\lambda\in \Lambda_Q^\perp} w\lda=w\rho_\phi.$$ Now the denominator can be written  as
\begin{equation}\begin{split}
d=& e^{\sum_{\lambda\in \Lambda_Q^\perp}<\lp\lambda,wt-wy>}
\prod_{\lambda\in \Lambda_Q}(1-e^{\lp<\lambda,wt-wy>})
D_0(\tau)D_b^w(w\tau)D_a^w(w\tau)\\
&=(-1)^{\sigma(w)+\sigma(w)}D_0(\tau)D_b(\tau)D_a(\tau)e^{\lp<\trho,t-y>}\prod_{\lambda\in \Lambda_Q}(e^
{-\lp<\lambda,wt-wy>}-1)).\\
 \end{split}
\end{equation}
The number of terms with a  fixed $Q$,  and varying $\Delta$ is exactly $|W_Q|$ which is the same number of $\Delta$ containing $Q$.
Gather them together and get rid of the term $1/|W_Q|$.
Now the summation over all the  $Q\subset \Delta=w\ft_+$ containing $\phi(F)$ to yield the following:
\begin{equation}\label{fpincutup}\begin{split}
&\sum_{ Q\subset \Delta}\frac{1}{|W_Q|}\int_F
\frac{\tTd( F)\tCh(L_F\oplus\Lambda^{\max}\no(V_Q,V_\Delta)|_F)}
{ \det_{\no(F,M_A^G}(1-s^{-1}e^{-\Omega})}\\
&=\sum_{ Q\subset \Delta}\int_F
\frac{\tTd( F)\tCh(L_F)}
{D_0(\tau)D_b(\tau)D_a(\tau)e^{\lp<\rho_\phi,t-y>})}
 \sum_{Q}
\frac{1}{\prod_{\lambda\in \Lambda_Q}(e^ {-\lp<w\lambda,t-y>}-1)}\\
&=\sum_{ Q\subset \Delta}\int_F
\frac{\tTd( F)\tCh(L_F)}
{D_0(\tau)D_b(\tau)D_a(\tau)e^{\lp<\rho_\phi,t-y>} }
\sum_{Q}
\frac{ (-1)^{\#\Lambda_Q} }{\prod_{\lambda\in \Lambda_Q}(1-e^ {-\lp<\lambda,wt-wy>})}.
\end{split}
\end{equation}
Summing  over various $\Delta$ in the above notations is the same as going
over  
$w\in W_\phi$, thus the total contributions of $F$ 
is 0 since
$$\sum_{w\in W_\phi} \sum_{Q}
\frac{(-1)^{\#\Lambda_Q} }
{\prod_{\lambda\in \Lambda_Q}(1-e^ {w\lambda}) }=0$$
by Eq.~(\ref{fund0}).

Thus only those $\{F\}$ on $V_Q$  with $\phi(F)$ not on the boundary contributes. They
 are exactly the   $T$-fixed point components on $V$. Therefore, the right
hand side is the same as the left one. QED

\section{Expression for the   remainder term $\RT$ }
If  for each fixed point set component $F$ of $Y$, $\FC_F$ has no pole   
on $$\explife, $$  then
the cancellation of Section 11 shows all the $T$-fixed points on the boundary
cancels on the subset.  However, such poles may occur, thus
the fixed points coming from compactification do not cancel readily.
The question is how they contribute.
Here we will find an explicit expression for the
 remainder term $\RT$.

First let's outline the steps in calculating the remainder term
$\RT$ as function on $\explife$. Instead of using the fixed points,
we will construct varieties which contain the fixed points on the
compactifying locus. Those varieties are obtained  using cuts which are transversal to $\partial wC$. Instead of transporting  the fixed
points  by translating elements in $W^\aff$, which fails to work if the
singular poles occur, we will move the varieties. 
  The varieties have   equivariant  Riemann-Roch which  are  well defined functions on $T$.

After transporting the varieties using the appropriate elements 
in $W_\aff$, one can cancel the  contributions of fixed points lying
on the compactifying locus, or on the boundary of certain Weyl chambers,  just as was  done in the last section.
Then we will apply Prop. \ref{cutup} and transport only $\{V_Q\}$ back   
to get a explicit formula of $\RT$.

\subsection{Partition of  the affine alcove}

The union $\bigcup_{w\in W^\aff}wC$ forms a tiling of $\ft$. Let $x$ be in the
interior of $C$, then $W^{\aff}(x)$ is the vertices of a dual
$W^\aff$-invariant
decomposition of $\ft$.

Each vertex $a$ of $C$ is now the center of a convex polytope with $W^\aff_v(x)$
as vertices.

 Let $l$
be  the rank of $G$, a little experiment shows that under the $W^\aff$
translations there are $l+1$ of different polytopes in the dual to
$\bigcup_{w\in W^\aff}wC$. Each one can be moved so that it contains 
exactly one of the vertices of $C$ as an interior point. 

Denote those $l+1$ polytopes by $\{R_a\}$.
Let $(LG)_a$ be the subgroup of $LG$ preserving  $$(a,1)\in \cft=\ft\times \R\subset \tlg,$$ under the co-adjoint action 
of the central extension.
By construction  $R_a$  is invariant under $W^\aff_a$ which is a subgroup of
$W^\aff$ and is the Weyl group of $(LG)_a$. The invariance comes from the fact
that the edges of $R_a$ are spanned by affine roots, hence invariant under 
the reflection by the corresponding  root.

\begin{figure}
\begin{center}
\setlength{\unitlength}{0.00041667in}
\begingroup\makeatletter\ifx\SetFigFont\undefined%
\gdef\SetFigFont#1#2#3#4#5{%
  \reset@font\fontsize{#1}{#2pt}%
  \fontfamily{#3}\fontseries{#4}\fontshape{#5}%
  \selectfont}%
\fi\endgroup%
{\renewcommand{\dashlinestretch}{30}
\begin{picture}(7010,6058)(0,-10)
\drawline(3912,4981)(3462,5806)
\drawline(4887,6031)(4887,6031)
\drawline(4212,2131)(5562,2131)
\drawline(4212,2131)(3837,1306)
\drawline(1812,1306)(3837,1306)
\drawline(3237,3856)(2337,3856)
\drawline(2337,3856)(1362,2131)
\drawline(5562,2131)(6537,3856)
\drawline(2798,4850)(2798,4850)
\drawline(2798,4850)(2798,4850)
\drawline(2787,4831)(2787,4831)(709,3631)
\drawline(2308,3861)(2308,3861)
\drawline(2308,3861)(2308,3861)
\drawline(2308,3861)(2308,3861)
\drawline(698,1212)(698,1212)
\drawline(698,1212)(698,1212)
\drawline(709,1231)(709,1231)(709,3631)
\drawline(1309,2131)(1309,2131)
\drawline(1309,2131)(1309,2131)
\drawline(687,1231)(2787,31)
\drawline(2776,4850)(2776,4850)
\drawline(2776,4850)(2776,4850)
\drawline(2787,4831)(2787,4831)(4865,3631)
\drawline(3266,3861)(3266,3861)
\drawline(3266,3861)(3266,3861)
\drawline(3266,3861)(3266,3861)
\drawline(4876,1212)(4876,1212)
\drawline(4876,1212)(4876,1212)
\drawline(4265,2131)(4265,2131)
\drawline(4265,2131)(4265,2131)
\drawline(4265,2131)(4265,2131)
\drawline(4265,2131)(4265,2131)(4190,2131)
\drawline(4898,1212)(4898,1212)
\drawline(4898,1212)(4898,1212)
\drawline(4876,1212)(4876,1212)(2798,12)
\drawline(3796,1283)(3796,1283)
\drawline(3796,1283)(3796,1283)
\drawline(3796,1283)(3796,1283)
\drawline(6987,2431)(6987,2431)
\drawline(6998,4850)(6998,4850)
\drawline(6998,4850)(6998,4850)
\drawline(6987,4831)(6987,4831)(4909,3631)
\drawline(6508,3861)(6508,3861)
\drawline(6508,3861)(6508,3861)
\drawline(6508,3861)(6508,3861)
\drawline(6987,2431)(6987,2431)
\drawline(4898,1212)(4898,1212)
\drawline(4898,1212)(4898,1212)
\drawline(4909,1231)(4909,1231)(4909,3631)
\drawline(5509,2131)(5509,2131)
\drawline(5509,2131)(5509,2131)
\drawline(5509,2131)(5509,2131)
\drawline(5509,2131)(5509,2131)(5584,2131)
\drawline(3237,3856)(4212,2131)
\drawline(65,2131)(65,2131)
\drawline(1362,2131)(1812,1306)
\drawline(12,2131)(1362,2131)
\drawline(3912,4906)(3237,3856)
\drawline(1812,1306)(1137,256)
\drawline(5562,2131)(6012,1306)(6012,1381)
\drawline(4512,181)(3837,1306)(1887,1306)
\drawline(6537,3856)(5862,4981)(3912,4981)
\drawline(2337,3856)(1737,4981)
\drawline(4876,1212)(4876,1212)
\drawline(4876,1212)(4876,1212)
\drawline(4898,1212)(4898,1212)(6976,12)
\drawline(5978,1283)(5978,1283)
\drawline(5978,1283)(5978,1283)
\drawline(5978,1283)(5978,1283)
\drawline(5978,1283)(5978,1283)(6015,1347)
\drawline(762,3631)(4887,1231)
\drawline(687,1231)(4887,3631)
\drawline(2787,4831)(2787,31)
\drawline(4887,6031)(2787,4831)
\drawline(4887,6031)(4887,3631)
\drawline(6987,2431)(4887,1231)
\drawline(6987,2431)(4887,3631)
\drawline(4887,6031)(6987,4831)(6987,31)
\put(3162,3931){\makebox(0,0)[lb]{\smash{{{\SetFigFont{6}{7.2}{\rmdefault}{\mddefault}{\updefault}$x$}}}}}
\put(2037,3856){\makebox(0,0)[lb]{\smash{{{\SetFigFont{6}{7.2}{\rmdefault}{\mddefault}{\updefault}$r_1(x)$}}}}}
\put(4287,2206){\makebox(0,0)[lb]{\smash{{{\SetFigFont{6}{7.2}{\rmdefault}{\mddefault}{\updefault}$r_2(x)$}}}}}
\put(3987,5131){\makebox(0,0)[lb]{\smash{{{\SetFigFont{6}{7.2}{\rmdefault}{\mddefault}{\updefault}$r_{\alpha_0}(x)$}}}}}
\put(1887,1456){\makebox(0,0)[lb]{\smash{{{\SetFigFont{6}{7.2}{\rmdefault}{\mddefault}{\updefault}$r_3(x)$}}}}}
\put(4812,3856){\makebox(0,0)[lb]{\smash{{{\SetFigFont{6}{7.2}{\rmdefault}{\mddefault}{\updefault}$a_2$}}}}}
\put(2637,4906){\makebox(0,0)[lb]{\smash{{{\SetFigFont{6}{7.2}{\rmdefault}{\mddefault}{\updefault}$a_1$}}}}}
\put(2862,2656){\makebox(0,0)[lb]{\smash{{{\SetFigFont{6}{7.2}{\rmdefault}{\mddefault}{\updefault}$a$}}}}}
\end{picture}
}
\end{center}
\caption{Partition of $\ft, wC$}
\end{figure}

Let $U_a=\mu^{-1}( (LG)_a(R_a))$, then $U_a$ is a symplectic $(LG)_a$-manifold
of finite  dimension, since $\mu$ is proper and $(LG)_a(R_a)$ is compact in
$\clg$. The
 proof that it is symplectic is
identical to the finite dimensional situation.

For $H\subset T$, the $H$-fixed points on $X_N$ have images under $\phi$
as linear subsets in $kC$, denote them by
$\{P_d\}$, $ 0\leq d\leq l-1$ is the  dimension of the connected component.
And each $P_d$ is a relatively open set in $kC$.
The following is easy to verify:
\begin{lemma}\label{transvers}
There is a point  $x\in C$, rational w.r.t. the weight lattice and in the
interior
or of $C$,
 such that the corresponding decomposition of $\ft$ as described above is
transversal to all the linear subsets $\{P_d/k\}$.
\end{lemma}

\subsection{Varieties corresponding to the partition}.

 For each $R_a$, $W^\aff_a(R_a\cap C)=R_a$.
So $R_a\cap C$ acts as the fundamental domain of $W^\aff_a$ on $R_a$. 

Let the faces of $R_a$ be denoted by $\{\Box\}$, then each $\Box$ is preserved 
by a subgroup in $W^\aff_a$, $W_\Box$.  Because each face of $R_a$ is
spanned by roots, and the reflections defined by those roots preserve $\Box$.
Anther way to see $W_\Box$ is as follows:
 If $Q$ is the smallest face of $wC$ which meets $\Box$, 
then the intersection has to be a point, otherwise even a smaller face can be
found to meet $\Box$. Let $p=\Box\cap Q$, then $\Box$ is perpendicular to
$Q$, since it is defined by the roots vanishing on $Q$. 
Therefore the group $W_Q\subset W^\aff$ which fixes $Q$, preserves $\Box$.
\begin{definition}\label{Gbox}
Let $(LG)_\Box$ denote the group which fixes the smallest face $Q$
intersecting $\Box$. $W_\Box\subset W^\aff$ be the subgroup preserving
$\Box$.
\end{definition}

Next we shall use $\{\Box\}$ to define symplectic cuts on both $X, Y=G\times
_T X_N$. 

 For each face $\Box$, the intersection $\Box\cap C$ is a
convex polytope.  The map $\phi:X_N\rightarrow kC$ defines for each $\Box$,
$\phi^{-1}(k(C\cap \Box) )$. For simplicity of notations, let $k=1$. 

For each sub-face $B$ of
$\Box$ in the interior of $C$, let $\ft_B$ be the subalgebra of $\ft$
perpendicular to the linear set defined by the face $B$, $T_B$ be the group
generated by $\ft_B$.
\begin{definition}
Let $X_\nbox$ be the cut space associated with $\phi^{-1}(C\cap \Box )$,
i.e.,  $$X_\nbox=\cup_B  \phi^{-1}(B)/\simeq.$$ 
\end{definition}
 That $X_\nbox$ is an $T$-orbifold with an invariant almost complex structure
 follows from the same argument for  $X_N$ itself.
Although $X_\nbox$ has degenerate symplectic form, just as $X_N$ does.
The orbifold line bundle $L_N$ also induces one on $X_\nbox$,
denoted by $L_\nbox$ which can be defined using quotients by $T_B$ on
$L_N| \phi^{-1}(C\cap \Box )$. 

The proof of the following can be found in [M]:
\begin{lemma}
The $T$-equivariant Riemann-Roch satisfies the following 
$$\RR_T(X_N,L_N)=\sum_{\Box\cap C\neq \emptyset}(-1)^{\cdim \Box}\RR_T(X_\nbox,
L_\nbox).$$
\end{lemma}
Therefore the corresponding $G$-spaces $Y$,  $Y_\Box=G\times_T X_\nbox$ and
the line bundles $L_Y, L_\Box$  satisfy
\begin{equation}\label{msurgery}
\RR_G(Y,L )=\sum_{\Box\cap C\neq \emptyset}(-1)^{\cdim
\Box}\RR_G(Y_\Box, L_\Box).
\end{equation}
{\it Remark}: Readers should compare this  surgery formula with that in
Prop. \ref{cutup}.

There is a map $\phi:Y_\Gbox$ which serves as moment map, for the degenerate
symplectic form, just as we have shown for $Y=G\times_T X_N$. 

The image $\phi(Y_\Gbox)\cap \ft$ is contained in $W(k(\Box\cap C))$, since
the image of $X_\nbox$ is in $(k(\Box\cap C))$.

\begin{figure}
\begin{center}
\setlength{\unitlength}{0.00041667in}
\begingroup\makeatletter\ifx\SetFigFont\undefined%
\gdef\SetFigFont#1#2#3#4#5{%
  \reset@font\fontsize{#1}{#2pt}%
  \fontfamily{#3}\fontseries{#4}\fontshape{#5}%
  \selectfont}%
\fi\endgroup%
{\renewcommand{\dashlinestretch}{30}
\begin{picture}(2124,2439)(0,-10)
\drawline(50,95)(50,95)
\drawline(99,49)(99,49)
\drawline(34,2412)(34,2412)
\drawline(2112,1212)(12,12)
\drawline(12,2412)(2112,1212)
\drawline(12,2412)(12,12)(12,12)
\drawline(462,1437)(12,1437)
\drawline(537,1437)(837,1887)
\drawline(500,1417)(950,592)
\put(162,1812){\makebox(0,0)[lb]{\smash{{{\SetFigFont{6}{7.2}{\rmdefault}{\mddefault}{\updefault}$\Box$}}}}}
\put(612,987){\makebox(0,0)[lb]{\smash{{{\SetFigFont{6}{7.2}{\rmdefault}{\mddefault}{\updefault}$\Box_1$}}}}}
\end{picture}
}
\end{center}
\caption{The ranges of the images of $X_{\nbox}$}
\end{figure}

\subsection{Another space $Z_\Box$ associated with $\Box$}
Define $Z_\Box$ to be the $(LG)_\Box$-symplectic   orbifold which is the cut
space  associated with
$\mu^{-1}(LG_\Box(\Box))$.
There are two ways of defining it. The first one is as in [C3],
where a holomorphic $LG\times (LG)_\Box$ symplectic orbifold $\calM_\Box$ was constructed, 
then $Z_\Box$ can be defined as the symplectic reduced space of the product
$X\times \calM_\Box$.

The other approach is given in [M].

The space $Z_\Box$ is compact, has a moment map $\Phi$ whose image $\Phi(Z_\Box\cap
\ft)\subset\Box$. 

If $LG_\Box\subseteq G$, i.e., $\Box\cap C^\aff=\emptyset$,
the two space $G_\Box\times_T X_\nbox, Z_\Box$ are related as twins in the sense of the
previous section, where $M, M_G$ would be $Z_\Box$ and $ G_\Box\times_T
X_\nbox$ respectively.
In particular they share the same Riemann-Roch.

On the other hand if $\Box\cap C^\aff\notempty$, then we will see below how the
Riemann-Roch are related.
\subsection{Riemann-Roch of $X_\Gbox$ and $Z_\Box$}

For $w\in W$, clearly $w(\Box\cap C)\subset wC$. If $W^\aff_\Box$  preserves $\Box$,
then $Ad_{w}W^\aff_\Box$ preserves $w\Box$, $w\in W$.

 In the following $W_\Box$ will be used in place of $W^\aff_\Box$.

For $W_\Box\subset W^\aff$, there is the isomorphic subgroup  $W_\Box^0\subset
W$. The proof of this simple fact is identical to that of the last statement
in  Prop. 7.2.

For each  $w$ is in $W_\Box^0$, let  $w'$ denote  the corresponding element in $W_\Box$.
\begin{lemma}\label{translation}
For every pair $w,w'$, there is a translating  element $v$ in the long root
lattice  
so that $w=vw'$. 
\end{lemma}
{\it Pf:}
This is a fairly simple fact. Since I can not find  a reference for
that, the proof is included. The proof is based on induction of the length of
$w $.

 By applying an element in $W$, we may assume that $\Box\cap C\notempty$ since
$W$ preserves the long root lattice.  For such a $\Box$,
all the simple roots of $(lg)_\Box$ are all simple, or contain $\alpha_0$.

Suppose $w$ is $r_i$, then $$ r_i r_i'=r_i(r_i')^{-1}$$ is either $I$, or 
$-\theta$ depending on  whether $r'_i$ is defined by a simple root of $\fg$ or
by $\alpha_0$. Thus the assertion holds for elements of length 1.

Assume it's proved for elements with length less $n$. Suppose  $w$ has length $n$, $w=r w_1$ with $w_1$ of length less $n$, and $r$ is one of the
generating reflections.  By  induction assumption $w_1=v_1 w_1'$, and $v_1$ is a 
translation.   One can write  $w'=r'w_1'$,  then 
$$w=rv_1 w_1'= rv_1r^{-1} rr' w_1'=rv_1r^{-1}rr' w'=vw',$$
where both $rv_1r^{-1}$,  $rr'$ are translations by  elements in the long root
lattice.
 So is  $v$ as the composition of the two.
 QED

The following is  essential, it explains in geometric context the appearance
of  the lattice $\latk$. By the definition of $M^*$, we have
$e^{(k+h)v}=1 \onlatk$.  The first part generalizes Prop.
7.2.
\begin{proposition}\label{rule}
1). Let $\rho_\Box, \rho$ be the half sum of positive roots of $(lg)_\Box,
\fg$ respectively. Let $D_\Box, D$ be their Weyl denominators, and $w=vw'$
with $w,w',v$ as in the previous lemma.

Suppose $\lambda$ is of level $k$, then
\begin{equation} 
w\frac{e^\lambda}{D}=e^{(k+h^ \vv)v}\frac{D_\Box}{D}w'\frac{e^\lambda}{D_\Box}.
\end{equation}
In particular,
\begin{equation} 
w\frac{e^\lambda}{D}=
\frac{D_\Box}{D}w'\frac{e^\lambda}{D_\Box}
 \onlatk.
\end{equation}
where $v$ is the translation in the previous lemma.

2). As function on $T$, the following holds
 \begin{equation} w\int_{F}\frac{\tTd(F)\tCh(L_F\oplus H)}{\det(1-t^{-1}e^{-\Omega} )}
=
 e^{(k+h^ \vv)v}\frac{D_\Box}{D}w'
\int_{F}\frac{\tTd(F)\tCh(L_F\oplus H)}{\det(1-t^{-1}e^{-\Omega})}
\end{equation}
where $H$ is a bundle of level $0$ on which $W^\aff$ acts.
If there is no pole on $\latk$, then  
\begin{equation} 
w\int_F\frac{\tTd(F)\tCh(L_F\oplus H)}
{\det(1-t^{-1}e^{-\Omega})}
= \frac{D_\Box}{D}w'\int_F 
\frac{\tTd(F)\tCh(L_F\oplus H)}{\det(1-t^{-1}e^{-\Omega})}  \onlatk.
\end{equation}

\end{proposition}
{\it Pf:} By conjugation, we may assume without loss of generality that 
$\Box\cap C\notempty$. 

  The following is  equivalent to the assertion, after applying the well-known
transformation rule on the Weyl denominator: 
$$ e^{w\lambda+\rho-w\rho}=e^{(k+h^ \vv)v}e^{w'\lambda+\rho_\Box-w'\rho_\Box}$$
which can be verified the same way as Prop. 7.2, replacing $\rho_\mu$ there
by $\rho_\Box$.

To see the second part, lift $w,w$  to  $N(T)$, they act
on     $F$, the  normal bundle and $L_F$. Notice that the line bundle has level
$k$, by assumption on the moment map $\mu$, i.e. the central part of $\clg$
acts with weight $k$ on the fiber. But the central part acts trivially on $X$,
the action on the normal bundle is of weight $0$. Thus if $\lambda$ is a
weight on the normal bundle, we have $w\lambda=w'\lambda$. 
Hence, we have the desired identity as a consequence of part 1). 
QED

\subsection{ Moving $X_\nbox$ and the consequence}

As we have mentioned earlier when $\FC_F$ has  poles  on $ \latk$, we
can not replace $ w\FC_F$ by $ \frac{D_\Box}{D}w'\FC_F$ on $\latk$. On the other
hand, the function $\RR_T(X_\nbox,L_\nbox)$ is a polynomial,  thus it can be
evaluated everywhere. For that function, the following holds: 

The above transformation rule yields the following 
\begin{corollary}\label{moving}
\begin{equation}
 w\frac{\RR_T(X_\nbox,L_\nbox)}{D}=e^{(k+h^ \vv)v}\frac{D_\Box}{D}w'\frac{\RR_T(X_\nbox,L_\nbox)}{D_\Box};
\end{equation}
Furthermore,
\begin{equation}
 w\frac{\RR_T(X_\nbox,L_\nbox)}{D}=
\frac{D_\Box}{D}w'\frac{\RR_T( X_\nbox,L_\nbox)}{D_\Box} \onlatk
\end{equation}
For $X_\Gbox=G\times_T X_\nbox$, one has
\begin{equation}
\RR(X_\Gbox,L_\Gbox)=
\sum_{u\in W/W_\Box^0}u \frac{D_\Box}{D}\cdot  u\RR(Z_\Box,L_\Box)\onlatk.
\end{equation}

\end{corollary}

{\it Pf:} The first one is an immediate consequence of Prop.\ref{rule}, after
one writes both sides in terms of the $T$-fixed points contributions.

Since the function $\RR_T(X_\nbox,L_\nbox)$
is the equivariant index of a $spin\C$ complex defined by the pair
$X_\nbox,L_\nbox$, it is a well defined function everywhere.
Hence we can evaluate on the lattice $\latk$ to  get the second formula.

To see the next identity,  expand $\RR(X_\Gbox,L_\Gbox)$ in terms of
$\RR(X_\nbox,L\nbox)$, then apply the first and second  identity to yield
\begin{equation}
\begin{split}
\RR(X_\Gbox,L_\Gbox)&=\sum_{w\in W}w\frac{\RR_T(X_\nbox,L_\nbox)}{D}\\
&=\sum_{u\in W/W_\Box^0, w\in W_\Box^0}uw \frac{\RR_T(X_\nbox,L_\nbox)}{D}\\
&=\sum_{u\in W/W_\Box^0}u\big(\frac{D_\Box}{D}\sum_{w'\in W_\Box} w'\frac{\RR_T(X_\nbox,L_\nbox)}{D_\Box}\big)\onlatk,\\
\end{split}
\end{equation}
it is easy to recognize the sum $\sum_{w'\in W_\Box} 
w'\frac{\RR_T(X_\nbox,L_\nbox)}{D_\Box}$ is simply the Riemann-Roch of the space
$(LG)_\Box\times _T X_\nbox$. How is it related to $Z_\Box$? 
They are twin-pairs as discussed in the last section, replacing 
$M,G$ there by $Z_\Box, (LG)_\Box$.   
Thus one has 
\begin{equation}
\RR(X_\Gbox,L_\Gbox)=\sum_{u\in W/W_\Box^0}u\big(\frac{D_\Box}{D} \RR(Z_\Box,L_\Box)\big)\onlatk.
\QED
\end{equation}
So the above relates the Riemann-Roch of $G$-space $X_\Gbox$ and the
$LG_{u\Box}$-space $\{Z_{u\Box}:=uZ_\Box\}$.

The fig~\ref{133} illustrates the relations between the intersections with $\ft$
of  the  images of $X_\Gbox$ and three $Z_{u\Box}$. The three separate regions, inside the middle hexagon, with dotted lines
 are
associated with $X_\Gbox$.  In this  case $W_\Box^0\simeq W$, the image of
$Z_\Box$ meeting $\ft$ inside one of  the regions filled with dashed lines.
What are the other two identical regions? If one starts with $u\Box$,
 in this case $u\neq I, r_3$, then $W_{u\Box}^0 \simeq W$ holds as well. 
And   one of the other two regions will  contain  the
intersection of the image with $\ft$ of  $uZ_\Box$.

In the second figure, $W_\Box\simeq \Z_2$, thus $W_\Box^0$ is not the same as
$W$, $W/W_\Box^0$ has three elements. The union of the six short segments   
inside the hexagon contains $\mu(X_\Gbox)\cap \ft$. The long segments 
contain the intersections with $\ft$ of the  images of $\{uZ_\Box\}$. 

\begin{figure}\label{133}
\begin{center}
\setlength{\unitlength}{0.00029167in}
\begingroup\makeatletter\ifx\SetFigFont\undefined%
\gdef\SetFigFont#1#2#3#4#5{%
  \reset@font\fontsize{#1}{#2pt}%
  \fontfamily{#3}\fontseries{#4}\fontshape{#5}%
  \selectfont}%
\fi\endgroup%
{\renewcommand{\dashlinestretch}{30}
\begin{picture}(14000,6591)(0,-10)
\drawline(2483,5095)(2483,5095)
\dottedline{90}(3008,6004)(1658,3666)
\dottedline{90}(4617,5941)(5863,3749)
\dottedline{90}(2490,2406)(5141,2348)
\dottedline{90}(3297,1854)(4171,1869)
\dashline{60.000}(5168,2246)(3893,37)
\dashline{60.000}(4799,2805)(3186,12)
\dashline{60.000}(4166,2910)(2816,572)
\dashline{60.000}(3459,2886)(2447,1132)
\drawline(5187,3710)(4812,2885)
\drawline(4212,5435)(3312,5435)
\drawline(3762,4010)(3762,4010)
\drawline(3762,4010)(3762,4010)
\drawline(1673,2791)(1673,2791)
\drawline(1673,2791)(1673,2791)
\drawline(1662,2810)(3762,1610)
\drawline(5851,2791)(5851,2791)(3773,1591)
\drawline(2337,3710)(2787,2885)
\drawline(1662,5210)(5862,2810)
\drawline(3762,6410)(3762,1610)
\drawline(5884,2810)(5884,2810)(5884,5210)
\drawline(2141,1765)(2141,1765)
\drawline(2141,1765)(2141,1765)
\drawline(5383,1765)(5383,1765)
\drawline(5383,1765)(5383,1765)
\drawline(5383,1765)(5383,1765)
\drawline(4384,35)(4384,35)
\drawline(4384,35)(4384,35)
\drawline(41,5440)(41,5440)
\drawline(41,5440)(41,5440)
\drawline(3762,1610)(4962,935)
\drawline(3762,1610)(2562,935)
\drawline(3762,6410)(1662,5210)
\drawline(5862,5210)(1662,2810)
\drawline(3762,6410)(5862,5210)
\drawline(5383,1765)(5383,1765)
\drawline(5383,1765)(5383,1765)
\drawline(3087,35)(3761,35)
\drawline(4437,35)(3761,35)
\drawline(12,5435)(987,3710)
\drawline(41,5440)(41,5440)
\drawline(41,5440)(41,5440)
\drawline(41,5440)(41,5440)
\drawline(987,3710)(12,5435)
\drawline(41,5440)(41,5440)
\drawline(41,5440)(41,5440)
\drawline(41,5440)(41,5440)
\drawline(612,6560)(12,5435)
\drawline(1662,6560)(612,6560)(612,6560)
\drawline(1662,3710)(987,3710)
\drawline(1684,2810)(1684,2810)(1684,5210)
\drawline(1662,5210)(537,4535)
\drawline(1662,5210)(312,5960)
\drawline(5862,5210)(7212,5960)
\drawline(5862,5210)(6987,4535)
\drawline(2112,1760)(3087,35)
\drawline(4812,2885)(5412,1760)
\drawline(3762,1610)(3761,35)
\drawline(5412,1760)(4437,35)
\drawline(7483,5440)(7483,5440)
\drawline(7483,5440)(7483,5440)
\drawline(7512,5435)(6537,3710)
\drawline(5862,5210)(5862,6560)
\drawline(6837,6560)(5862,6560)
\drawline(6837,6560)(7512,5435)
\drawline(5862,5210)(5862,6560)
\drawline(4887,6560)(5862,6560)
\drawline(1040,3710)(1040,3710)
\drawline(1040,3710)(1040,3710)
\drawline(1662,5210)(1662,6560)
\drawline(2712,6560)(3312,5435)
\drawline(1662,6560)(2712,6560)(2712,6560)
\dashline{60.000}(312,6035)(3012,6035)
\dashline{60.000}(86,5510)(3236,5510)
\dashline{60.000}(312,4910)(3012,4909)
\dashline{60.000}(687,4310)(2637,4310)
\dashline{60.000}(4587,6035)(7137,6035)
\dashline{60.000}(4287,5511)(7437,5511)
\dashline{60.000}(4513,4910)(7137,4910)
\dashline{60.000}(4887,4235)(6838,4235)
\drawline(6536,3706)(5224,3728)
\drawline(4886,6564)(4249,5417)
\drawline(4837,2898)(2753,2908)
\drawline(3309,5394)(2303,3687)(1596,3712)
\drawline(2753,2908)(2116,1761)
\drawline(4249,5417)(5224,3728)
\dottedline{90}(2937,2135)(4737,2135)
\dottedline{90}(2637,2660)(4962,2660)
\dottedline{90}(3162,5735)(2037,3710)
\dottedline{90}(4437,5735)(5562,3710)
\dottedline{90}(5037,5735)(5862,4235)
\dottedline{90}(5862,4760)(5487,5435)
\dottedline{90}(2487,5660)(1662,4310)
\dottedline{90}(2037,5435)(1662,4910)
\drawline(13651,2791)(13651,2791)
\drawline(13651,2791)(13651,2791)
\drawline(13662,2810)(11562,1610)
\drawline(9473,2791)(9473,2791)(11551,1591)
\drawline(13662,5210)(9462,2810)
\drawline(11562,6410)(11562,1610)
\drawline(9440,2810)(9440,2810)(9440,5210)
\drawline(11562,6410)(13662,5210)
\drawline(9462,5210)(13662,2810)
\drawline(11562,6410)(9462,5210)
\drawline(13640,2810)(13640,2810)(13640,5210)
\drawline(11562,6410)(11562,1610)
\drawline(11562,6410)(11562,1610)
\drawline(12612,2885)(13212,1760)
\drawline(10062,3710)(8862,3710)
\drawline(12686,6564)(12049,5417)
\drawline(10812,5960)(10812,5960)
\drawline(10737,5960)(11112,5435)
\drawline(13062,3710)(13662,3710)
\drawline(10512,2885)(10137,2360)
\put(12312,5585){\makebox(0,0)[lb]{\smash{{{\SetFigFont{5}{6.0}{\rmdefault}{\mddefault}{\updefault}$\Box$}}}}}
\put(11037,5660){\makebox(0,0)[lb]{\smash{{{\SetFigFont{5}{6.0}{\rmdefault}{\mddefault}{\updefault}$r_1(\Box)$}}}}}
\put(12612,6110){\makebox(0,0)[lb]{\smash{{{\SetFigFont{5}{6.0}{\rmdefault}{\mddefault}{\updefault}$r'(\Box)$}}}}}
\put(10437,2585){\makebox(0,0)[lb]{\smash{{{\SetFigFont{5}{6.0}{\rmdefault}{\mddefault}{\updefault}$r(\Box)$}}}}}
\put(8212,3860){\makebox(0,0)[lb]{\smash{{{\SetFigFont{5}{6.0}{\rmdefault}{\mddefault}{\updefault}$r_3r'(\Box)$}}}}}
\put(9837,3860){\makebox(0,0)[lb]{\smash{{{\SetFigFont{5}{6.0}{\rmdefault}{\mddefault}{\updefault}$r_3r_1(\Box)$}}}}}
\put(12987,3785){\makebox(0,0)[lb]{\smash{{{\SetFigFont{5}{6.0}{\rmdefault}{\mddefault}{\updefault}$r_2(\Box)$}}}}}
\put(12762,2735){\makebox(0,0)[lb]{\smash{{{\SetFigFont{5}{6.0}{\rmdefault}{\mddefault}{\updefault}$r_3r_2(\Box)$}}}}}
\put(13062,2135){\makebox(0,0)[lb]{\smash{{{\SetFigFont{5}{6.0}{\rmdefault}{\mddefault}{\updefault}$r_3r_2r'(\Box)$}}}}}
\end{picture}
}
\end{center}
\caption{ Locations of the images of $X_\Gbox$ and $\{Z_{u\Box}\}$.}
\end{figure}
\subsection{The final step}

{\it Proof of the Main Theorem:}
{\it Step 1:} By cutting and applying  Cor. 13.1,  we  first cancel the contributions to $\RR_G(Y)$,   from  the fixed point set $F$ with
 $\phi(F)\in w(\partial C\setminus C^\aff)$. The cancellation below
is implicit when applying Cor. 13.1.

By the fundamental property of cutting, one has
\begin{equation*}
\begin{split}
\RR(Y)&=\sum_{w\in W}w\frac{\RR(X_N)}{D}\\
&=\sum_{w\in W}w\frac{1}{D}\sum_{\Box\cap C\neq \emptyset}(-1)^{\cdim
\Box}\RR_T(X_\nbox) \\
&=\sum_{\Box\cap C\neq \emptyset}(-1)^{\cdim \Box}
\sum_{w\in W}w\frac{1}{D}\RR_T(X_\nbox)\\
&=\sum_{\Box\cap C\neq \emptyset}(-1)^{\cdim \Box}
\sum_{u\in W/W_\Box^0}u\frac{D_\Box}{D}u\RR(Z_\Box)\onlatk\\
\end{split}
\end{equation*}
where the last step uses Cor. 13.1. 
See the figure for the translations of the images of $\phi(X_\nbox)$, similar
actions are taken place for the varieties with images lying on the lower
dimensional $\Box$ which are not  illustrated.
\begin{figure}\label{134}
\begin{center}
\setlength{\unitlength}{0.00041667in}
\begingroup\makeatletter\ifx\SetFigFont\undefined%
\gdef\SetFigFont#1#2#3#4#5{%
  \reset@font\fontsize{#1}{#2pt}%
  \fontfamily{#3}\fontseries{#4}\fontshape{#5}%
  \selectfont}%
\fi\endgroup%
{\renewcommand{\dashlinestretch}{30}
\begin{picture}(6549,5983)(0,-10)
\drawline(3912,4981)(3462,5806)
\drawline(4212,2131)(5562,2131)
\drawline(4212,2131)(3837,1306)
\drawline(1812,1306)(3837,1306)
\drawline(3237,3856)(2337,3856)
\drawline(2337,3856)(1362,2131)
\drawline(5562,2131)(6537,3856)
\drawline(2787,2431)(2787,2431)
\drawline(2798,4850)(2798,4850)
\drawline(2798,4850)(2798,4850)
\drawline(2308,3861)(2308,3861)
\drawline(2308,3861)(2308,3861)
\drawline(2308,3861)(2308,3861)
\drawline(2787,2431)(2787,2431)
\drawline(698,1212)(698,1212)
\drawline(698,1212)(698,1212)
\drawline(709,1231)(709,1231)(709,3631)
\drawline(1309,2131)(1309,2131)
\drawline(1309,2131)(1309,2131)
\drawline(687,1231)(2787,31)
\drawline(2776,4850)(2776,4850)
\drawline(2776,4850)(2776,4850)
\drawline(3266,3861)(3266,3861)
\drawline(3266,3861)(3266,3861)
\drawline(3266,3861)(3266,3861)
\drawline(4876,1212)(4876,1212)
\drawline(4876,1212)(4876,1212)
\drawline(4265,2131)(4265,2131)
\drawline(4265,2131)(4265,2131)
\drawline(4265,2131)(4265,2131)
\drawline(4265,2131)(4265,2131)(4190,2131)
\drawline(4898,1212)(4898,1212)
\drawline(4898,1212)(4898,1212)
\drawline(4876,1212)(4876,1212)(2798,12)
\drawline(3796,1283)(3796,1283)
\drawline(3796,1283)(3796,1283)
\drawline(3796,1283)(3796,1283)
\drawline(6508,3861)(6508,3861)
\drawline(6508,3861)(6508,3861)
\drawline(6508,3861)(6508,3861)
\drawline(4898,1212)(4898,1212)
\drawline(5509,2131)(5509,2131)
\drawline(5509,2131)(5509,2131)
\drawline(5509,2131)(5509,2131)
\drawline(5509,2131)(5509,2131)(5584,2131)
\drawline(3237,3856)(4212,2131)
\drawline(65,2131)(65,2131)
\drawline(1362,2131)(1812,1306)
\drawline(12,2131)(1362,2131)
\drawline(1812,1306)(1137,256)
\drawline(4512,181)(3837,1306)(1887,1306)
\drawline(6537,3856)(5862,4981)(3912,4981)
\drawline(2337,3856)(1737,4981)
\drawline(3912,4906)(3237,3856)
\drawline(4887,3631)(687,1231)
\drawline(687,3631)(4887,1231)
\drawline(2787,4831)(2787,31)
\drawline(4909,1231)(4909,1231)(4909,3631)
\drawline(2787,4831)(2787,5956)
\drawline(3462,5806)(2112,5806)(1737,4981)
\drawline(4887,3631)(4887,4981)
\drawline(4887,3631)(6237,4381)
\dottedline{90}(3162,781)(5187,4381)
\drawline(5154.316,4261.703)(5187.000,4381.000)(5102.021,4291.119)
\dottedline{90}(2112,1081)(4212,4381)
\drawline(4172.885,4263.654)(4212.000,4381.000)(4122.265,4295.867)
\dottedline{90}(1287,1231)(3462,4906)
\drawline(3426.699,4787.451)(3462.000,4906.000)(3375.064,4818.010)
\dottedline{90}(987,1756)(3012,5506)
\drawline(2981.379,5386.157)(3012.000,5506.000)(2928.585,5414.666)
\dottedline{90}(987,2731)(5337,2731)
\drawline(5217.000,2701.000)(5337.000,2731.000)(5217.000,2761.000)
\dottedline{90}(1287,3631)(5862,3781)
\drawline(5743.048,3747.084)(5862.000,3781.000)(5741.081,3807.052)
\dottedline{90}(4662,1906)(2487,5506)
\drawline(2574.731,5418.804)(2487.000,5506.000)(2523.376,5387.777)
\dottedline{90}(4212,1156)(2037,4906)
\drawline(2123.157,4817.248)(2037.000,4906.000)(2071.255,4787.145)
\drawline(687,3631)(3687,5356)
\drawline(2787,4831)(6012,2956)
\drawline(2787,4831)(1887,5356)
\put(3762,3631){\makebox(0,0)[lb]{\smash{{{\SetFigFont{6}{7.2}{\rmdefault}{\mddefault}{\updefault}$\Box$}}}}}
\put(4137,4456){\makebox(0,0)[lb]{\smash{{{\SetFigFont{6}{7.2}{\rmdefault}{\mddefault}{\updefault}$r'(\Box)$}}}}}
\put(1962,856){\makebox(0,0)[lb]{\smash{{{\SetFigFont{6}{7.2}{\rmdefault}{\mddefault}{\updefault}$r(\Box)$}}}}}
\end{picture}
}
\end{center}
\caption{Comparing varieties with images differed by translation elements in
$W^\aff$}
\end{figure}
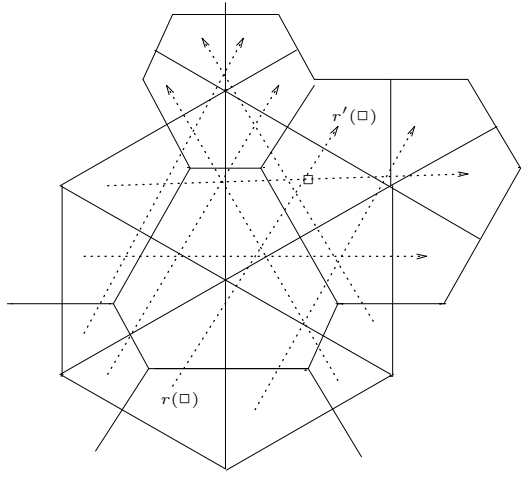

{\it Step 2: Localize to $V$. }
It should be clear that $u\RR(Z_\Box)=\RR(Z_{u\Box})$.  
While the intersection of the  image of $Z_\Box$ with $\ft$ is in
$W_\Box^\aff(\Box)$, the image of $uZ_\Box$ is in $$uW_\Box^\aff(\Box)=
\Ad_u(W_\Box^\aff)( u\Box)=W_{u\Box}^\aff(u\Box).$$ 
Each 
$\tau\in \explife$ is a generic element in $T$, hence 
the connected components of $\tau$-fixed points in $Z_{u\Box}$, denoted by
$\{V\}$, must have image under moment map
in $\ft$. Each $V$ is an orbifold, and is symplectic, by general theory on
fixed point sets and the fact that $Z_{u\Box}$ is a symplectic  orbifold. 
 We may apply Prop.~\ref{cutup}, replacing  $M,G, V$ there   by $Z_{u\Box},
K, V$ where  $K=(LG)_{u\Box}$. The collection $\{\Delta\}$ there  is now replaced by  $\{uw' (\Box\cap C)|w'\in W_\Box^\aff\}$ which is denoted by $\{\Delta'\}$ in the following. Denote by  $Q_1$ a  face of $ u(\Box\cap \partial C)$, and
$$Q=\Ad_uw(Q_1),\quad Q'=\Ad_uw'(Q_1)$$ 
where $w'\in W_\Box^\aff$ and $w$ is the corresponding element in $W_\Box^0$. 
We obtain
\begin{equation*}
\begin{split}
&u\big(\RR(Z_\Box)\big)(\tau)\\
&=\RR(Z_{u\Box})(\tau)\\
&=\sum_{\Delta'}\sum_{Q'\subset \Delta'}\frac{1 }{|W_{Q'}||I_{V_{Q'} } |}\sum_{t\in \tau I_{V_{Q'} }} \int_{V_{ Q'}} \frac{\tTd(V_{Q'})\tCh(L_{V_{Q'}}\oplus
\Lambda^{\max}\no(V_{Q'},V_{\Delta'} ) )}{
\det_{\no(V_{Q'},X_{Q'}^ K)}(1-t^{-1}e^{-\Omega})}.
\end{split}\end{equation*}
The varieties $\{Z_{u`\Box}\}$ have images as in the following figure:
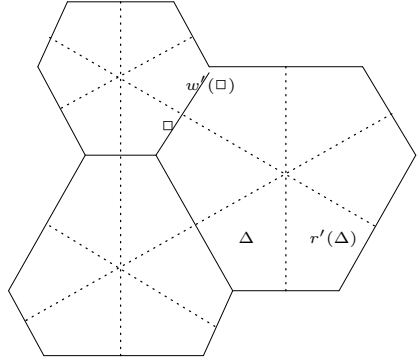
\begin{figure}\label{135}
\begin{center}
\setlength{\unitlength}{0.00041667in}
\begingroup\makeatletter\ifx\SetFigFont\undefined%
\gdef\SetFigFont#1#2#3#4#5{%
  \reset@font\fontsize{#1}{#2pt}%
  \fontfamily{#3}\fontseries{#4}\fontshape{#5}%
  \selectfont}%
\fi\endgroup%
{\renewcommand{\dashlinestretch}{30}
\begin{picture}(5199,4539)(0,-10)
\drawline(2862,837)(4212,837)
\drawline(2862,837)(2487,12)
\drawline(462,12)(2487,12)
\drawline(1887,2562)(987,2562)
\drawline(987,2562)(12,837)
\drawline(1887,2562)(2862,837)
\drawline(2562,3612)(1887,2562)
\drawline(5187,2562)(4512,3687)(2562,3687)
\drawline(987,2562)(387,3687)
\drawline(762,4512)(2112,4512)
\drawline(4212,837)(5187,2562)
\drawline(237,462)(237,462)
\drawline(2562,3687)(2112,4512)
\drawline(12,837)(462,12)
\drawline(762,4512)(387,3687)
\dottedline{90}(237,462)(4887,3087)
\dottedline{90}(1437,4512)(1437,12)
\dottedline{90}(537,4062)(4662,1662)
\dottedline{90}(687,3162)(2337,4062)
\dottedline{90}(3537,3687)(3537,837)
\dottedline{90}(537,1662)(2637,387)
\put(2262,3387){\makebox(0,0)[lb]{\smash{{{\SetFigFont{6}{7.2}{\rmdefault}{\mddefault}{\updefault}$w'(\Box)$}}}}}
\put(1962,2862){\makebox(0,0)[lb]{\smash{{{\SetFigFont{6}{7.2}{\rmdefault}{\mddefault}{\updefault}$\Box$}}}}}
\put(3837,1437){\makebox(0,0)[lb]{\smash{{{\SetFigFont{6}{7.2}{\rmdefault}{\mddefault}{\updefault}$r'(\Delta)$}}}}}
\put(2937,1437){\makebox(0,0)[lb]{\smash{{{\SetFigFont{6}{7.2}{\rmdefault}{\mddefault}{\updefault}$\Delta$}}}}}
\end{picture}
}
\end{center}
\caption{After translation and the twin pair comparison}
\end{figure}

{\it Step 3: Transporting $V_{Q'}$.}  
The transporting here are in the directions opposite to the arrows in the figure~\ref{134}. 
One does not simply get back the  result in the beginning, since the
cancellation has already taken place. This is an important realization.

Two observations can be made here: Each integral  $\int_{V_{Q'} }$ in the above  is finite, when evaluated
at $\tau$, because $V_{\Ad_uw' (Q)}$ is a connected component of $\tau$-fixed points in
$Z_Q^ K$, so the action by $\tau$ on the normal bundle is non-trivial, thus
the denominator is well defined at $\tau$. Second
observation is that one can apply Prop. \ref{rule} 2) to this situation.
To make this more explicit,
let $\Ad_uw=s\Ad_uw'$ where $s$ is a translation defined by an element in the
long root lattice, according to  Lemma 13.3. 
For each $\Ad_uw' (\Box\cap C)$, there is the corresponding 
$\Ad_uw (\Box\cap C)$, obtained by shifting $- s$.
Likewise $\Ad_uw ( Q_1)$ is obtained from $\Ad_uw' (Q_1)$ by shifting $- s$.
The corresponding varieties $V_Q\subset X_Q$ have  a similar 
relation $s^{-1}(V_Q)=V_{-s(Q)}\subset s^{-1}(X_Q)=X_{-s(Q)}$, and $s^{-1}(V_Q)$ is a connected
component of $\tau$-fixed points, since the translation commutes with
$T$-action.

The denominator is given by 
$$\det_{\no(V_{Q'},X_{Q'} ^K)} (1-t^{-1}e^{-\Omega})=
uw'D_{\Box}\det_{\no(V_Q,X_Q)} (1-t^{-1 }e^{-\Omega}),
$$ thus we obtain the following relation as a consequence of Prop. \ref{rule}:
\begin{equation*}
\begin{split}
& e^{(k+h^\vv)s}u\frac{D_\Box}{D}\int_{V_{Q'}} \frac{\tTd(V_{Q'})\tCh(L_{V_{
Q'}}\oplus
\Lambda^{\max}\no(V_{Q'},V_{\Delta'} ) )}{
\det_{\no(V_{Q'},X_{Q'}^ K)}(1-t^{-1}e^{-\Omega})}\\
&=\int_{V_{Q}} \frac{\tTd(V_{Q})\tCh(L_{V_{Q}}\oplus
\Lambda^{\max}\no(V_{Q},V_{\Delta} ) )}{
\det_{\no(V_{Q},X_{Q}^ G)}(1-t^{-1}e^{-\Omega})}
\end{split}\end{equation*}
where $Q'=uw'(Q_1), Q=uw(Q_1)$ with $Q_1$ a face of $\Box\cap C$.

Both sides have no poles at $\tau$ since $t\in \tau I_{V_{Q'}}$ acts 
on the normal bundles with no 0 eigenvalue. 
 Evaluate them  at  $\tau\in\exp(\latk)$, i.e. $e^{(k+h^\vv)s}=1$,
for $Q=s Q'$, $  \Delta=s\Delta'$,  one has the following 
\begin{equation}\label{final6}
\begin{split}
&\RR(Y)(\tau)\\
&=\sum_{\Box\cap C\neq \emptyset}(-1)^{\cdim \Box}
\sum_{u\in W/W_\Box^0}
 u \sum_{w\in \Ad_uW_\Box}
\frac{1}{|W_Q||I_{V_{Q}}|}\sum_{t\in \tau I_{V_{Q}} } \calI_{V_{Q}},\\
&\calI_{V_{Q}}=
\int_{V_{Q}} \frac{\tTd(V_{Q})\tCh(L_{V_{Q}}\oplus
\Lambda^{\max}\no(V_{Q},V_{\Delta} ) )}{
\det_{\no(V_{Q},X_{Q}^ G)}(1-t^{-1}e^{-\Omega})}.
\end{split}
\end{equation}

{\it Step 4: Eliminate the extra things.} 
Next, we will write the above in a concise form.
To see further  cancellation, it is  easiest to write the integrals above
in terms of the $T$-fixed points contribution:
\begin{equation*}\label{final5}
\begin{split}
&\frac{1}{|I_{V_Q}|}\sum_{t\in \tau I_{V_Q}} \int_{V_{Q}} \frac{\tTd(V_{Q})\tCh(L_{V_{Q}}\oplus
\Lambda^{\max}\no(V_{Q},V_{\Delta} ) )}{
\det_{\no(V_{Q},X_{Q}^ G)}(1-t^{-1}e^{-\Omega})}\\
&=\sum_{F\subset V_{Q}}\frac{1}{|I_F|}
\sum_{t\in \tau I_F}\int_{F}
 \frac{\tTd( F)\tCh(L_F \oplus \Lambda^{\max}\no(V_{Q},V_{\Delta}|_F ) )}
{ \det_{\no(F,X_{Q}^ G)}(1-t^{-1}e^{-\Omega})}  .
\end{split}\end{equation*}
Substitute the above into  Eq. (\ref{final6}), then we will treat those terms
 according to
the type of $F$:

a). $\phi(F)$ is on a cut defined by $\Box_1$ with $\dim\Box_1<l=\dim \ft$.  

b). $\phi(F)$ is on some $wQ, Q\subset \partial  C\setminus C^\aff,  $; but not on $\Box_1$ with $\dim\Box_1<l=\dim \ft$.

c). $\phi(F)$ is on some $wQ $ with $ Q\subset  C^\aff$; but not
on $\Box_1$ with $\dim\Box_1<l=\dim \ft$.

d). $\phi(F)$ is neither on any  $wQ , \forall w\in W, \forall Q\subset
\partial C$
nor on $\Box_1$ with $\dim\Box_1<l=\dim \ft$.

Obviously, the above covers all possibilities for $\phi(F)$. We claim
the contributions of $F$ of the first two kinds amount to  0 in Eq.~(\ref{final5}).

In the first case, suppose $u(\Box_1\cap C), u(Q_1)$ are  the smallest among 
$ u(\Box\cap C), u(Q) $   containing
$\phi(F)$. Fixing  a $u(Q)\supseteq u(Q_1)$,  and $u(\Box)\supseteq 
u(\Box_1)$, what is the integrand in the above integral?

For convenience of notations, assume below that $u=I$.

 Let us take a moment to  discuss the weights on  $\no(F, V_{Q\cap \Box})$ and the bundle 
$\nor(V_{Q\cap \Box }, V_{\Box})|_F$.
Suppose that $\bAbox=\{\ta\}$ is a basis of  weights in $\Box$ normal to 
$\Box_1$. 
By the construction of the cut space we know in $\mu^{-1}(\ft)$, there is the set $\tF$ whose quotient under the subgroups
$T_Q \times  T_{\Box_1}$ is $F$, here $\Lie T_Q$ is  normal to $Q$, $\Lie
T_{\Box_1}$ is normal to $\Box_1$. 
Now assume that $T_P$ stabilizes a generic point on $\tF$, then
$F$ is fixed by $T$ iff 
$$\Lie T_P\oplus \Lie T_{\Box_1} \oplus \Lie T_Q=\ft.$$
The proof of the above follows from  that of  Lemma 3.1. 
Let $\{\tlambda\}, \{\tgamma\}$ be weights of $\ft_Q, \ft_P$  respectively,
as in Prop. 4.3, (replacing $\ft_z$ there by $\ft_Q$),
and let $\tbeta$ be the roots of the group $(LG)_Q$ which is the stabilizer in
$LG$ of
$Q$. 

 Each $t\in T$ comes from a product of a triple $$t_P\cdot t_{\Box_1}\cdot t_Q\in T_P\cdot T_{\Box_1}\cdot T_Q,$$ the finite ambiguity of the choice of the
triple causes  $F$ to
be an orbifold singularity. 
So the orbifold weights on the normal bundle of $F$ in $X_{Q\cap\Box}$  can be described as follows:
\begin{equation}
\begin{split}
\lambda(t)=\tlambda(t_Q) ;\quad \gamma(t)=\tgamma(t_P); \\
\quad a(t)=\ta(t_{\Box_1});\quad \beta(t)=\tbeta(t_P)
\end{split}
\end{equation}
the verification of the above is the same as Prop. 4.3. 

The normal bundle $\nor(V_{Q\cap \Box},V_{  C\cap \Box})|_F$ has weights
given by $\lambda(t)=\tlambda(t_Q), \tlambda\notin \Lambda_Q$ where
$\Lambda_Q$ is the subset of $\{\tlambda\}$ not parallel to $Q$.

One realizes in the above that for different $\Box\supseteq \Box_1$, the 
only difference in the integrand is the term $\prod_{a\in \bAbox}(1-e^{-\lp<a,t+\fp dA>})
$ appearing in the denominator, here $B$ is the form representing the Chern
class of the principle bundle corresponding to the orbit of $T_{\Box_1}$.

Let $D_a(s), D_b(s)$ be the same as in Section 11, and $D'_0$ be defined
the same as $D_0$ except only those $\gamma$ tangent to  $\Box_1$ will be involved.
 We have the following representation of the integrand:
\begin{equation}
\begin{split}
&(-1)^{\cdim\Box}\sum_{\Box\supseteq \Box_1}\int_{F}
 \frac{\tTd( F)\tCh(L_F \oplus \Lambda^{\max}
\no(V_{Q\cap \Box},V_{C\cap \Box} )|_F  )}
{ \det_{\no(F,X_{Q\cap\Box}^ G)}(1-s^{-1}e^{-\Omega}) }  \\
&=\sum_{\Box\supseteq \Box_1}
\int_F \frac{\tTd( F)\tCh(L_F)}{e^{\sum_{\tlambda\notin \Lambda_Q} 2\pi i<\lambda, t-t_P-\fp dA> }
D_a(s)D_b(s)D'_0(s)}\\
&\quad
\times{\sum_{\bAbox}(-1)^{\cdim \Box}
\frac{1}{\prod _{a\in \bAbox}(1-e^{-<a,t_P+\fp dA>})} }\\
&=0
\end{split}
\end{equation}
since 
$$ \sum_{\bAbox}(-1)^{\cdim \Box}\frac{1}{\prod _{a\in
\bAbox}(1-e^{-<a,t_P+\fp dA>})} =0.$$

As for $F$ of the second type, $\phi(F)$ is in the interior of $\cup_ww(C)$,
but on a wall of Weyl chamber $w(Q)$. 
The  sum of the integrals  over various $Q'$ with 
$Q\subset Q'$ was already  shown to be 0  in    the proof of the Prop. \ref{cutup}.
 The argument there for those fixed points
appearing on the boundary of the Weyl chamber shows the same  cancellation here.

Therefore, only $F$ in c) or d) survives. We do recognize these two types,
$F$ of type c) are those produced by the intersection of $\tau$ fixed point
set $V$,
  with the compactification locus  corresponding to $ w(C^\aff)$, i.e.
the affine walls. Those of type d) are  the fixed point sets in $X$ with
images in $W(C^\inte)$.  
Thus we have the following form for $\RR(Y)$ when evaluated at $\tau$:
\begin{equation}\label{final}
\begin{split}
&\RR(Y)=\sum_{\{F|\phi(F)\in W(C^\inte)\}} \FC(F) \\
&\quad\quad  +\sum_{\{F|\phi(F)\in W(C^\aff)\}} \frac{1}{|W_Q||I_F|}\sum_{t\in \tau I_F}\int_{F}\frac{\tTd( F)\tCh\left(L_F
\oplus
\Lambda^{\max}\no(Y_Q,Y )|_F \right)}{ \det_{\no(F,Y_Q)}(1-s^{-1}e^{-\Omega})}.
\end{split}
 \end{equation}
 The first sum is over the true fixed points on $X\cap \mu^{-1}(W(C))$,
the absence of the isotropy  group $I_F$ is due to the smoothness of $Y$ at 
the interior fixed points, since no cutting passes such  $F$.
The second sum is over those fixed points lying on the intersection of the  compactification locus
$Y_Q$ and a $\tau$-fixed point set component.
 The only explanation needed here
is that 
$$\tCh(\no(Y_Q,Y _\Delta)|_F)=\tCh(\no(V_{Q\cap \Box_1},V_{  C\cap \Box} )|_F)$$
which is obvious in terms of these  weights  $\lambda$ with $\tlambda\in Q$.

Thus we obtain the second expression for $\RR(Y)(\tau)$ in the main theorem.
The first one can be obtained easily now by localize the integrals over
$V_\Delta$, $V_Q$ to their fixed points.
QED.

We emphasize that there are plenty of $T$-fixed  point sets on the
compactified locus, they do not contribute to the formula in Eq. (\ref{final})  
unless they are on the compactified $\tau$-fixed points.

\section{References}
\noindent [AB1] M .\,F.\,Atiyah and R.\,Bott.{\it  A Lefschetz fixed
point formula for elliptic complexes. I} {\em Ann. Math.}  86
(1967), 374--407.

\noindent [AB2] M.\,F.\,Atiyah and R.\,Bott. {\it A Lefschetz fixed
point formula for elliptic complexes. II} Applications {\em Ann. Math.}
{ 87} (1968), 451--491.

\noindent [AB] M. Atiyah \& R. Bott, {\it The Yang-Mills equations over
Riemann
surfaces}, Phil. Trans. R. Soc. London, A308 (1982) 523.

\noindent[AS] M.\,F.\,Atiyah and G.\,B.\,Segal. {\it The index
of elliptic operators II}. {\em Ann. Math.} { 87} (1968), 531--545.

\noindent [AS1] M.\,F.\,Atiyah and I.\,M.\,Singer. {\it The
index of elliptic operators I.} {\em Ann. Math.} {87} (1968),
484--530.

\noindent [AS3]  M.\,F.\,Atiyah and I.\,M.\,Singer. {\it The
index of elliptic operators III.} {\em Ann. Math.} { 87} (1968),
546--604.

\noindent [AS4] M.\,F.\,Atiyah and I.\,M.\,Singer. {\it The index of
elliptic operators IV.}  {\em Ann. Math.} { 93} (1971), 119--138.

\noindent [BGV] N. Berline, E. Getzler \& M. Vergne, {\it Heat kernels and
Dirac operators}, preliminary version (1990).

\noindent [Be] A. Beauville, {\it Conformal blocks, fusion rules and the
Verlinde formula}. Preprint, (1994)

\noindent [Bi] J.-M.\,Bismut. {\it The infinitesimal Lefschetz
formulas: a heat equation proof}. {\em J. Funct. Anal.} { 62} (1985),
435--457.

\noindent [B] R.\,Bott. {\it  Homogeneous vector bundles.} {\em Ann. Math.} {
66} (1957), 203--248.

\noindent [C]  S.\,Chang. {\it A fixed point formula on orbifolds}.
MIT preprint 1995.

\noindent [C1] S.\,Chang. {\it Geometric  loop group actions and the highest weight modules}.  MIT preprint 1996.

\noindent [F] G. Faltings, {\it A proof for the Verlinde formula}, J.
Algebraic
 Geom. 3, 347-374 (1994).

\noindent [K] V. Kac.{\it Infinite dimensional Lie algebras} 3rd ed. Cambridge
University Press, 1990.

\noindent [KP] V. Kac \& D. Peterson, {\it Infinite-dimensional Lie algebras,
theta functions and modular forms,} Advances in Math. 53 (1984), 125-264

\noindent [KW] V. Kac \& M. Wakimoto, {\it Modular and conformal invariance
constraints in representation theory of affine algebras}, Advances in Math. 70
(1988), 156-234.

\noindent [KS] M. Kreck \& S. Stolz, {\it Nonconnected moduli spaces of positive
curvature metrics.} JAMS, Vol. 6, No. 4, 825-850 (1993).

\noindent [GS] V. Guillemin \& S. Sternberg, {\it Geometric quantization and
multiplicity of group representations}, Invent. Math. 67 (1982) 515-538.

\noindent [Ma] I. Macdonald. {\it  Affine Lie algebras and modular forms}, 
258-276, Seminaire
Bourbaki Vol. 1980/81 Exposes 561-578, Lect. Notes in Math. 901.

\noindent [M] E.\,Meinrenken. {\it Symplectic Surgery and the Spin-C Dirac
operator}, dg-ga/9504002. 

\noindent [MS] G. Moore \& N. Seiberg, {\it Classical and quantum conformal
field theory} Comm. Math. Phys. 123 (1989) 177-254.

\noindent [PS] A. Pressley \&  G. Segal, {\it Loop group}, Oxford Sciece
Publications, 1986.

\noindent [S] G. Segal, {\it Two dimensional conformal field theories and
modular functors}, IXth International Congress on Math. Phys. 1988, edited by B.
Simon
 {\it et al} 22-37.

\noindent [T] C. Teleman, {\it Lie algebra cohomology and the fusion rules.}
Comm. Math. Phys. 173 (1995), no. 2, 265--311.

\noindent [TZ] Y. Tian \& W. Zhang, {\it
Symplectic reduction and analytic localization }, Courant Institute Preprint (1996)

\noindent [TUY] A. Tsuchiya, K. Ueno \& Y. Yamada, {\it Conformal field theory
on universal family of stable curves with gauge symmetris}, Advanced Studies in
Pure Math. 19 (1982)

\noindent [Ve]  M. Vergne, {\it  Multiplicity formula for geometric
quantization
 Part I, II},
Duke Math. J. 82 (1996), 143-179; Duke Math. J. 82 (1996), 181-194.

\noindent [V] E. Verlinde, {\it Fusion rules and modular transformations in 2D
conformal field theory,} Nuclear Phys. B300 (1989) 360-376.

\vskip .1in
\noindent Department of Mathematics 2-271\\
 M. I. T.\\
Cambridge, MA 02139

\noindent schang@@math.mit.edu

\end{document}